\documentclass[twoside,
 12pt, a4paper]{article}
\usepackage{amsmath,amsthm,amssymb,amscd}
\usepackage{bm}
\usepackage{latexsym}
\usepackage{indentfirst}
\usepackage{graphics}
\usepackage{amsmath}
\usepackage{amssymb}
\usepackage{bm}
\usepackage{latexsym}
\usepackage[mathscr]{eucal}
\makeatletter \@addtoreset{equation}{section}

\makeatother
\newcommand{\wth}{\widetilde{h}}
\newcommand{\wt}{\widetilde}
\newcommand{\mP}{\mathsf{P}}
\newcommand{\mQ}{\mathsf{Q}}
\newcommand{\mO}{\mathsf{O}}
\newcommand{\bi}{\mathbf{i}}
\newcommand{\bj}{\mathbf{j}}
\newcommand{\bk}{\mathbf{k}}
\newcommand{\ba}{{\mbox{\boldmath{$\alpha$}}}}
\newcommand{\sba}{{\mbox{\scriptsize\boldmath{$\alpha$}}}}
\newcommand{\bb}{{\mbox{\boldmath{$\beta$}}}}
\newcommand{\sbb}{{\mbox{\scriptsize\boldmath{$\beta$}}}}
\newcommand{\bzero}{\mathbf{0}}
\newcommand{\mi}{\mathrm{i}}
\newcommand{\me}{\mathrm{e}}
\begin{document}
\title{Perturbation of self-similar sets and
 some
 regular configurations
  and comparison of fractals
 }
\author{Junyang Yu
\thanks{Partially supported by NSFC grant 10271077.}
 }                     
%
%
\date{}
%
\maketitle
\begin{abstract}
We consider several distances between two sets of points, which
are modifications of the Hausdorff metric, and apply them to
describe some fractals such as $\delta$-quasi-self-similar sets,
and some other geometric notions in Euclidean space, such as
tilings with quasi-prototiles and patterns with quasi-motifs. For
the $\delta$-quasi-self-similar sets satisfying the open set
condition we obtain the same result as a classical theorem due to
P.~A.~P.~Moran. In this paper we try to gaze on fractals in an
aspect of their ``form'' and suggest a few of related questions.
Finally, we attempt to inquire an issue
--- what nature and behavior do non-crystalline solids that
approximate to crystals show?
\end{abstract}
\section{Introduction}
\label{intro} \nopagebreak
 Generally fractals are considered to
possess three important features: {\it form}, {\it chance}\/ and
{\it dimension}, just as indicated in the titles
of~\cite{Mandelbrot-1975} and~\cite{Mandelbrot-1977}. The
dimension has been a very important and fundamental subject in
researches of fractals (see e.g. \cite{Falconer-2003}
and~\cite{Mandelbrot-1982}). Random fractals as more natural
description of things in nature are also extensively investigated
by experts in many subjects (see e.g.
\cite{Barnsley-Hutchinson-Stenflo-2005}, \cite{Falconer-1986},
\cite{Falconer-2003}, \cite{Graf-Mauldin-Williams-1988},
\cite{Hutchinson-Ruschendorf-1998}, \cite{Mandelbrot-1982} and
\cite{Mauldin-Williams-1986}). The research of fractals is also
closely connected to geometric measure theory (refer to
\cite{Falconer-1985}, \cite{Federer-1969}, \cite{Mattila-1995} and
\cite{Rogers-1998}, etc) and other scientific subjects (refer to
\cite{Feder-1988}, \cite{Mandelbrot-1977}, \cite{Mandelbrot-1982},
\cite{Mandelbrot-1984} and \cite{Mandelbrot-1989}, etc).

Usually fractals are also considered to possess recursive or
recurrent structure. Self-similarity is one of simple and
important natures of fractals, where the self-similarity may often
be approximate or statistical.

A mathematical (strict) self-similar set as an extension of the
classical Cantor set has been investigated deeply and extensively
(see \cite{Falconer-2003}, \cite{Hutchinson-1981} and
\cite{Moran-1946}, etc). In this paper we will describe an
approximate self-similar set in a quantitative respect,
considering it as a perturbation of a strict self-similar set.

However, for the fine structure of a (mathematical) self-similar
set, any small perturbation will probably destroy its fractal
details. For example, assume that $F$ is a self-similar set
(fractal) in $\mathbb R^n$ ($n$-dimensional Euclidean space)
satisfying
$$
h(F,E)\leqslant\varepsilon,
$$
where $\varepsilon>0$, $E$ is a subset of $\mathbb R^n$ and  $h$
is the Hausdorff metric. Then no matter how small $\varepsilon$
is, $E$ can be chosen as a usual Euclidean figure.

The research on tilings, patterns and packings has a long history,
which was once advanced by Hilbert's 18th problem and developments
of other subjects, especially crystallography. For the context and
introduction, let us refer to \cite{Brass-Moser-Pach-2005},
\cite{Conway-Sloane-1999}, \cite{Goodman-O'Rourke-2004},
\cite{Gruber-2007}, \cite{Gruber-Lekkerkerker1987},
\cite{Grunbaum-Shephard-1980}, \cite{Grunbaum-Shephard-1987},
\cite{Hilbert-1902} and \cite{Senechal-1995}, etc.

In this paper, first we consider several modified Hausdorff
metrics, which are called {\it shape differences}\/ and are proved
to be complete metrics (in appropriate spaces)(see Section~2), and
then we perturb self-similar sets scale by scale using them. We
also get the Hausdorff dimension of the perturbed self-similar set
by a classical method (\cite{Hutchinson-1981} and
\cite{Moran-1946}) (Section~3). Furthermore, we approach the form
of fractals by comparison and suggest the notion of {\it splines
of fractals}, the {\it fractal index}\/ and the {\it similarity
index}\/ to indicate inside structures and complexity of fractals
(Section~4). In the last part (Section~5), we try to modify some
classical concepts in tiling, pattern, packing and crystallography
to lead valuable investigations from other people. As an example
we extend a basic result about engulfing in the research of
patterns (see \cite[5.1.1]{Grunbaum-Shephard-1987}). In the paper
we suggest a few of related questions for consideration.

In fact, the modified Hausdorff metric has been studied by experts
in computational geometry for quite a long time (see e.g.
\cite{Alt-Brab-Godau-Knauer-Wenk-2003}, \cite{Brab-2002},
\cite{Chew-Goodrich-Huttenlocher-Kedem-Kleinberg-Kravets-1997},
\cite{Huttenlocher-Kedem-Sharir-1993}, etc). But according to my
limited knowledge the experts mainly make investigations in the
aspects of algorithms and their time and so on.
 A modified Hausdorff metric (the Hausdorff-Chabauty distance)
  has been also defined in~\cite{Tan-1990} and some
wonderful observations about the Mandelbrot set and Julia sets
were proved there.

The introduction to researches on fractals (in many subjects of
science), tilings, patterns, packings and crystals can not be
included here for their tremendous amount. For example, only in
the investigation of the Mandelbrot and Julia sets the researches
have been so plentiful that perhaps one ordinary book can not
include them all. Here we just mention one result (about
``topological form'' of a fractal) that the Mandelbrot set of
$f_c(z)=z^2+c$ is connected, which was observed by
B.~B.~Mandelbrot (see \cite{Mandelbrot-1982}
and~\cite{Mandelbrot-1983}) and proved by A.~Douady and J.~Hubbard
(see~\cite{Douady-Hubbard-1982}).

In this paper by the term ``fractal'' we mean not only an
irregular object but also sometimes a regular one.

In the preface of \cite{Senechal-1995} M.~Senechal said:
``$\cdots$ Some of these tools were new to me, and although I have
enjoyed the adventure of learning how to use them, I am also aware
that I may have made errors or am ignorant of the relevant
literature. I will be grateful for any criticisms, comments, and
suggestions: the adventure continues.'' This is also my feeling
while writing this paper. It is my hope that this paper might play
a role of ``casting a brick to attract a gem''.

\section{Shape differences}
\label{sec:3} \nopagebreak
 Let $\boldsymbol{X}$ be a complete
metric space with a metric~$d$. Define
$$d(A,x)=d(x,A):=\inf\{d(x,a):\: a\in A\}$$
for $\boldsymbol{X}$ and a subset~$A$ of~$\boldsymbol{X}$. For
$\delta\geqslant0$ we denote
\begin{equation*}
\mathscr{P}_\delta(A)=\mathscr{P}(A,\delta):=\{x:\:
d(x,A)\leqslant\delta\},
\end{equation*}
 which is  called a {\it $\delta$-parallel
body}\/ of~$A$, and
\begin{equation*}
\mathscr{N}_\delta(A)=\mathscr{N}(A,\delta):=\{x:\:
d(x,A)<\delta\}
\end{equation*}
 is called a {\it $\delta$-neighborhood}\/ of~$A$ or an {\it open $\delta$-parallel
body}\/ of~$A$.

Suppose $A$ and~$B$ are two nonempty subsets of~$\boldsymbol{X}$.
Define
\begin{equation}
h(A,B):=\max\left\{\sup_{a\in A}d(a,B),\ \sup_{b\in
B}d(A,b)\right\},\label{tag2.3}
\end{equation}
which is called the {\it Hausdorff distance} ({\it metric})
between $A$ and~$B$. It also follows that
\begin{equation}
h(A,B)=\inf\left\{\delta:\: A\subseteq\mathscr{P}_\delta(B),\
B\subseteq\mathscr{P}_\delta(A) \right\}.\label{tag2.4}
\end{equation}
Let
\begin{equation*}
\mathcal C(\boldsymbol{X}):=\left\{C:\: C \textrm{ is a nonempty
closed bounded subset of }~\boldsymbol{X}
\right\}.
\end{equation*}
Then $\left(\mathcal C(\boldsymbol{X}),h\right)$ is a complete
metric space (see \cite[2.10.21]{Federer-1969} and
\cite[\S~2.6]{Rogers-1998}). Here and hereafter $\mathbb R$
indicates the real numbers and $\mathbb P$ indicates the set of
positive integers.

\paragraph{{\bf 2.1. Shape differences.}}

We consider $n$-dimensional Euclidean space~$\mathbb R^n$ with
usual Euclidean distant~$d$.
\paragraph{{\bf 2.1.1.}}
 We say that two nonempty subsets
$A$ and~$B$ of~$\mathbb R^n$ are {\it isometrically equivalent} if
there exists an isometry $\varphi:\: \mathbb R^n\rightarrow\mathbb
R^n$ such that $B=\varphi(A)$, denoted by $A\sim B$. Obviously
this relation is an equivalence relation. The equivalence class
of~$A$ is denoted~$\wt{A}$.

\bigskip
 From~(\ref{tag2.3}) or~(\ref{tag2.4}) it follows easily that
\paragraph{{\bf 2.1.2. Lemma.}} 
 {\it
 Let\/ $\varphi$ be an isometry in\/~$\mathbb R^n$ and let\/ $A$
and\/~$B$ be nonempty subsets of\/~$\mathbb R^n$. Then
$h(\varphi(A),\varphi(B))=h(A,B)$.\hfill $\Box$
 }

\paragraph{{\bf 2.1.3. Definition.}} 
 For nonempty subsets $A$
and~$B$ of~\/$\mathbb R^n$ we define
$$
\wth(A,B):=\inf\left\{h(\varphi(A),\psi(B):\: \textrm{$\varphi$
and $\psi$ are isometries in~$\mathbb R^n$}\right\}.
$$
By Lemma~2.1.2 we have the following expressions:
\begin{align*}
\wth(A,B):&=\inf\left\{h(A_1,B_1):\: A_1\sim A,\ B_1\sim B\right\}\\
 &=\inf\left\{h(A_0,B_1):\: B_1\sim B\right\}\quad\textrm{(where $A_0\sim A$)}\\
 &=\inf\left\{h(A_1,B_0):\: A_1\sim A\right\}\quad\textrm{(where $B_0\sim B$)}.
\end{align*}

If $A_1\sim A_2$ and $B_1\sim B_2$, then
$\wth(A_1,B_1)=\wth(A_2,B_2)$. So we may define
$$
\wth(\wt{A},\wt{B}):=\wth(A,B).
$$
We call $\wth(\wt{A},\wt{B})=\wth(A,B)$ the ({\it absolute}) {\it
shape difference} between $\wt A$ and~$\wt B$ or between $A$
and~$B$.

\paragraph{{\bf 2.1.4.}}
 Denote
$$
\wt{\mathcal{C}}(\mathbb{R}^n):=\left\{\wt{C}:\: C \textrm{ is a
nonempty compact subset of }~\mathbb{R}^n\right\}.
$$
Then $\wt{\mathcal{C}}(\mathbb{R}^n)=\left\{\wt{C}:\: C
\in\mathcal{C}(\mathbb{R}^n)\right\}.$

\paragraph{{\bf 2.2. Theorem.}} 
 {\it
$\bigl(\wt{\mathcal{C}}(\mathbb{R}^n),\wth\bigr)$ is a complete
metric space.
 }
\par\noindent{\it Proof.}   At first we show that $\wth$ is a distance
function on~$\wt{\mathcal{C}}(\mathbb{R}^n)$. Let $A$, $B$ and
$C\in\wt{\mathcal{C}}(\mathbb{R}^n)$.

(i) It is trivial that $\wth(\wt A,\wt B)=\wth(\wt B,\wt A)$.

(ii) The triangle inequality $\wth(\wt A,\wt B)\leqslant\wth(\wt
A,\wt C)+\wth(\wt C,\wt B)$ follows from that for any
$\varepsilon>0$ we have
\begin{align*}
\wth(\wt A,\wt B)&\leqslant h(A,B_1)\leqslant
h(A,C_0)+h(C_0,B_1)\\
&<\left(\wth(\wt A,\wt C)+\frac\varepsilon2\right)+\left(\wth(\wt
C,\wt B)+\frac\varepsilon2\right)\\
&=\wth(\wt A,\wt C)+\wth(\wt C,\wt B)+\varepsilon,
\end{align*}
where $C_0\sim C$ is chosen to satisfy
$$
h(A,C_0)<\wth(\wt A,\wt C)+\frac\varepsilon2
$$
and then $B_1\sim B$ is chosen such that
$$
h(C_0,B_1)<\wth(\wt C,\wt B)+\frac\varepsilon2.
$$

(iii) Obviously $\wth(\wt A,\wt B)\geqslant0$ and if $\wt A=\wt B$
then $\wth(\wt A,\wt B)=0$. Below we will prove that if $\wth(\wt
A,\wt B)=0$ then $\wt A=\wt B$.


\paragraph{{\bf 2.2.1. Definition.}} 
Denote the set of $l\times m$ matrices with entries in~$\mathbb R$
by $\mathbb R^{l\times m}$ ($l$, $m\in\mathbb P$). If
$\mP=\left(p_{_{\scriptstyle ij}}\right)_{l\times m}\in\mathbb
R^{l\times m}$, then define
$$
\|\mP\|:=\sum_{i=1}^l\sum_{j=1}^m|p_{_{\scriptstyle ij}}|.
$$

\paragraph{{\bf 2.2.2. Lemma and Definition.}} 
{\it
 Let\/ $\mP_k=\left(p_{_{\scriptstyle ij}}^{(k)}\right)_{l\times m}$ $(k\in\mathbb
 P)$, $\mP=\left(p_{_{\scriptstyle ij}}\right)_{l\times m}\in\mathbb R^{l\times
 m}$. Then
 $$
 p_{_{\scriptstyle ij}}^{(k)}\rightarrow p_{_{\scriptstyle ij}}\quad (k\rightarrow +\infty)
$$
$($for $i=1$, $\dots$, $l$ and $j=1$, $\dots$, $m)$ if and only if
$$
\|\mP_k-\mP\|\rightarrow0\quad (k\rightarrow +\infty),
$$
where\/ $\mP$ is unique. } We say that the sequence $\{\mP_k\}$ of
matrices has the {\it limit}\/ $\mP$ or $\{\mP_k\}$ {\it
approaches}\/ $\mP$, denoted by $\displaystyle\lim_{k\to
+\infty}\mP_k=\mP$ or $\mP_k\rightarrow\mP$
($k\rightarrow+\infty$).

\paragraph{{\bf 2.2.3. Lemma.}} 
(1) {\it If\/ $\mP$, $\mQ\in\mathbb R^{l\times
 m}$, then
 $$
 \|\mP+\mQ\|\leqslant\|\mP\|+\|\mQ\|.
 $$
 }
 (2) {\it
 If\/ $\mP=\left(p_{_{\scriptstyle ij}}\right)_{l\times m}\in\mathbb R^{l\times m}$ and
 $\mQ=\left(q_{_{\scriptstyle \alpha\beta}}\right)_{m\times s}\in\mathbb R^{m\times
 s}$, then
 $$
 \|\mP\mQ\|\leqslant\|\mP\|\|\mQ\|.
 $$
 }
\par\noindent{\it Proof.}   (1) is obvious and (2) follows from
\begin{align*}
 \sum_{i=1}^l\sum_{\beta=1}^s\left|\sum_{j=1}^mp_{_{\scriptstyle ij}}q_{_{j\beta}}\right|
 &\leqslant\sum_{i=1}^l\sum_{\beta=1}^s\sum_{j=1}^m|p_{_{ij}}||q_{_{j\beta}}|\\
 &\leqslant\sum_{i=1}^l\sum_{j=1}^m\sum_{\alpha=1}^m\sum_{\beta=1}^s
 |p_{_{\scriptstyle ij}}||q_{_{\scriptstyle \alpha\beta}}|\\
 &=\left(\sum_{i=1}^l\sum_{j=1}^m|p_{_{\scriptstyle ij}}|\right)
 \left(\sum_{\alpha=1}^m\sum_{\beta=1}^s|q_{_{\scriptstyle \alpha\beta}}|\right).
\end{align*}

\paragraph{{\bf 2.2.4. Lemma.}} 
 {\it
 Let\/ $\mP_k$ $(k\in\mathbb P)$, $\mP\in\mathbb R^{l\times m}$.
 If\/ $\mP_k\rightarrow\mP$ and\/ $\mQ_k\rightarrow\mQ$, then\/
 $\mP_k\mQ_k\rightarrow\mP\mQ$.
 }
\par\noindent{\it Proof.}   From $\mQ_k\rightarrow\mQ$ and
$\|\mQ_k\|\leqslant\|\mQ_k-\mQ\|+\|\mQ\|$ we know that$\|\mQ_k\|$
is bounded, i.e., there exists a constant $M>0$ such that
$\|\mQ_k\|\leqslant M$. By Lemma~2.2.3 we have
\begin{align*}
 \|\mP_k\mQ_k-\mP\mQ\|&=\|(\mP_k-\mP)\mQ_k+\mP(\mQ_k-\mQ)\|\\
 &\leqslant \|\mP_k-\mP\|\|\mQ_k\|+\|\mP\|\|\mQ_k-\mQ\|\\
 &\leqslant M\|\mP_k-\mP\|+\|\mP\|\|\mQ_k-\mQ\|\rightarrow0,
\end{align*}
which implies $\mP_k\mQ_k\rightarrow\mP\mQ$.


\paragraph{{\bf 2.2.5. Lemma.}} 
 {\it
 Let\/ $\mP_k=\left(p_{_{\scriptstyle \alpha\beta}}^{(k)}\right)_{l\times m}\in\mathbb R^{l\times m}$
 and\/ $\|\mP_k\|\leqslant M$, where\/ $M>0$ is a constant. Then there
 exists a subsequence\/ $\{k_j\}_{j=1}^{\infty}$ of the natural
 number sequence and\/ $\mP\in\mathbb R^{l\times m}$ such that\/
 $\mP_{k_j}\rightarrow\mP$ $(j\rightarrow+\infty)$.
 }
\par\noindent{\it Proof.}   From $\|\mP_k\|\leqslant M$ we know that
for all $\alpha$ and $\beta$ ($\alpha=1$, $\dots$, $l$; $\beta=1$,
$\dots$, $m$) the number sequences $\left\{p_{_{\scriptstyle
\alpha\beta}}^{(k)}\right\}_{k=1}^{\infty}$ are bounded. When
$(\alpha,\beta)=(1,1)$, we can get a number sequence
$\left\{p_{_{\scriptstyle 11}}^{(k_j)}\right\}_{j=1}^{\infty}$,
which has a limit~$p_{_{\scriptstyle 11}}$. Then consider
$(\alpha,\beta)=(1,2)$, i.e., a number sequence
$\left\{p_{_{\scriptstyle 12}}^{(k_j)}\right\}_{j=1}^{\infty}$,
which is of course bounded. Hence there exists a subsequence
$\{k_{j_{_i}}\}_{j=1}^{\infty}$ of $\{k_j\}_{j=1}^{\infty}$ such
that $p_{_{\scriptstyle 12}}^{(k_{j_{_i}})}\rightarrow
p_{_{\scriptstyle 12}}\in\mathbb R$ ($i\to+\infty$). Here we still
have $p_{_{\scriptstyle 11}}^{(k_{j_{_i}})}\rightarrow
p_{_{\scriptstyle 11}}$. So we might write $k_j$ instead of
$k_{j_{_i}}$ for simplicity. By induction we can obtain
$p_{_{\scriptstyle \alpha\beta}}$ for all $\alpha$ and $\beta$
($\alpha=1$, $\dots$, $l$; $\beta=1$, $\dots$, $m$) so that
$p_{_{\scriptstyle \alpha\beta}}^{(k_j)}\rightarrow
p_{_{\scriptstyle \alpha\beta}}$ ($j\to+\infty$), where
$\{k_j\}_{j=1}^{\infty}$ is some subsequence of the natural number
sequence. Let $\mP=\left(p_{_{\scriptstyle
\alpha\beta}}\right)_{l\times m}$. Then $\mP_{k_j}\rightarrow\mP$
$(j\rightarrow+\infty)$.


\paragraph{{\bf 2.2.6. Definition.}} 
Let $B$ be a nonempty subset of~$\mathbb R^n$. If
$\{A_k\}_{k=1}^{\infty}$ is a sequence of nonempty subsets
of~$\mathbb R^n$ that satisfies $A_k=\{a_k(b):\ b\in B\}$ ($k=1$,
$2$, $\dots$), then $\{A_k\}_{k=1}^{\infty}$ is said to be a {\it
$B$-index sequence}. If further there exists a positive number~$K$
so that when $k>K$ we have $d(a_k(b), f(b))<\varepsilon$ for all
$b\in B$, then we say that $\{a_k(b)\}_{k=1}^{\infty}$ is {\it
uniformly convergent} to $f(b)$ on~$B$, which is denoted by
$\displaystyle d(a_k(b), f(b))\overset{B}{\rightrightarrows} 0$
($k\rightarrow+\infty$) or $\displaystyle
a_k(b)\overset{B}{\rightrightarrows} f(b)$
($k\rightarrow+\infty$).

\paragraph{{\bf 2.2.7. Lemma.}} 
 {\it
 Let\/ $\{A_k\}_{k=1}^{\infty}$ be a\/ $B$-index sequence, where\/ $B$ is a nonempty subset of\/~$\mathbb
 R^n$ and\/ $A_k=\{a_k(b):\: b\in B\}$. If\/  $\displaystyle a_k(b)\overset{B}{\rightrightarrows} f(b)$
 $(k\rightarrow+\infty)$, then
 $$
 h(A_k,f(B))\rightarrow0\quad (k\rightarrow+\infty).
 $$
 }
\par\noindent{\it Proof.}   Given any $\varepsilon>0$, there exists
$K>0$ such that when $k>K$ we have
$$d(a_k(b), f(b))<\varepsilon$$
 for all $b\in B$. Hence for any $b\in B$ we deduce
 $$
 d(a_k(b), f(B))=\inf_{x\in B}\{d(a_k(b), f(x))\}
 \leqslant d(a_k(b), f(b))<\varepsilon
 $$
and
 $$
 d(A_k, f(b))=\inf_{x\in B}\{d(a_k(x), f(b))\}
 \leqslant d(a_k(b), f(b))<\varepsilon.
 $$
 These imply
 $$
 \sup_{b\in B} d(a_k(b), f(B))\leqslant\varepsilon
 $$
 and
 $$
 \sup_{b\in B} d(A_k, f(b))\leqslant\varepsilon.
 $$
Therefore
$$
h(A_k,f(B))=\max\left\{\sup_{b\in B} d(a_k(b), f(B)),\ \sup_{b\in
B} d(A_k, f(b))\right\}\leqslant\varepsilon.
$$


By Linear Algebra it follows that
\paragraph{{\bf 2.2.8. Lemma.}} 
 {\it
 The transformation\/ $\varphi:\: \mathbb R^n\rightarrow\mathbb R^n$
 is an isometry if and only if there are an orthogonal
 transformation\/ $\sigma:\: \mathbb R^n\rightarrow\mathbb R^n$ and a
 translation\/ $\tau:\: \mathbb R^n\rightarrow\mathbb R^n$
 satisfying\/
 $\varphi=\sigma\nolinebreak{\circ}\nolinebreak\tau$
 $(\sigma\nolinebreak{\circ}\nolinebreak\tau$ denotes the composite of\/ $\sigma$ and\/ $\tau)$.
 }

\paragraph{{\bf 2.2.9. Lemma.}} 
 {\it
 If\/ $\mO\in\mathbb R^{n\times n}$ is an orthogonal matrix then
 $$|\mO|\leqslant n^2.$$
 }
\par\noindent{\it Proof.}   Let $\mO=(a_{_{\scriptstyle ij}})_{n\times
n}$. Then by $\mO\mO'=\mathsf I_n$ (the identity matrix), where
$\mO'$ denotes the transpose of $\mO$, we get
$$
\sum_{j=1}^n a_{_{\scriptstyle ij}}^2=1\quad(\text{$i=1$, $\dots$,
$n$}).
$$
 Hence
$a_{_{\scriptstyle ij}}^2\leqslant 1$, i.e., $|a_{_{\scriptstyle
ij}}|\leqslant 1$. Therefore
$$
\|\mO\|=\sum_{i=1}^n\sum_{j=1}^n |a_{_{\scriptstyle ij}}|\leqslant
n^2.
$$


({\it Continuation of the proof of Theorem~2.2.}) Suppose that
$A$, $B\in\mathcal C(\mathbb R^n)$ and $\wth(\wt A,\wt B)=0$. Then
for $k\in\mathbb P$ there are $B_k\sim B$ such that
$$
H(A,B_k)<\frac1k.
$$
 So there exist isometries $\varphi_{_{\scriptstyle k}}
 =\sigma_{_{\scriptstyle k}}\nolinebreak\circ\nolinebreak\tau_{_{\scriptstyle k}}$
 with $B_k=\varphi_{_{\scriptstyle k}}(B)$, where $\sigma_{_{\scriptstyle k}}$
 are orthogonal transformations and $\tau_{_{\scriptstyle k}}$ are translations. Hence
$$
B_k=\{\varphi_{_{\scriptstyle k}}(b):\: b\in B\},
$$
and $\{B_k\}_{k=1}^{\infty}$ is a $B$-index sequence. Since
$$
\varphi_{_{\scriptstyle k}}(b)=\sigma_{_{\scriptstyle
k}}(\tau_{_{\scriptstyle k}}(b))=(b+t_k)\mO_k,
$$
where $t_k\in\mathbb R^n$ and $\mO_k$ is an orthogonal matrix, we
get
\begin{align*}
|b+t_k|&=|\sigma_{_{\scriptstyle
k}}(b+t_k)|=|\varphi_{_{\scriptstyle k}}(b)|
\leqslant|\varphi_{_{\scriptstyle k}}(b)-a|+|a|\\
&\leqslant1+d(\varphi_{_{\scriptstyle k}}(b),A)+|a|\leqslant
1+h(A,B_k)+|a|<2+|a|,
\end{align*}
where $|x|$ denotes the length of the vector corresponding to
$x\in\mathbb R^n$ and $a\in A$ is chosen suitably. Therefore
$$
|t_k|\leqslant|-b|+|b+t_k|\leqslant |b|+2+|a|\leqslant 2R+2,
$$
where $R=\max\{|a|,\, |b|:\: a\in A,\ b\in B\}$. Thus there is a
subsequence $\{t_{k_j}\}_{j=1}^{\infty}$
of~$\{t_{k}\}_{k=1}^{\infty}$ with $t_{k_j}\to t\in\mathbb R^n$
($j\to +\infty$). By Lemma~2.2.9 we have $\|\mO_{k_j}\|\leqslant
n^2$. Again by Lemma~2.2.5 we obtain a subsequence
of~$\{\mO_{k_j}\}_{j=1}^{\infty}$, still denoted
by~$\{\mO_{k_j}\}_{j=1}^{\infty}$, which approaches $\mO\in\mathbb
R^{n\times n}$. Hence $\mO_{k_j}'\to\mO'$ ($j\to+\infty$). By
Lemma~2.2.4 it follows
$$
\mO_{k_j}\mO_{k_j}'\to\mO\mO'.
$$
As $\mO_{k_j}\mO_{k_j}'=\mathsf I_n$, we have $\mO\mO'=\mathsf
I_n$, which means $\mO$ is orthogonal. Let
$\varphi=\sigma\nolinebreak\circ\nolinebreak\tau$, where
$\sigma(x)=x\mO$ and $\tau(x)=x+t$. Then $\varphi$ is an isometry.

Since
\begin{align*}
d(\varphi_{_{\scriptstyle k_j}}(b),\varphi(b))&=|\varphi_{_{\scriptstyle k_j}}(b)-\varphi(b)|\\
&\leqslant\|\varphi_{_{\scriptstyle k_j}}(b)-\varphi(b)\|
=\|\tau_{_{\scriptstyle k_j}}(b)\mO_{k_j}-\tau(b)\mO\|\\
&\leqslant\|\tau_{_{\scriptstyle k_j}}(b)-\tau(b)\|\|\mO_{k_j}\|+\|\tau(b)\|\|\mO_{k_j}-\mO\|\\
&\leqslant n^2\|t_{k_j}-t\|+M\|\mO_{k_j}-\mO\|\to
0\quad(j\to+\infty),
\end{align*}
where $M$ relies on $A$ and~$B$ but does not rely on $b\in B$, we
have
$$
 \varphi_{_{\scriptstyle k_j}}(b)\overset{B}{\rightrightarrows}\varphi(b)\quad(j\to+\infty).
$$
By Lemma~2.2.7 it follows
$$
 h(\varphi_{_{\scriptstyle k_j}}(B),\varphi(b))\to0\quad(j\to+\infty).
$$
Therefore
$$
h(A,\varphi(B))\leqslant h(A,B_{k_j})+h(\varphi_{_{\scriptstyle
k_j}}(B),\varphi(b))\to0\quad(j\to+\infty).
$$
So $h(A,\varphi(B))=0$. This implies $A=\varphi(B)$, i.e., $\wt
A=\wt B$.

Now let us prove that the metric~$\wth$ is complete. Suppose
$\{\wt A_k\}_{k=1}^\infty$ is a Cauchy sequence in~$\wt{\mathcal
C}(\mathbb R^n)$. Then for an arbitrary $\varepsilon>0$, there
exists $K>0$ such that when $j$, $k\geqslant K$ one has $\wth(\wt
A_j,\wt A_k)<\varepsilon$.

\paragraph{{\bf (a) Claim.}}
 {\it
 One may select a subsequence\/ $\{\wt A_{k_i}\}_{i=1}^\infty$ of\/ $\{\wt
 A_k\}_{k=1}^\infty$ such that
 $$
 \wth(\wt A_{k_i},\wt A_{k_{i+1}})<2^{-i}
 $$
 for $i=1$, $2$, $\dots$.
 }
\par\noindent{\it Proof.}   We can easily know that there exists a
subsequence $\{\wt A_{k_i}\}_{i=1}^\infty$
($k_1<k_2<\dots<k_i<\dotsb$) so that
$$
\wth(\wt A_{k_i},\wt A_k)<2^{-i}
$$
hold for all $k>k_i$, where $i=1$, $2$, $\dots$. Taking
$k=k_{i+1}$ the claim follows.

\paragraph{{\bf (b) Claim.}}
 {\it
 There exist\/ $B_j\in\mathcal C(\mathbb R^n)$ so that\/ $B_j\sim
 A_{k_j}$ and
 $$
 h(B_j,B_{j+1})<2^{-j}\quad\text{$(j=1$, $2$, $\dots)$}.
 $$
 }
\par\noindent{\it Proof.}   Let $B_1=A_{k_1}$. From $\wth(\wt
A_{k_1},\wt A_{k_2})<2^{-1}$ we get that there exists $B_2\sim
A_{k_2}$ such that
$$
h(A_{k_1}, B_2)<2^{-1}.
$$
Suppose there exist $B_i\sim A_{k_i}$ ($i=1$, $2$, $\dots$, $j$)
such that
 $$
 h(B_i,B_{i+1})<2^{-i}\quad\text{($i=1$, $2$, $\dots$, $j-1$)}.
 $$
 Then by Claim~(a) we have
 $$
 \wth(\wt B_j,\wt A_{k_{j+1}})=\wth(\wt A_{k_j},\wt
 A_{k_{j+1}})<2^{-j}.
 $$
So we may choose $B_{j+1}\sim A_{j+1}$ such that
 $$
 h(B_j,B_{j+1})<2^{-j}.
 $$
By induction Claim~(b) holds.

\paragraph{{\bf (c) Claim.}}
 {\it
Let\/ $C_k=\mathrm{cl}\left(\bigcup_{j\geqslant k}B_j\right)$
$(\mathrm{cl}\,(A)$ indicates the closure of\/~$A)$ and\/
$C=\bigcap_{k=1}^\infty C_k$. Then
$$
h(C_k,B_k)\to0\quad\text{and}\quad h(C_k,C)\to0\quad(k\to+\infty).
$$
 }
\par\noindent{\it Proof.}   Let $k\in\mathbb P$. It is obvious that
$C_k\supseteq C_{k+1}$
 and $C$ is a nonempty closed set. By Claim~(b) it follows that
$$
h(B_j,B_{j+p})\leqslant\sum_{i=j}^{j+p-1}
h(B_i,B_{i+1})<\sum_{i=j}^{j+p-1} 2^{-i}<2^{-j+1}
$$
for all $p\in\mathbb P$. Hence $\{B_j\}_{j=1}^\infty$ is a Cauchy
sequence. As $h(B_1,B_j)<1$
 we know that $B_j\subseteq\mathcal
P(B_1,1)$ ($j\in\mathbb P$). So
$$
C_k=\mathrm{cl}\left(\bigcup_{j\geqslant
k}B_j\right)\subseteq\mathscr P(B_1,1),
$$
which implies $C_k\in\mathcal C(\mathbb R^n)$ ($k\in\mathbb P$)
and $C\in\mathcal C(\mathbb R^n)$. Similarly we have
$B_j\subseteq\mathscr P(B_k, 2^{-k+1})$ ($j\geqslant k$). Thus
$$
C_k=\mathrm{cl}\left(\bigcup_{j\geqslant
k}B_j\right)\subseteq\mathscr P(B_k, 2^{-k+1}).
$$
This implies
$$
h(C_k,B_k)\leqslant 2^{-k+1}\to0.
$$
Therefore
\begin{align*}
h(C_k,C_{k+p})&\leqslant
h(C_k,B_k)+h(B_k,B_{k+p})+h(B_{k+p},C_{k+p})\\
&<2^{-k+1}+2^{-k+1}+2^{-k-p+1}<2^{-k+3}.
\end{align*}
So we have
$$
C_k\subseteq\mathscr P(C_{k+p},2^{-k+3})\quad(p\in\mathbb P).
$$
Suppose $x\in C_k$. Then there is $x_p\in C_{k+p}$ so that
$d(x,x_p)\leqslant 2^{-k+3}$. We choose a subsequence
$\{x_{p_{_i}}\}_{i=1}^\infty$ of $\{x_p\}_{p=1}^\infty$ such that
$x_{p_{_i}}\to x_0\in\mathbb R^n$ ($i\to+\infty$). Then $x_0\in
C_{k+p_{_i}}$. Hence $x_0\in C$ and
$$
d(x,x_0)\leqslant 2^{-k+3},
$$
which means $x\in\mathscr P(C,2^{-k+3})$. Consequently
$C_k\subseteq\mathscr P(C,2^{-k+3})$. Therefore
$$
h(C_k,C)\leqslant 2^{-k+1}\to0,
$$
and Claim~(c) follows.

Given $\varepsilon>0$, by Claim~(c) there exists $K>0$ so that we
may take a sufficiently great~$j$, when $k>K$ we have
\begin{align*}
\wth(\wt A_k,\wt C)&\leqslant\wth(\wt A_k,\wt A_{k_j})+\wth(\wt
B_j,\wt C_j)+\wth(\wt C_j,\wt C)\\
 &\leqslant\wth(\wt A_k,\wt
A_{k_j})+h(\wt B_j,\wt C_j)+h(\wt C_j,\wt
C)<\frac{\varepsilon}3+\frac{\varepsilon}3+\frac{\varepsilon}3=\varepsilon.
\end{align*}
Finally $\wth(\wt A_k,\wt C)\to0$ ($k\to+\infty$). This completes
the proof of Theorem~2.2. \hfill $\Box$


\paragraph{{\bf 2.3.} {\it Remark}.}
 (1) The metric space
$\bigl(\wt{\mathcal{C}}(\mathbb{R}^n),\wth\bigr)$
  is separable.

 (2) We may regard $\bigl({\mathcal{C}}(\mathbb{R}^n),\wth\bigr)$
 as $\bigl(\wt{\mathcal{C}}(\mathbb{R}^n),\wth\bigr)$
and ``$\sim$'' as ``$=$'' in the definition of metric spaces.
Under this convention, we may say that
$\bigl(\mathcal{C}(\mathbb{R}^n),\wth\bigr)$ is a complete metric
space.

 (3) By a similar reasoning to that in the proof of Theorem~2.2 we may show
$$
\wth(A,B)=\min\left\{h(A_1,B_1):\: A_1\sim A,\ B_1\sim B\right\}.
$$
Starting from this conclusion we can also deduce Theorem~2.2.

\paragraph{{\bf 2.4. Rigid shape differences.}}
\paragraph{{\bf 2.4.1. Definition.}}
An isometry~$\psi$ in~$\mathbb R^n$ is said to be {\it rigid\/} if
$\psi=\sigma\nolinebreak\circ\nolinebreak\tau$, where $\tau$ is a
translation and $\sigma$ is a rigid orthogonal transformation or
rotation, i.e., $\det\sigma=1$ ($\det\sigma:=\det\mO$ where $\mO$
is a matrix of~$\sigma$ under some orthonormal basis of~$\mathbb
R^n$). We also call a rigid isometry a {\it rigid motion}.

\paragraph{{\bf 2.4.2. Definition.}}
Two nonempty subsets $A$ and~$B$ of~$\mathbb R^n$ are {\it rigidly
equivalent}\/ if there exists a rigid isometry~$\psi$ in~$\mathbb
R^n$ such that $B=\psi(A)$, denoted by $A\simeq B$.

\paragraph{{\bf 2.4.3. Lemma.}}
{\it
 Let\/ $a\in\mathbb R^n$ and define\/ $\tau_{_{\scriptstyle a}}:=x+a$ for\/ $x\in\mathbb
 R^n$. Suppose\/ $\sigma$ is an orthogonal transformation, $\psi_{_{\scriptstyle 1}}$
 and\/ $\psi_{_{\scriptstyle 2}}$ are rigid isometries and\/ $b\in\mathbb R^n$. Then
\begin{itemize}
 \item[{\rm(1)}] $\tau_{_{\scriptstyle a}}\circ\tau_{_{\scriptstyle b}}=\tau_{_{\scriptstyle a+b}}$;
 \item[{\rm(2)}] $\tau_{_{\scriptstyle a}}^{-1}=\tau_{_{\scriptstyle -a}}$;
 \item[{\rm(3)}] $\tau_{_{\scriptstyle a}}\circ\sigma=\sigma\circ\tau_{_{\scriptstyle \sigma^{-1}(a)}}$;
 \item[{\rm(4)}] $\psi_{_{\scriptstyle 1}}^{-1}$ and\/ $\psi_{_{\scriptstyle 2}}\circ\psi_{_{\scriptstyle 1}}$
  are also rigid isometries.
\end{itemize}
  }
\par\noindent{\it Proof.}  (1) and (2) are obvious. For (3) and
(4), letting $x\in\mathbb R^n$ we have
$$
\tau_{_{\scriptstyle
a}}\circ\sigma(x)=\sigma(x)+a=\sigma(x+\sigma^{-1}(a))=\sigma\circ\tau_{_{\scriptstyle
\sigma^{-1}(a)}}(x).
$$
Let $\psi_{_{\scriptstyle 1}}=\sigma_{_{\scriptstyle
1}}\circ\tau_{_{\scriptstyle a}}$ and $\psi_{_{\scriptstyle
2}}=\sigma_{_{\scriptstyle 2}}\circ\tau_{_{\scriptstyle b}}$,
where $\sigma_{_{\scriptstyle 1}}$ and $\sigma_{_{\scriptstyle
2}}$ are rigid orthogonal transformations. Then
$\sigma_{_{\scriptstyle 1}}^{-1}$ and $\sigma_{_{\scriptstyle
2}}\circ\sigma_{_{\scriptstyle 1}}$ are rigid orthogonal
transformations,
$$
\psi_{_{\scriptstyle 1}}^{-1}=\tau_{_{\scriptstyle
a}}^{-1}\circ\sigma_{_{\scriptstyle 1}}^{-1}=\tau_{_{\scriptstyle
-a}}\circ\sigma_{_{\scriptstyle 1}}^{-1}=\sigma_{_{\scriptstyle
1}}^{-1}\circ\tau_{_{\scriptstyle {-\sigma_{_1}}(a)}}
$$
and
\begin{align}
\psi_{_{\scriptstyle 2}}\circ\psi_{_{\scriptstyle
1}}&=(\sigma_{_{\scriptstyle 2}}\circ\tau_{_{\scriptstyle
b}})\circ(\sigma_{_{\scriptstyle 1}}\circ\tau_{_{\scriptstyle a}})
=\sigma_{_{\scriptstyle 2}}\circ(\tau_{_{\scriptstyle b}}
\circ\sigma_{_{\scriptstyle 1}})\circ\tau_{_{\scriptstyle a}}\notag\\
&=\sigma_{_{\scriptstyle 2}}\circ(\sigma_{_{\scriptstyle
1}}\circ\tau_{_{\scriptstyle
\sigma_{_1}^{-1}(b)}})\circ\tau_{_{\scriptstyle a}}
=(\sigma_{_{\scriptstyle 2}}\circ\sigma_{_{\scriptstyle
1}})\circ\tau_{_{\scriptstyle
a+\sigma_{_1}^{-1}(b)}}.\tag*{$\Box$}
\end{align}


\paragraph{{\bf 2.4.4.}}
By Lemma~2.4.3(4) we see that the relation~``$\simeq$'' is an
equivalent relation. The equivalence class of~$A$ is
denoted~$\overline{A}$. And we denote
$$
\overline{\mathcal{C}}(\mathbb{R}^n)=\left\{\overline{C}:\: C
\in\mathcal{C}(\mathbb{R}^n)\right\}.
$$

\paragraph{{\bf 2.4.5. Definition.}}
For nonempty subsets $A$ and~$B$ of~\/$\mathbb R^n$ we define
$$
\overline{h}(A,B):=\inf\left\{h(\varphi(A),\psi(B):\:
\textrm{$\varphi$ and $\psi$ are rigid isometries in~$\mathbb
R^n$}\right\}.
$$
By Lemmas~2.1.2 and~2.4.3(4) we have the following expressions:
\begin{align*}
\overline{h}(A,B):&=\inf\left\{h(A_1,B_1):\: A_1\simeq A,\ B_1\simeq B\right\}\\
 &=\inf\left\{h(A_0,B_1):\: B_1\simeq B\right\}\quad\textrm{(where $A_0\simeq A$)}\\
 &=\inf\left\{h(A_1,B_0):\: A_1\simeq A\right\}\quad\textrm{(where $B_0\simeq B$)}.
\end{align*}

If $A_1\simeq A_2$ and $B_1\simeq B_2$, then
$\overline{h}(A_1,B_1)=\overline{h}(A_2,B_2)$. So we may define
$$
\overline{h}(\overline{A},\overline{B}):=\overline{h}(A,B),
$$
which is called the ({\it absolute}) {\it rigid shape difference}
between $\wt A$ and~$\wt B$ or between $A$ and~$B$.

\paragraph{{\bf 2.5. Theorem.}} 
 {\it
$\bigl(\overline{\mathcal{C}}(\mathbb{R}^n),\overline{h}\bigr)$ is
a complete metric space.
 }

\bigskip
The proof of Theorem~2.5 is just similar to that of Theorem~2.2,
where we use Definition~2.4.1 instead of Lemma~2.2.8.\hfill $\Box$

\paragraph{{\bf 2.6. Relative shape differences.}} 

\paragraph{{\bf 2.6.1. Definition.}} 
Let $r>0$.

(1) A transformation $S:\: \mathbb R^n\to\mathbb R^n$ is a {\it
similitude} or {\it $r$-similitude} if
$$
d(S(x),S(y))=rd(x,y)
$$
for all $x$, $y\in\mathbb R^n$, and $r$ is called the {\it
Lipschitz constant} of~$S$.

(2) Two nonempty subsets $A$ and~$B$ of~$\mathbb R^n$ are {\it
similar} if there exists a similitude $S:\: \mathbb
R^n\rightarrow\mathbb R^n$ such that $B=S(A)$, denoted by
$A\widehat{\sim} B$. We easily know that this relation is an
equivalence relation. The equivalence class of~$A$ is
denoted~$\widehat{A}$.

(3) The transformation $\mu_r(x)=rx$ ($x\in\mathbb R^n$) is called
a {\it homothety} or {\it $r$-homothety}.

\paragraph{{\bf 2.6.2. Definition.}} 
Let $A$ be a nonempty subset of~$\mathbb R^n$. The {\it diameter}
of~$A$ is
$$
|A|:=\sup\{d(x,y):\: x,\, y\in A\}.
$$
 We define the {\it
radius}\/ $\mathrm{r}(A)$ of~$A$ by
$$
\mathrm{r}(A):=\inf_{x\in\mathbb R^n}\left\{\sup_{a\in
A}d(x,a)\right\}.
$$
We also denote $\mathrm{r}(A)$ by $\mathrm{r}_{_{\scriptstyle
A}}$. Denote
$$
\rho A:=\{\rho a:\: a\in A\}
$$
for $\rho\geqslant0$. Generally $\mathrm{r}(A)\neq\frac12|A|$.

\bigskip
For the need later on, we list the following conclusions.

\paragraph{{\bf 2.6.3. Proposition.}} 
{\it
 Let\/ $A$ and\/~$B$ be nonempty bounded subsets of\/~$\mathbb R^n$. 
\begin{itemize}
 \item[{\rm(1)}] $\mathrm{r}(S(A))=r\,\mathrm{r}(A)$
  for a similitude $S$ of Lipschitz constant\/~$r$;
   specially\/ $\mathrm{r}(\sigma(A))=\mathrm{r}(A)$
   for an isometry\/~$\sigma$.
 \item[{\rm(2)}] $\left||A|-|B|\right|\leqslant 2\wth(A,B)$.
 \item[{\rm(3)}] $|\mathrm{r}(A)-\mathrm{r}(B)|\leqslant
 h(A,B)$.
 \item[{\rm(4)}] $|\mathrm{r}(A)-\mathrm{r}(B)|\leqslant\wth(A,B)$.
 \item[{\rm(5)}] There exists\/
 $x_{_{\scriptstyle 0}}\in\mathbb R^n$ such that
 $$
 \sup_{a\in A}d(x_{_{\scriptstyle 0}},a)=\mathrm{r}(A).
 $$
 {\rm We call $x_{_{\scriptstyle 0}}$ a {\it center}\/ of~$A$ (it is possible that $x_0\notin A$).}
 \item[{\rm(6)}] If\/ $0<\rho<+\infty$ then $h(A,\rho
 A)\leqslant|\rho-1|\,\mathrm{r}(A)$.
\end{itemize}
 }
\par\noindent{\it Proof.}  (1) is obvious and (4) follows from (1)
and~(3). The proof of~(2) is similar to that of~(4) and easier.
Now we prove (3), (5) and~(6).

(i) Let $\varepsilon$ be an arbitrary positive number. Given any
$x\in\mathbb R^n$ we have
$$
\mathrm{r}(A)\leqslant\sup_{a\in A}d(x,a).
$$
Thus there exists $a\in A$ such that
$$
\mathrm{r}(A)-\varepsilon<d(x,a).
$$
Choose $b_0\in B$ so that
$$
d(a,b_0)<d(a,B)+\varepsilon.
$$
Then
\begin{align*}
\mathrm{r}(A)-\varepsilon&<d(x,b_0)+d(b_0,a)\\
&<\sup_{b\in B}d(x,b)+d(a,B)+\varepsilon\leqslant\sup_{b\in
B}d(x,b)+h(A,B)+\varepsilon.
\end{align*}
It follows that
$$
\mathrm{r}(A)-2\varepsilon\leqslant\inf_{x\in\mathbb
R^n}\left\{\sup_{b\in
B}d(x,b)\right\}+h(A,B)=\mathrm{r}(B)+h(A,B).
$$
Therefore
$$
\mathrm{r}(A)\leqslant\mathrm{r}(B)+h(A,B).
$$
Hence (3) is true.

(ii) For any positive integers~$k$ there exist $x_{_{\scriptstyle
k}}\in\mathbb R^n$ such that
$$
\sup_{a\in A}d(x_{_{\scriptstyle k}},a)<\mathrm{r}(A)+\frac1k.
$$
Thus there exist $x_{_{\scriptstyle 0}}\in\mathbb R^n$ and a
subsequence $\{x_{_{\scriptstyle k_j}}\}_{j=1}^\infty$ of
$\{x_{_{\scriptstyle k}}\}_{k=1}^\infty$ such that
$$
d(x_{_{\scriptstyle k_j}},x_{_{\scriptstyle
0}})\to0\quad(j\to+\infty).
$$
Therefore
$$
\sup_{a\in A}d(x_{_{\scriptstyle 0}},a)\leqslant\mathrm{r}(A).
$$

(iii) Suppose $x_{_{\scriptstyle 0}}$ is a center of~$A$. Let
$$
B=A-x_{_{\scriptstyle 0}}:=\{a-x_{_{\scriptstyle 0}}:\:a\in A\}.
$$
Then
$$
\rho B=\rho A-\rho x_{_{\scriptstyle 0}}:=\{\rho a-\rho
x_{_{\scriptstyle 0}}:\:a\in A\}.
$$
Since
$$
d(b,\rho B),\ d(\rho b,B)\leqslant|b-\rho
b|=|\rho-1|\,|b|\leqslant|\rho-1|\,\mathrm{r}(A),
$$
we have
\begin{align}
\wth(A,\rho A)&=\wth(B,\rho B)\notag\\
&\leqslant\sup\{d(b,\rho B),\,d(\rho b,B):\: b\in
B\}\leqslant|\rho-1|\,\mathrm{r}(A).\tag*{$\Box$}
\end{align}


\paragraph{{\bf 2.6.4. Lemma.}} 
{\it
 Let\/ $A$ and\/~$B$ be nonempty bounded subsets of\/~$\mathbb R^n$.
 Then\/ $\widehat{A}=\widehat{B}$
  if and only if\/ $\displaystyle \frac A{\mathrm{r}(A)}\sim \frac B{\mathrm{r}(B)}$,
   where if\/ $A$ is a singleton then we treat
   $\displaystyle \frac A{\mathrm{r}(A)}$ as\/~$A$.
 }
\par\noindent{\it Proof.}   We easily know that an isometry is a
$1$-similitude, an $r$-homothety is an $r$-similitude ($r>0$), and
if $S_1$ and~$S_2$ are an $r_1$-similitude and an $r_2$-similitude
respectively then $S_1\circ S_2$ is an $r_1r_2$-similitude. When
$A$ or~$B$ is a singleton, the lemma is obviously true. Now
suppose neither $A$ nor~$B$ is a singleton.

If $\displaystyle \frac A{\mathrm{r}(A)}\sim \frac
B{\mathrm{r}(B)}$, then  $\displaystyle \frac
B{\mathrm{r}(B)}=\varphi\left(\frac A{\mathrm{r}(A)}\right)$,
where $\varphi$ is an isometry. Hence
$$
B={\mathrm{r}(B)}\,\varphi\left(\frac1{\mathrm{r}(A)}
A\right)=\left(\mu_{_{\scriptstyle\mathrm{r}(B)}}\circ\varphi\circ
\mu_{_{\scriptstyle(\mathrm{r}(A))^{-1}}}\right)(A),
$$
which means $A\widehat{\sim} B$.

If $\widehat{A}=\widehat{B}$, i.e., $A\widehat{\sim} B$, then
$B=S(A)$, where $S$ is an $r$-similitude ($r>0$), so
$\mathrm{r}(B)=r\,\mathrm{r}(A)$. Consequently
$$
\frac
B{\mathrm{r}(B)}=\frac1{r\,\mathrm{r}(A)}\,S(A)=\left(\mu_{_{\scriptstyle(r\,\mathrm{r}(A))^{-1}}}\circ
S \circ \mu_{_{\scriptstyle(\mathrm{r}(A))}}\right)\left(\frac
A{\mathrm{r}(A)}\right),
$$
where $\displaystyle
\mu_{_{\scriptstyle(r\,\mathrm{r}(A))^{-1}}}\circ S \circ
\mu_{_{\scriptstyle(\mathrm{r}(A))}}$ is an isometry.\hfill $\Box$


\paragraph{{\bf 2.6.5. Definition.}} 
Let $A$ and~$B$ be nonempty bounded subsets of~$\mathbb R^n$.
Define
$$
\widehat{h}(A,B):=\wth\left(\frac A{\mathrm{r}(A)},\,\frac
B{\mathrm{r}(B)}\right)
$$
and
$$
\widehat{h}(\widehat{A},\widehat{B}):=\widehat{h}(A,B),
$$
which is called the {\it relative shape difference} between
$\widehat{A}$ and~$\widehat{B}$ or between $A$ and~$B$.

\paragraph{{\it Remark}.} By Lemma~2.6.4 we see that if $A\widehat{\sim} A_1$
and $B\widehat{\sim} B_1$ then
$$
\frac A{\mathrm{r}(A)}\sim \frac {A_1}{\mathrm{r}(A_1)}
\quad\text{and}\quad
 \frac B{\mathrm{r}(B)}\sim \frac {B_1}{\mathrm{r}(B_1)},
$$
so
$$
\widehat{h}(A,B)=\wth\left(\frac A{\mathrm{r}(A)},\,\frac
B{\mathrm{r}(B)}\right)=\wth\left(\frac
{A_1}{\mathrm{r}(A_1)},\,\frac
{B_1}{\mathrm{r}(B_1)}\right)=\widehat{h}(A_1,B_1).
$$
This implies that the above definition of
$\widehat{h}(\widehat{A},\widehat{B})$ is well-defined.

\paragraph{{\bf 2.6.6.}} 
Denote
$$
\widehat{\mathcal{C}}(\mathbb{R}^n)=\left\{\widehat{C}:\: C
\in\mathcal{C}(\mathbb{R}^n)\right\}.
$$

\paragraph{{\bf 2.7. Theorem.}} 
 {\it
$\bigl(\widehat{\mathcal{C}}(\mathbb{R}^n),\widehat{h}\bigr)$ is a
complete metric space. }
\par\noindent{\it Proof.}   Let $A$, $B$,
$C\in\mathcal{C}(\mathbb{R}^n)$. It is obvious that
$$
\widehat{h}(\widehat{A},\widehat{B})=\widehat{h}(\widehat{B},\widehat{A})\geqslant0.
$$
By Lemma~2.6.4 we see that $\widehat{A}=\widehat{B}$ if and only
if
$$
 \frac A{\mathrm{r}(A)}\sim \frac
B{\mathrm{r}(B)},
$$
 which is equivalent to
$$
\wth\left(\frac A{\mathrm{r}(A)},\,\frac
B{\mathrm{r}(B)}\right)=0,
$$
 i.e., $\widehat{h}(\widehat{A},\widehat{B})=0$.

According to Definition~2.6.5 and Theorem~2.2 we have
\begin{align*}
\widehat{h}(\widehat{A},\widehat{B})&=\wth\left(\frac
A{\mathrm{r}(A)},\,\frac B{\mathrm{r}(B)}\right)\\
 &\leqslant
\wth\left(\frac A{\mathrm{r}(A)},\,\frac C{\mathrm{r}(C)}\right)+
\wth\left(\frac C{\mathrm{r}(C)},\,\frac B{\mathrm{r}(B)}\right)
=\widehat{h}(\widehat{A},\widehat{C})+\widehat{h}(\widehat{C},\widehat{B}).
\end{align*}

 Now suppose $\{\widehat{A}_k\}_{k=1}^\infty$ is a Cauchy sequence
 in $\bigl(\widehat{\mathcal{C}}(\mathbb{R}^n),\widehat{h}\bigr)$.
 Then $\{\widetilde{B}_k\}_{k=1}^\infty$, where $\displaystyle B_k=\frac {A_k}{\mathrm{r}(A_k)}$,
  is a Cauchy sequence
 in $\bigl(\widetilde{\mathcal{C}}(\mathbb{R}^n),\wth\bigr)$ by
 Definition~2.6.5, which implies by Theorem~2.2 that there exists
 $A\in\mathcal{C}(\mathbb{R}^n)$ such that
 $$
 \widetilde{h}(\widetilde{B}_k,\widetilde{A})\to0\quad(k\to+\infty).
 $$
For any $\varepsilon>0$, there exists $K>0$ such that when $k>K$
we have
$$
\widetilde{h}(\widetilde{B}_k,\widetilde{A})<\varepsilon.
 $$
 By Proposition~2.6.3(4) it follows that
 $$
 1-\varepsilon<\mathrm{r}(A)<1+\varepsilon.
 $$
 Thus $\mathrm{r}(A)=1$. Consequently
\begin{equation}
\widehat{h}(\widehat{A}_k,\widehat{A})=\widetilde{h}(\widetilde{B}_k,\widetilde{A})\to0\quad(k\to+\infty).
\tag*{$\Box$}
\end{equation}


\paragraph{{\bf 2.8. Relative rigid shape differences.}}

By a normal reasoning we know that an $r$-similitude $S:\: \mathbb
R^n\to\mathbb R^n$ can just be expressed to be
$S=\mu_r\circ\sigma\circ\tau_{_{\scriptstyle a}}$, where $\mu_r$
is an $r$-homothety, $\sigma$ is an orthogonal transformation and
$\tau_{_{\scriptstyle a}}$ is a translation (see
\cite[Proposition~2.3(1)]{Hutchinson-1981}). Now let us give the
following

\paragraph{{\bf 2.8.1. Definition.}}
A similitude $S:\: \mathbb R^n\to\mathbb R^n$ is called a {\it
rigid $r$-similitude}\/ ($r>0$) if
$S=\mu_r\circ\sigma\circ\tau_{_{\scriptstyle a}}$, where $\sigma$
is a rotation, i.e., $\det\sigma=1$ (cf. Definition~2.4.1).

\paragraph{{\bf 2.8.2. Lemma.}}
 {\it Let\/ $r$, $r_{_1}$, $r_{_2}$ be positive real numbers.
\begin{itemize}
\item[{\rm(1)}]  Let\/ $\mu_r$ be an\/ $r$-homothety and
let\/ $\sigma$, $\sigma_{_{\scriptstyle 1}}$,
$\sigma_{_{\scriptstyle 2}}$ be orthogonal transformations. Let\/
$\tau_{_{\scriptstyle a}}$ $(a\in\mathbb R^n)$ be a translation
defined in Lemma~2.4.3. Then
\begin{enumerate}
\item [{\rm(i)}] $\mu_r\circ\sigma=\sigma\circ\mu_r$;
\item [{\rm(ii)}] $\mu_r\circ\tau_{_{\scriptstyle a}}=\tau_{_{\scriptstyle ra}}\circ\mu_r$;
\item [{\rm(iii)}] $\mu_{r_{_1}}\circ\mu_{r_{_2}}=\mu_{{r_{_1}}{r_{_2}}}$;
\item [{\rm(iv)}] $\mu_r^{-1}=\mu_{r^{-1}}$;
\item [{\rm(v)}] $\sigma_{_{\scriptstyle 1}}\circ\sigma_{_{\scriptstyle 2}}$
    is an orthogonal transformation, and\/ $\det\sigma_{_{\scriptstyle 1}}=\det\sigma_{_{\scriptstyle 2}}=1$
    implies\/ $\det(\sigma_{_{\scriptstyle 1}}\circ\sigma_{_{\scriptstyle 2}})=1$;
\item [{\rm(vi)}] $\sigma^{-1}$ is an orthogonal transformation, and\/
     $\det\sigma=1$ if and only if\/ $\det\sigma^{-1}=1$.
\end{enumerate}

\item[{\rm(2)}]  Let $S_i$ be a rigid $r_{_{\scriptstyle i}}$-similitude
in\/ $\mathbb R^n$ $(i=1$, $2)$. Then
\begin{enumerate}
\item [{\rm(i)}] identical mapping\/ $\mathrm{id}:\: \mathbb
R^n\to\mathbb R^n$
 is a rigid\/ $1$-similitude;
\item [{\rm(ii)}] $S_1^{-1}$ is a rigid $r_{_1}^{-1}$-similitude;
\item [{\rm(iii)}] $S_1\circ S_2$ is a rigid $r_{_1}r_{_2}$-similitude.
\end{enumerate}
\end{itemize}
 }
\par\noindent{\it Proof.}   (1) and (2)(i) are clear. Suppose
$S_1=\mu_{r_{_1}}\circ\sigma_{_{\scriptstyle
1}}\circ\tau_{_{\scriptstyle a_{_1}}}$ and
$S_2=\mu_{r_{_2}}\circ\sigma_{_{\scriptstyle
2}}\circ\tau_{_{\scriptstyle a_{_2}}}$. Then by (1) and
Lemma~2.4.3,
\begin{align*}
 S_1^{-1}&=\tau_{_{\scriptstyle a_{_1}}}^{-1}\circ\sigma_{_{\scriptstyle 1}}^{-1}\circ\mu_{r_{_1}}^{-1}
 =\tau_{_{\scriptstyle -a_{_1}}}\circ\sigma_{_{\scriptstyle 1}}^{-1}\circ\mu_{r_{_1}^{-1}}\\
 &=\tau_{_{\scriptstyle -a_{_1}}}\circ\mu_{r_{_1}^{-1}}\circ\sigma_{_{\scriptstyle 1}}^{-1}
 =\mu_{r_{_1}^{-1}}\circ\tau_{_{\scriptstyle -r_{_1}a_{_1}}}\circ\sigma_{_{\scriptstyle 1}}^{-1}\\
 &=\mu_{r_{_1}^{-1}}\circ\sigma_{_{\scriptstyle 1}}^{-1}\circ\tau_{_{\scriptstyle -r_{_1}\sigma_{_1}(a_{_1})}}
\end{align*}
is a rigid $r_{_1}^{-1}$-similitude; and
\begin{align*}
 S_1\circ S_2&=\mu_{r_{_1}}\circ\sigma_{_{\scriptstyle 1}}\circ\tau_{_{\scriptstyle a_{_1}}}
 \circ\mu_{r_{_2}}\circ\sigma_{_{\scriptstyle 2}}\circ\tau_{_{\scriptstyle a_{_2}}}
\\
 &=\mu_{r_{_1}}\circ\sigma_{_{\scriptstyle 1}}\circ\mu_{r_{_2}}\circ\sigma_{_{\scriptstyle 2}}
 \circ\tau_{_{\scriptstyle r_{_2}^{-1}\sigma_{_2}^{-1}(a_{_1})}}\circ\tau_{_{\scriptstyle a_{_2}}}
\\
 &=\mu_{r_{_1}r_{_2}}\circ(\sigma_{_{\scriptstyle 1}}\circ\sigma_{_{\scriptstyle 2}})
 \circ\tau_{_{\scriptstyle a_{_2}+r_{_2}^{-1}\sigma_{_2}^{-1}(a_{_1})}}
\end{align*}
is a rigid $r_{_1}r_{_2}$-similitude.\hfill $\Box$


\paragraph{{\bf 2.8.3. Definition.}}
Two nonempty subsets $A$ and~$B$ of~$\mathbb R^n$ are {\it rigidly
similar}\/ if there exists a rigid similitude~$S$ in~$\mathbb R^n$
such that $B=S(A)$, denoted by $A\check{\sim} B$.
 By Lemma~2.8.2(2) we know that the rigid similarity is an
equivalence relation. The equivalence class of~$A$ is
denoted~$\check{A}$. And we denote
$$
\check{\mathcal{C}}(\mathbb{R}^n)=\left\{\check{C}:\: C
\in\mathcal{C}(\mathbb{R}^n)\right\}.
$$

\bigskip
By Lemma~2.8.2 and similarly to Lemma~2.6.4 we may get
\paragraph{{\bf 2.8.4. Lemma.}} 
{\it
 Let $A$ and~$B$ be nonempty bounded subsets of\/~$\mathbb R^n$.
 Then\/ $\displaystyle \frac A{\mathrm{r}(A)}\simeq \frac B{\mathrm{r}(B)}$
  if and only if $\check{A}=\check{B}$.\hfill $\Box$
 }

\paragraph{{\bf 2.8.5. Definition.}} 
Let $A$ and~$B$ be nonempty bounded subsets of~$\mathbb R^n$.
Define
$$
\check{h}(A,B):=\overline{h}\left(\frac A{\mathrm{r}(A)},\,\frac
B{\mathrm{r}(B)}\right)
$$
and
$$
\check{h}(\check{A},\check{B}):=\check{h}(A,B),
$$
which is called the {\it relative rigid shape difference} between
$\check{A}$ and~$\check{B}$ or between $A$ and~$B$.

\paragraph{{\it Remark}.} By Lemma~2.8.4 the above definition of
$\check{h}(\check{A},\check{B})$ is well-defined.

\paragraph{{\bf 2.9. Theorem.}} 
 {\it
$\bigl(\check{\mathcal{C}}(\mathbb{R}^n),\check{h}\bigr)$ is a
complete metric space.
 }

\bigskip
The proof of this theorem is similar to that of Theorem~2.7.\hfill
$\Box$

\paragraph{{\bf 2.10.} {\it Remark}.}
 (1) The metric spaces
$\bigl(\overline{\mathcal{C}}(\mathbb{R}^n),\overline{h}\bigr)$,
$\bigl(\widehat{\mathcal{C}}(\mathbb{R}^n),\widehat{h}\bigr)$ and
$\bigl(\check{\mathcal{C}}(\mathbb{R}^n),\check{h}\bigr)$ are all
separable.

 (2) We may regard $\bigl({\mathcal{C}}(\mathbb{R}^n),\overline{h}\bigr)$,
 $\bigl({\mathcal{C}}(\mathbb{R}^n),\widehat{h}\bigr)$ and
  $\bigl({\mathcal{C}}(\mathbb{R}^n),\check{h}\bigr)$ as
$\bigl(\overline{\mathcal{C}}(\mathbb{R}^n),\overline{h}\bigr)$,
$\bigl(\widehat{\mathcal{C}}(\mathbb{R}^n),\widehat{h}\bigr)$ and
$\bigl(\check{\mathcal{C}}(\mathbb{R}^n),\check{h}\bigr)$
respectively;
 and regard
``$\simeq$'', ``$\widehat{\sim}$'' and ``$\check\sim$'' as ``$=$''
correspondingly in the definition of metric spaces. Under this
convention, we may say that
 $\bigl({\mathcal{C}}(\mathbb{R}^n),\overline{h}\bigr)$,
 $\bigl({\mathcal{C}}(\mathbb{R}^n),\widehat{h}\bigr)$ and
  $\bigl({\mathcal{C}}(\mathbb{R}^n),\check{h}\bigr)$ are complete
metric spaces.

 (3) In Definition~2.6.5 and~2.8.5 we may define the relative
shape difference and relative rigid shape difference using
diameters instead of radii, i.e., define
$$
\widehat{h}(A,B):=\wth\left(\frac A{|A|},\,\frac B{|B|}\right)
\quad\text{and}\quad \check{h}(A,B):=\overline{h}\left(\frac
A{|A|},\,\frac B{|B|}\right),
$$
and obtain similar results.

 (4) We may also define the translation shape difference by
using translation equivalence instead of isometric equivalence,
and obtain some similar results.

\section{Perturbation of self-similar sets}\nopagebreak

In 1946, P.~A.~P.~Moran considered a self-similar set as an
extension of Cantor's set, obtained its Hausdorff dimension and
proved it has a finite and positive Hausdorff measure at its
Hausdorff dimension when it satisfies the open set condition (see
\cite{Moran-1946}).

Researches into self-similar sets were once motivated by
Mandelbrot's work (see \cite{Mandelbrot-1977} and
\cite{Mandelbrot-1982}). In 1981, J.~Hutchinson
(\cite{Hutchinson-1981}) considered a self-similar set as an
invariant set of a finite set of contraction maps (similitudes)
(called an {\it iterated function system}\/) in a systematic
manner and a mathematical self-similar set is presented in a clear
and wonderful way. In \cite{Hutchinson-1981} he also considered an
invariant measure with respect to the iterated function system.
Later on a large number of researches have been done in the area
of self-similarity and related subjects.

Researches on self-similarities have been developed in many
directions. Separation properties for self-similar sets have been
considered in \cite{Schief-1994} and \cite{Zerner-1996}, etc. For
researches of iterated function systems and some related topics,
we refer to \cite{Bant-1989}, \cite{Barnsley-1993},
\cite{Barnsley-2006}, \cite{Barnsley-Demko-1985},
\cite{Barnsley-Hutchinson-Stenflo-2008}, \cite{Falconer-2003},
\cite{Hutchinson-1981} and \cite{Hutchinson-Ruschendorf-1998},
etc. If the iterated function system consists of affine
transformations then the invariant set is a self-affine set, which
is an extension of a self-similar set and has been investigated
extensively (see e.g. \cite{Barnsley-1993}, \cite{Bedford-1989},
\cite{Bedford-Urbanski-1990}, \cite{Falconer-1988},
\cite{Falconer-1992} and \cite{McMullen-1984}, etc). Infinite
iterated function systems have also been considered (see e.g.
\cite{Fernau-1994}, \cite{Mauldin-Urbanski-1996} and
\cite{Moran-1996}, etc). For researches on random cases we refer
to e.g. \cite{Barnsley-Hutchinson-Stenflo-2005},
\cite{Barnsley-Hutchinson-Stenflo-2008}, \cite{Falconer-1986},
\cite{Graf-1987}, \cite{Graf-Mauldin-Williams-1988},
\cite{Hutchinson-Ruschendorf-1998},
\cite{Hutchinson-Ruschendorf-2000}
and\cite{Mauldin-Williams-1986}, etc. For Moran sets, which are
extensions of self-similar sets, we refer to e.g.
\cite{Cawley-Mauldin-1992}, \cite{Feng-Wen-Wu-1997},
\cite{Hua-Rao-Wen-Wu-2000}, \cite{Li-Dekking-2000},
\cite{Liu-Wen-2005}, \cite[Chapter~8]{Wen-2000} and
\cite{Wen-2001}, etc. For graph directed constructions we refer to
e.g. \cite{Das-Ngai-2004} and \cite{Mauldin-Williams-1988}, etc.
For sub-self-similar sets we refer to e.g. \cite{Falconer-1995}
and \cite[Section~3.1]{Falconer-1997}, etc. We note that Dekking
(\cite{Dekking-1982}) once gave a recurrent structure method to
construct some fractals. A kind of quasi-self-similar sets has
been considered in e.g. \cite[Theorem~8.6]{Blanchard-1984},
\cite{Falconer-1989}, \cite{McLaughlin-1987} and
\cite[p.\nolinebreak\,742]{Sullivan-1983}, etc.

$\cdots\cdots$

In this section we consider a kind of fractals which can be
approximately regarded as self-similar sets and we deal with them
as perturbation of strict self-similar sets.

\paragraph{{\bf 3.1. Sequences of integers.}}

We call $\bi=\bi_k=i_1\cdots i_k$, where $i_j\in\mathbb P$ (the
set of positive integers) ($j=1$, $\dots$, $k$), a {\it word of
finite length}\/ (the {\it length}\/~$|\bi|=k$); and call
$\ba=i_1\cdots i_j\cdots$, where $i_j\in\mathbb P$ ($j=1$, $2$,
$\dots$), a {\it word of infinite length}.

Let $\{m_j\}_{j=1}^\infty$ be a sequence of positive integers and
usually $m_j\geqslant2$ ($j=1$, $2$, $\dots$). Let $\bi=i_1\cdots
i_k$. We write $\bj=\bi\, i_{k+1}\cdots i_l$ if $\bj=i_1\cdots i_k
i_{k+1}\cdots i_l$ ($l\geqslant k$) and write $\ba=\bi\,
i_{k+1}\cdots i_l\cdots$ if $\ba=i_1\cdots i_k i_{k+1}\cdots
i_l\cdots$. And then we denote $\bj|_k=\bi$ and $\ba|_k=\bi$.

Now let us define
$$
\hat\bi:=\{\bi\, i_{k+1}\cdots i_l\cdots:\: \text{$i_l=1$,
$\dots$, $m_l$; \, $l=k+1$, $k+2$, $\dots$}\},
$$
where $|\bi|=k\geqslant1$, and we still denote $\hat\bi$ by~$\bi$.
Define
\begin{align*}
\mathcal I_\infty &=\mathcal I_\infty(\{m_j\}) \\
&:=\{\ba=i_1\cdots i_k\cdots:\: \text{$i_k=1$, $\dots$, $m_k$; \,
$k=1$, $2$, $\dots$}\},
\end{align*}
 which is also denoted $\hat\bzero$
or~$\bzero$. Assume $\bi|_0:=\bzero$. Then $\bi\supseteq\bj$ if
and only if $|\bj|=l\geqslant k=|\bi|$ and $\bj|_k=\bi$
($k\geqslant 0$). It is obvious that $\ba\in\bi$ if $\ba|_k=\bi$
($k\geqslant 1$) and that all $\ba\in\bzero$. Let
\begin{align*}
\mathcal I&=\mathcal I(\{m_j\}) \\
&:=\{\bi=i_1\cdots i_k:\: \text{$i_k=1$, $\dots$, $m_k$; \, $k=1$,
$2$, $\dots$}\}\cup\{\bzero\},
\end{align*}
$$
\mathcal I_k=\mathcal I_k(\{m_j\}):=\{\bi\in\mathcal I:\: |\bi|=k
\} \quad (k=0,\,1,\,2,\, \dots)
$$
 and
$$
\mathcal I^{l)}=\mathcal I^{l)}(\{m_j\}):=\bigcup_{k=0}^l\mathcal
I_k(\{m_j\}) \quad (l=0,\,1,\,2,\, \dots).
$$
 If $m_k=m$ for all $k\in\mathbb P$
then we say $\bi$ and~$\ba$ to be {\it normal words} and,
$\mathcal I(\{m_j\})$, $\mathcal I_k(\{m_j\})$, $\mathcal
I^{l)}(\{m_j\})$ and~$\mathcal I_\infty(\{m_j\})$ are denoted
$\mathcal I(m)$, $\mathcal I_k(m)$, $\mathcal I^{l)}(m)$
and~$\mathcal I_\infty(m)$ respectively.

If some or all of ~$m_j$ ($j=1$, $2$, $\dots$) equal $+\infty$, we
may also give similar concepts to the above and we will use the
same notations to denote them.

\paragraph{{\bf 3.2. Perturbation of self-similar sets.}}

\paragraph{{\bf 3.2.1. Definition.}}
Let $\mathcal S=\{S_i:\:i=1,\,\dots,\,m\}$ be a family of
contraction similitudes, which is called an {\it iterated function
system of similitudes}\/ (abbreviated to {\it IFSS}\/). Let $E$ be
the compact invariant set determined by~$\mathcal S$, i.e.,
$\mathcal S(E)=E$, where $\mathcal S(E):=\bigcup_{i=1}^m S_i(E)$
(see \cite[Chapter~9]{Falconer-2003} and~\cite{Hutchinson-1981}).

Let $\mathcal F=\{F_\bi:\: \bi\in\mathcal I\}$ be a family of
compact sets in~$\mathbb R^n$ satisfying
$$
F_\bi=\bigcup_{i=1}^m F_{\bi i}
$$
for all $\bi\in\mathcal I$. Let $F=F_\bzero$. Then we say that
$\mathcal F$ is a {\it structure system}\/ of~$F$. Given a family
$\varDelta=\{\delta_\bi\geqslant0:\: \bi\in\mathcal I\}$ of
nonnegative real numbers, let $\delta=\sup\{\delta_\bi:\:
\delta_\bi\in\varDelta\}$. Assume $\mathrm{Lip}\, S_i=c_i$, i.e.,
$$
|S_i(x)-S_i(y)|=c_i|x-y|
$$
for $x$, $y\in\mathbb R^n$ ($i=1$, $\dots$, $m$). Denote
$c_\bi:=c_{i_1}\cdots c_{i_k}$, where $\bi=i_1\cdots i_k$, and
$c_\bzero:=1$.

Suppose
\begin{equation}
\wth(F_\bi,E_\bi)\leqslant\delta_\bi c_\bi\,
\mathrm{r}(E),\label{tag3.2-1}
\end{equation}
where $S_\bi:=S_{i_1}\circ\cdots\circ S_{i_k}$ and
$E_\bi:=S_\bi(E)$ for $\bi=i_1\cdots i_k$, $S_\bzero:=\mathrm{id}$
(the identical mapping), $\mathrm{r}(E)$ is the radius of~$E$.
Then $F$ is called a {\it $\varDelta$-perturbation of the
self-similar set}\/~$E$ or a {\it $\varDelta$-quasi-self-similar
set}\/ with the IFSS~$\mathcal S$. If $\delta<+\infty$, then $F$
is also called a {\it $\delta$-perturbation of the self-similar
set}\/~$E$ or a {\it $\delta$-quasi-self-similar set}\/ with the
IFSS~$\mathcal S$.

\paragraph{{\it Remark}.} The inequality~(\ref{tag3.2-1}) implies
\begin{equation*}
\wth(F_\bi,E_\bi)\leqslant\delta_\bi c_\bi\,
|E|,\tag{\ref{tag3.2-1}$'$}\label{tag3.2-1$'$}
\end{equation*}
where $|E|$ is the diameter of~$E$. In the definition we may use
(\ref{tag3.2-1$'$}) to replace~(\ref{tag3.2-1}) and obtain similar
results.

\paragraph{{\bf 3.2.2. Open set condition.}}
$F$ is said to satisfy the {\it open set condition}\/ if there
exists a family $\{V_\bi:\: \bi\in\mathcal I\}$ of open sets such
that
\begin{itemize}
\item[(1)] $F_\bi\subseteq\mathrm{cl}\,(V_\bi)$ for $\bi\in\mathcal I$;
\item[(2)] $V_\bi\cap V_\bj=\emptyset$ (the empty set) for $\bi$,
$\bj\in\mathcal I$ and $\bi\cap\bj=\emptyset$;
\item[(3)] there exist two positive constants $a_1$ and~$a_2$ so
that each~$V_\bi$ contains a ball of radius~$a_1c_\bi$ and is
contained in a ball of radius~$a_2c_\bi$ ($\bi\in\mathcal I$).
\end{itemize}

\paragraph{{\bf 3.2.3. Lemmas.}}
(1) {\it Let\/ $a_1$ and\/~$a_2$ be two positive constants and\/
$r>0$. Suppose\/ $\{V_\bi:\: \bi\in\mathcal I\}$ is a family of
disjoint open sets. If each\/~$V_\bi$ contains a ball of
radius\/~$a_1 r$ and is contained in a ball of radius\/~$a_2 r$
$(\bi\in\mathcal I)$, then any closed ball\/~$B$ of radius\/~$r$
meets at most\/ $(1+2a_2)^n a_1^{-n}$ of the
closures\/~$\mathrm{cl}\,(V_\bi)$ } (see
\cite[Lemma~9.2]{Falconer-2003} or
\cite[Lemma~5.3(a)]{Hutchinson-1981}).

(2) {\it Let\/ $s$ be the similarity dimension\/ $\dim_S E$
of\/~$\mathcal S$ or\/~$E$, i.e., a unique solution of equation
$$
\sum_{i=1}^m c_i^s=1.
$$
Define
$$
\hat{\mu}(\bi)=c_\bi^s
$$
for\/ $\bi\in\mathcal I$. Then\/ $\hat{\mu}(\bi)$ can be expanded
into a measure or a mass distribution on\/~$\mathcal I_\infty$
with\/ $\hat{\mu}(\mathcal I_\infty)=1$}.

{\it For\/ $A\subseteq\mathbb R^n$, let
$$
I_A:=\{\ba\in\mathcal I_\infty:\: x_\sba\in A\cap F\},
$$
where\/ $\{x_\sba\}:=\bigcap _{\bi\ni\sba} F_\bi$, and
$$
\mu(A):=\hat{\mu}(I_A).
$$
Then\/ $\mu$ is a\/ $($an outer\/$)$ measure on\/~$F$, i.e.,
\begin{itemize}
\item[{\rm(i)}] $\mu(\emptyset)=0$ and\/ $\mu(A)\geqslant0$ for\/ $A\subseteq\mathbb R^n$;
\item[{\rm(ii)}] $\mu(A)\leqslant\mu(B)$ if\/ $A\subseteq B$;
\item[{\rm(iii)}] If\/ $A=\bigcup_{i=1}^{+\infty} A_i$ then
$$
\mu(A)\leqslant\sum_{i=1}^{+\infty} A_i.
$$
\end{itemize}
}
\noindent(see
 \cite[Section~1.3]{Falconer-1997}, \cite[the proof of
Theorem~9.3]{Falconer-2003} and \cite[\S3.2]{Wen-2000}).

(3) {\bf Mass distribution principle} (see
\cite[Section~4.1]{Falconer-2003} and
\cite[Theorem~I]{Moran-1946}).
 {\it
  Suppose that $\mu$ is a mass distribution on\/~$F$ $($a measure
  on~\/$F$ satisfying $0<\mu(F)<+\infty)$ and for some positive
  constants~$s$, $c$ and~$\varepsilon$ we have
  $$
  \mu(U)\leqslant c|U|^s
  $$
  for any subset~$U$ of~$F$ with $|U|\leqslant\varepsilon$. Then
  $$
  \mathcal H^s(F)\geqslant c^{-1}\mu(F)
  $$
  and
  $$
  s\leqslant\dim_H F\leqslant\underline{\dim}_B F\leqslant\overline{\dim}_B F,
  $$
where $\mathcal H^s(F)$ denotes $s$-dimensional Hausdorff measure
of~$F$ and, $\dim_H F$, $\underline{\dim}_B F$ and
$\overline{\dim}_B F$ denote the Hausdorff dimension, lower and
upper box dimensions of~$F$ respectively
 }\/ (see \cite[Chapters~3
and~4]{Falconer-2003}).\hfill $\Box$

\paragraph{{\bf 3.2.4. Theorem.}}
{\it
 If $F$ is a $\delta$-quasi-self-similar set with
IFSS~$\mathcal S$ satisfying the open set condition\/~$(3.2.2)$,
then
$$
\dim_H F=\dim_B F=s,
$$
where\/ $s$ is the similarity dimension of\/~$\mathcal S$
$($or\/~$F)$, and\/ $0<\mathcal H^s(F)<+\infty$.
 }
\par\noindent{\it Proof.}   We follow a normal method introduced by
P.~A.~P.~Moran (\cite{Moran-1946})(refer to
\cite[Section~9.2]{Falconer-2003} and
\cite[Section~5]{Hutchinson-1981}).

From~(\ref{tag3.2-1$'$}) it follows that
\begin{align*}
|F_{\bi}|&\leqslant|E_{\bi}|+2\delta_{\bi}c_{\bi}|E| \\
 &\leqslant (1+2\delta)c_{\bi}|E|\leqslant (1+2\delta)c_{\max}^k |E|
\end{align*}
by Proposition~2.6.3(2), where $c_{\max}=\max\{c_i:\:
i=1,\,\dots,\, m\}$. Hence given $\varepsilon>0$ there is a
$k\in\mathbb P$ such that $|F_{\bi}|\leqslant\varepsilon$ for all
$\bi\in\mathcal I_k$. Since
$$
F=\bigcup_{\bi\in\mathcal I_k}F_\bi,
$$
we have
\begin{align*}
\mathcal H_\varepsilon^s(F)&\leqslant\sum_{\bi\in\mathcal I_k}|F_{\bi}|^s \\
 &\leqslant (1+2\delta)^s|E|^s\sum_{\bi\in\mathcal I_k}
 c_{\bi}^s=(1+2\delta)^s|E|^s.
\end{align*}
Therefore $\mathcal H^s(F)\leqslant(1+2\delta)^s|E|^s$.

Now let $B=B(r)$ be a closed ball of radius~$r>0$. For any
$\ba\in\mathcal I_\infty$ choose the smallest~$k$ such that
$c_\bi\leqslant r$, where $\bi=\ba|_k$. Then $c_{\min} r<c_\bi$,
where $c_{\min}:=\min\{c_i:\: i=1,\,\dots,\,m\}$. Let $I(r)$
denote the set of all such~$\bi$. Then by the open set
condition~(3.2.2),
$$
F=\bigcup_{\bi\in I(r)} F_\bi\subseteq\bigcup_{\bi\in I(r)}
\mathrm{cl}\,(V_\bi),
$$
where each $V_\bi$ ($\bi\in I(r)$) contains a ball of radius~$a_1
c_{\min}r$ and is contained in a ball of radius~$a_2 r$. Let
$I^*(r)=\{\bi\in I(r):\: B\cap V_\bi\neq\emptyset\}$. Then
$$
\sharp\,(I^*(r))\leqslant a=(1+2a_2)^n a_1^{-n} c_{\min}^{-n}
$$
 by Lemma~3.2.3(1), where $\sharp\, I$ denotes the number of elements
in~$I$, and
$$
I_{B\cap F}\subseteq\bigcup_{\bi\in I^*(r)}\bi.
$$
Therefore
\begin{align*}
\mu(B)&=\mu(B\cap F)=\hat{\mu}(I_{B\cap
F})\leqslant\hat{\mu}\left(\bigcup_{\bi\in I^*(r)}\bi\right) \\
&\leqslant\sum_{\bi\in I^*(r)}\hat{\mu}(\bi)=\sum_{\bi\in
I^*(r)}c_\bi^s\leqslant a r^s.
\end{align*}
For a subset~$U$ of~$\mathbb R^n$ let $B=B(r)$ be a closed ball of
radius~$r=|U|$ centered at a point of~$U$. Then $U\subseteq B$,
and consequently
$$
\mu(U)\leqslant\mu(B)\leqslant a r^s=a|U|^s.
$$
By Lemma~3.2.3(3) we obtain
$$
\mathcal H^s(F)\geqslant a^{-1}\mu(F)=a^{-1}.
$$

Let $q(r)=\sharp\, I(r)$. Then
$$
q(r) c_{\min}^s r^s\leqslant\sum_{\bi\in I(r)} c_\bi^s=1.
$$
Thus $q(r)\leqslant c_{\min}^{-s} r^{-s}$.
From~(\ref{tag3.2-1$'$}) we deduce
\begin{align*}
|F_{\bi}|&\leqslant|E_{\bi}|+2\delta_{\bi}c_{\bi}|E| \\
 &\leqslant (1+2\delta_\bi)c_{\bi}|E|\leqslant (1+2\delta)|E|r=br,
\end{align*}
where $b=(1+2\delta)|E|$ is a positive constant. So finally
\begin{align*}
\overline{\dim}_B F&=\limsup_{r\to 0^+}\frac{\log
N(br)}{-\log(br)}\leqslant\limsup_{r\to 0^+}\frac{\log
q(br)}{-\log(br)} \\
&\leqslant\limsup_{r\to 0^+}\frac{\log
(c_{\min}^{-s}b^{-s}r^{-s})}{-\log(br)}=s,
\end{align*}
where $N(r)$ denotes the smallest number of sets of diameter at
most~$r$ which cover~$F$.
 \hfill $\Box$

\paragraph{{\bf 3.2.5. Corollary.}}
{\it
 Suppose\/ $F$ is a\/ $\delta$-quasi-self-similar set with
IFSS\/ $\mathcal S$ satisfying the following condition:
 There exists\/ $\varepsilon_0>0$ such that
\begin{equation}
d(F_\bi,F_\bj)\geqslant\varepsilon_0 c_{\bi|_k}\label{tag3.2-2}
\end{equation}
for any pair\/ $\bi$ and\/~$\bj$ satisfying\/
$\bi\cap\bj\neq\emptyset$, $\bi|_k=\bj|_k$ but\/ $i_{k+1}\neq
j_{k+1}$ $(k$ is some nonnegative integer associated with\/ $\bi$
and\/~$\bj)$. Then
$$
\dim_H F=\dim_B F=s,
$$
where\/ $s$ is the similarity dimension of\/~$\mathcal S$
$($or\/~$F)$, and\/ $0<\mathcal H^s(F)<+\infty$.
 }
\par\noindent{\it Proof.}   We only need to prove that
(\ref{tag3.2-1}) and~(\ref{tag3.2-2}) imply the open set
condition.

Let
$$
V_{\bi}=\mathscr{N}\left(F_\bi,\frac13\varepsilon_0 c_\bi \right).
$$
By~(\ref{tag3.2-1}) it follows that
$$
|F_\bi|\leqslant(1+2\delta)|E|c_\bi.
$$
Hence
$$
|V_\bi|\leqslant |F_\bi|+\frac23\varepsilon_0 c_\bi
\leqslant\left[(1+2\delta)|E|+\frac23\varepsilon_0\right]c_\bi.
$$
Let $a_1=\frac13{\varepsilon_0}$ and
$a_2=(1+2\delta)|E|+\frac23\varepsilon_0$. Then
$$
B_1(a_1 c_\bi)\subseteq V_\bi\subseteq B_2(a_2 c_\bi),
$$
where $B_1(a_1 c_\bi)$ and $B_2(a_2 c_\bi)$ are balls of radii
$a_1 c_\bi$ and~$a_2 c_\bi$ respectively. Denote $\bk=\bi|_k$
in~(\ref{tag3.2-2}). Then
\begin{align*}
d(V_\bi,V_\bj)&\geqslant d(F_\bi,F_\bj)-\frac13\varepsilon_0
c_\bi-\frac13\varepsilon_0 c_\bj \\
&\geqslant\varepsilon_0 c_\bk-\frac13\varepsilon_0
c_\bk-\frac13\varepsilon_0 c_\bk=\frac13\varepsilon_0 c_\bk>0.
\end{align*}
Hence $V_\bi\cap V_\bj=\emptyset$.
 \hfill $\Box$

\paragraph{{\bf 3.2.6.} {\it Remark}.}
We may follow an ordinary way below to construct a compact set
(fractal) in~$\mathbb R^n$ and its structure system.

Suppose $\{G_\bi:\: \bi\in\mathcal I\}$ is a family of nonempty
compact sets in~$\mathbb R^n$ such that $G_\bi\subseteq G_\bj$ if
$\bi\subseteq\bj$ and $|G_\bi|\to 0$ ($|\bi|\to+\infty$). Then for
$\ba\in\mathcal I_\infty$ the set $\bigcap_{\sba\in\bi}G_\bi$ is a
singleton, whose member is denoted~$x_\sba$. Let
$$
G^{(p)}(\bi):=\bigcup_{\bj\subseteq\bi,\; |\bj|=p} G_\bj
$$
 for
$p\geqslant |\bi|$ and
$$
F_\bi:=\bigcap_{p=|\bi|}^\infty G^{(p)}(\bi).
$$
Then $G^{(p)}(\bi)$ ($p\geqslant|\bi|$) are also nonempty compact
sets and $G^{(p)}(\bi)\subseteq G^{(q)}(\bi)$ if $p\geqslant q$,
hence $F_\bi$ is a nonempty compact set. Now we have constructed a
nonempty compact set~$F=F_\bzero$.

\paragraph{{\bf Proposition.}}
{\it
 $F_\bi=\mathrm{cl}\,(X_\bi)$, where
$X_\bi:=\{x_\sba:\:\ba\in\bi\}$.
 }
\par\noindent{\it Proof.}   Obviously $X_\bi\subseteq F_\bi$, thus
$\mathrm{cl}\,(X_\bi)\subseteq F_\bi$.

If $x\in F_\bi$ then $x\in G^{(p)}(\bi)$ for all
$p\geqslant|\bi|$. Hence for each $p\geqslant|\bi|$ there exists
$\bj\subseteq\bi$ so that $|\bj|=p$ and $x\in G_\bj$. Therefore
for any given neighborhood~$N(x)$ of~$x$ we can find $\bj$
($|\bj|\geqslant|\bi|$) such that $G_\bj\subseteq N(x)$. Now we
have $x_\sbb\in N(x)$ for $\bb\in\bj\subseteq\bi$.
 \hfill $\Box$

\paragraph{{\bf Corollary.}}
{\it
 $\displaystyle F_\bi=\bigcup_{i=1}^m F_{\bi i}.$
 }
\par\noindent{\it Proof.}   It is easy to see that
\begin{equation*}
X_\bi=\bigcup_{i=1}^m X_{\bi i}. \tag*{$\Box$}
\end{equation*}

By the corollary above we know that $\mathcal F=\{F_\bi:\:
\bi\in\mathcal I\}$ is a structure system of~$F$.

If $\{H_\bi:\: \bi\in\mathcal I\}$ is a family of nonempty compact
sets then we can let $G_{\bi}:=\bigcup_{\bj\subseteq\bi} H_\bj$.
If moreover $|G_\bi|\to0$ ($|\bi|\to+\infty$) then we return to
the above steps to construct a nonempty compact set (fractal) and
its structure system.

One perhaps more useful way to perturb a self-similar set (and
other similar structures) will be introduced in some concrete
examples below.

\paragraph{{\bf 3.2.7.} {\it Examples}.}
\paragraph{{\bf (1) Perturbation of the Cantor set.}}
Let $\mathscr C$ denote Cantor's ternary set, which is the
invariant set of $S_1(x)=\frac13x$ and $S_2(x)=\frac13x+\frac23$
in~$\mathbb R$. Let $\mathscr C_\bi:=S_\bi(\mathscr C)$.

(i) Given two families $\{a_\bi:\:\bi\in\mathcal I(2)\}$ and
$\{b_\bi:\:\bi\in\mathcal I(2)\}$ of real numbers satisfying
\begin{equation}
a_0\leqslant a_\bi<b_\bi\leqslant b_0,\label{tag3.2-3}
\end{equation}
 where $a_0$ and~$b_0$ are two fixed real numbers such that
 $a_0<b_0$, we assume $H_\bi=S_\bi([a_\bi,b_\bi])$ and
 $G_\bi=\bigcup_{\bj\subseteq\bi}H_\bj$.

Following Remark~3.2.6 we obtain a nonempty compact
set~$F=F_\bzero$ and its structure system  $\mathcal F=\{F_\bi:\:
\bi\in\mathcal I(2)\}$. It is easy to see that
$$
\wth(F_\bi,\mathscr C_\bi)\leqslant \frac{\delta_0}{3^k}\,\frac12,
$$
where $\delta_0=\max\{|b_0-a_0-1|,\,1\}$ and $k=|\bi|$. Thus $F$
is a $\delta_0$-perturbation of~$\mathscr C$, denoted $\mathscr C
\{a_\bi,b_\bi\}$. Specially we choose $a_\bi=0$ and $\mathscr C
\{a_\bi,b_\bi\}$ is written by~$\mathscr C\{b_\bi\}$. Then we may
let $a_0=0$. Assume $b_0<2$. We deduce
$$
d(G_\bi,G_\bj)\geqslant\frac1{3^k}
d(S_\bi([0,b_0]),\,S_\bj([0,b_0]))\geqslant\frac{\varepsilon_0}{3^k}
$$
for $\bi$, $\bj\in\mathcal I(2)$ and $\bi\cap\bj=\emptyset$, where
$\varepsilon_0=\min\left\{\frac23-\frac{b_0}3,\,\frac29\right\}>0$,
$\bi|_k=\bj|_k$ but $i_{k+1}\neq j_{k+1}$. By Corollary~3.2.5 we
have
$$
\dim_H\mathscr C\{b_\bi\}=\dim_B\mathscr C\{b_\bi\}=s,
$$
where $s=\log2/\log3$, and $0<\mathcal H^s(\mathscr
C\{b_\bi\})<+\infty$ ($0<b_0<2$).

We note that $\mathscr C \{a_\bi,b_\bi\}$ is also a Moran set
(see~\cite{Wen-2001}).

(ii) Now let us consider the perturbation of~$\mathscr C$
in~$\mathbb R^2$.

If we translate, rotate, stretch or contract each~$\mathscr C_\bi$
($\bi\in\mathcal I(2)$) on a reasonable small scale in~$\mathbb
R^2$ we may get a quasi-self-similar set. For example, let us
rotate each~$\mathscr C_\bi$ ($\bi\in\mathcal I(2)$) to perturb
$\mathscr C$ in~$\mathbb R^2$.

Let $L_\bi=S_\bi([0,1])\times\{0\}$. At first, we rotate
each~$L_i$ ($i=1$, $2$) around some point of itself to get a line
segment, denoted~$C_i$ ($i=1$, $2$). Thus each~$L_\bi$
($\bi\in\mathcal I(2)$, $|\bi|\geqslant2$) is moved to another
position. Let $L'_\bi$ denote $L_\bi$ in the new position
($\bi\in\mathcal I(2)$). Note $L'_\bi=C_\bi$ ($|\bi|=1$). Then we
rotate each~$L'_\bi$ ($|\bi|=2$) around some point of itself to
get a line segment~$C_\bi$ ($|\bi|=2$). In this way we may get
all~$C_\bi$ ($\bi\in\mathcal I(2)$). Let
$$
C_\bi^{(p)}=\bigcup_{\bj\subseteq\bi,\;|\bj|=|\bi|+p} C_\bj
\quad(p=0,\,1,\,2,\,\dots).
$$
Then
$$
C_\bi^{(p+1)}=\bigcup_{i=1}^2 C_{\bi i}^{(p)},
$$
\begin{align}
h(C_\bi^{(p)},C_\bi^{(q)})&=h\left(\bigcup_{\bj\subseteq\bi,\;|\bj|=|\bi|+p}
C_\bj,\,\bigcup_{\bj\subseteq\bi,\;|\bj|=|\bi|+p}
C_\bj^{(q-p)}\right)\notag \\
&\leqslant\sup
\left\{h(C_\bj,C_\bj^{(q-p)}):\:{\bj\subseteq\bi,\,|\bj|=|\bi|+p}\right\}\notag \\
 &\leqslant\frac12\,\frac1{3^{k+p}} \label{tag3.2-4}
\end{align}
if $q=p$, $p+1$, $\dots$, and
\begin{equation}
\wth(C_\bi^{(p)},\mathscr
C_\bi)\leqslant\frac12\,\frac1{3^k},\label{tag3.2-5}
\end{equation}
where $k=|\bi|$. By the completeness
of~$\bigl(\mathbb{R}^2,h\bigr)$ it follows that $C_\bi^{(p)}$
approaches a nonempty compact set~$F_\bi$ ($\bi\in\mathcal I(2)$)
in~$\mathbb R^2$ as $p\to+\infty$. By~(\ref{tag3.2-4}) we have
$$
h(C_\bi^{(p)},F_\bi)\leqslant\frac12\,\frac1{3^{k+p}}.
$$
Hence
\begin{align*}
&h(F_\bi,F_{\bi1}\cup F_{\bi2})\leqslant
h(F_\bi,C_\bi^{(p+1)})+h\left(\bigcup_{i=1}^2 C_{\bi
i}^{(p)},\,\bigcup_{i=1}^2 F_{\bi i}\right) \\
&\leqslant h(C_\bi^{(p+1)},F_\bi)+\max\left\{h(C_{\bi i}^{(p)},\,
F_{\bi i}):\:i=1,\,2\right\} 
\to0\quad(p\to+\infty).
\end{align*}
Therefore
$$
F_\bi=F_{\bi1}\cup F_{\bi2}
$$
for $\bi\in\mathcal I(2)$. From~(\ref{tag3.2-5}) we have
$$
\wth(F_\bi,\mathscr
C_\bi)\leqslant\frac12\,\frac1{3^k}\quad(k=|\bi|).
$$
Thus $\mathscr C'=\mathscr C'_2:=F_\bzero$ is a $1$-perturbation
of~$\mathscr C$ with a structure system $\mathcal
F=\{F_\bi:\:\bi\in\mathcal I(2)\}$.

Suppose that in the above procedure of constructing the structure
system~$\mathcal F$, the centers of rotation are restricted to
parts of~$C_\bi$ that are line segments of radii~$r_\bi\leqslant
r_0|C_\bi|$ ($r_0$ is some fixed positive number less
than~$\frac16$) centered at centers of~$C_\bi$ ($\bi\in\mathcal
I$) or absolute values~$\theta$ of angles of rotation of~$C_\bi$
are less than or equal to $\theta_0=\frac{\pi}2$ (which can
actually be replaced by some larger $\theta_0<\frac23\pi$). Then
we may deduce
$$
d(F_\bi,F_\bj)\geqslant\frac{\varepsilon_0}{3^k},
$$
where
$\varepsilon_0=\min\left\{\frac{1-6r_0}3,\,\frac1{12}\right\}$
 (for
$\theta_0=\frac{\pi}2$) and $k$ is a nonnegative integer such that
$\bi|_k=\bj|_k$ but $i_{k+1}\neq j_{k+1}$. By Corollary~3.2.5 we
know that
\begin{equation}
\dim_H\mathscr C'=\dim_B\mathscr
C'=\frac{\log2}{\log3}\label{tag3.2-6}
\end{equation}
and
\begin{equation}
0<\mathcal H^s(\mathscr C')<+\infty,\label{tag3.2-7}
\end{equation}
where $s=\log2/\log3$. In fact if we do not impose the above extra
limitations on the rotation of~$C_\bi$ ($\bi\in\mathcal I(2)$),
(\ref{tag3.2-6}) and~(\ref{tag3.2-7}) may still be true often.

The preceding discussion may similarly be conducted in~$\mathbb
R^3$ or even in~$\mathbb R^n$ to get compact sets $\mathscr C'_3$
or~$\mathscr C'_n$.

\paragraph{{\bf (2) Perturbation of the von Koch curve.}}
Let $\mathscr K$ denote the von Koch curve
(see~\cite{Falconer-2003} and~\cite{Hutchinson-1981}), which is
the invariant set of $\mathcal S=\{S_1,S_2,S_3,S_4\}$ in~$\mathbb
R^2$, where  $S_1(z)=\frac13z$,
$S_2(z)=\frac13\me^{\frac\pi3\mi}z+\frac13$,
$S_3(z)=\frac13\me^{-\frac\pi3\mi}z+\frac12+\frac{\sqrt3}6\mi$ and
$S_4(z)=\frac13z+\frac23$ ($z\in\mathbb C=\mathbb R^2$, $\mathbb
C$ denotes the set of all complex numbers). Let $\mathscr
K_\bi:=S_\bi(\mathscr K)$.

Let $\hat{\mathscr K}$ denote a random von Koch curve
in~\cite[Chapter~15]{Falconer-2003}. Let $K_0=[0,1]\times\{0\}$,
$\hat{S}_2(z)=\frac13\me^{-\frac\pi3\mi}z+\frac13$ and
$\hat{S}_3(z)=\frac13\me^{\frac\pi3\mi}z+\frac12-\frac{\sqrt3}6\mi$.
Suppose $T_i=S_i$ for $i=1$, $4$ and $T_i=S_i$ for $i=2$, $3$ or
$T_i=\hat{S}_i$ for $i=2$, $3$. Let
$$
K^{(k)}=\bigcup_{\bi\in\mathcal I_k(4)} S_\bi(K_0)
$$
and
$$
\hat{K}^{(k)}=\bigcup_{\bi\in\mathcal I_k(4)} T_\bi(K_0),
$$
where we choose corresponding~$T_\bi$, which may be different in
each step to construct~$T_\bi$. Then $K^{(k)}\to\mathscr K$ and
$\hat{K}^{(k)}\to\hat{\mathscr K}$ ($k\to+\infty$) in the sense of
Haudorff metric~$h$. At first let us note that $\hat{\mathscr K}$
is a $\frac19$-perturbation of~$\mathscr K$.

Let $V$ be the interior of the rhombus whose vertices are $0$,
$1$, $\frac12\pm\frac{\sqrt3}6\mi$. It is easy to see that
$\hat{\mathscr K}$ satisfy the open set condition~3.2.2 for
$V_\bi=S_\bi(V)$. Therefore
\begin{equation*}
\dim_H\hat{\mathscr K}=\dim_B\hat{\mathscr K}=\frac{\log4}{\log3}
\end{equation*}
and
\begin{equation*}
0<\mathcal H^s(\hat{\mathscr K})<+\infty,
\end{equation*}
where $s=\log4/\log3$ (cf. \cite[Exercise~15.3]{Falconer-2003}).

If $T_i=S_i$ for $i=1$, $4$ and $T_i=S_i$ or~$\hat{S}_i$ for
$i=2$, $3$, then we obtain another random von Koch
``curve''~$\tilde{\mathscr K}$ (the limit of
$\tilde{K}^{(k)}=\bigcup_{\bi\in\mathcal I_k(4)} T_\bi(K_0)$ in
the Hausdorff metric) (we call it a {\it random von Koch set}\/),
which is a $\frac13$-perturbation of~$\mathscr K$ satisfy the open
set condition~3.2.2.

Generally let
$$
T_{\theta,t}(z):=\frac13\me^{\mi\theta}z+t\quad(\theta\in\mathbb
R,\,t\in\mathbb C).
$$
Suppose $T_\bi:=T_{\theta_\bi,t_\bi}$, where $\theta_\bi\in\mathbb
R$, $|t_\bi|\leqslant\lambda$ ($\lambda$ is a fixed nonnegative
real number), and $K_\bi:=T_{\bi|_1}\circ
T_{\bi|_2}\circ\cdots\circ T_{\bi|_k}(L_0)$
($L_0=[0,1]\times\{0\}$ and $k=|\bi|$). Let
$$
K_\bi^{(p)}=\bigcup_{\bj\subseteq\bi,\;|\bj|=|\bi|+p} K_\bj
\quad(p=0,\,1,\,2,\,\dots).
$$
Then similarly to~(\ref{tag3.2-4}) it follows
$$
h(K_\bi^{(p)},K_\bi^{(q)})\leqslant\frac{\delta_1(\lambda)}{3^{k+p}}
\quad(q\geqslant p),
$$
where $\delta_1(\lambda)$ is a nonnegative real number related
to~$\lambda$. Let $K_\bi^{(p)}\to F_\bi$ ($p\to+\infty$) for
$\bi\in\mathcal I(4)$. We have
$$
F_\bi=\bigcup_{i=1}^4 F_{\bi i}.
$$
Similarly to~(\ref{tag3.2-5}) we can also deduce
$$
\wth(K_\bi^{(p)},\mathscr
K_\bi)\leqslant\frac{\delta_0(\lambda)}{3^k}\,\frac1{\sqrt3},
$$
where $\delta_0(\lambda)$ is a nonnegative real number related
to~$\lambda$. Letting $p\to+\infty$ we get
$$
\wth(F_\bi,\mathscr
K_\bi)\leqslant\frac{\delta_0(\lambda)}{3^k}\,\mathrm{r}(\mathscr
K).
$$
Thus $\mathscr K(\{\theta_\bi\},\{t_\bi\}):=F_\bzero$ (
$\theta_\bi\in\mathbb R$ and $|t_\bi|\leqslant\lambda$) is a
$\delta_0(\lambda)$-perturbation of~$\mathscr K$ with a structure
system $\mathcal F=\{F_\bi:\:\bi\in\mathcal I(4)\}$.

Now suppose $T_\bi:=T_{\theta_\bi,t_\bi,\rho_\bi}$, where
$\theta_\bi\in\mathbb R$, $|t_\bi|\leqslant\lambda$, $\rho_\bi>0$
and
$$
T_{\theta,t,\rho}(z):=\frac\rho3\me^{\mi\theta}z+t\quad(\theta\in\mathbb
R,\,t\in\mathbb C,\,\rho>0).
$$
Then following the preceding procedure we still get a
$\delta$-perturbation
 $$
 \mathscr K(\{\theta_\bi\},\{t_\bi\},\{\rho_\bi\})
 $$
 of $\mathscr K$, if we
put some suitable restrictive condition $\mathcal R$ on
$\{\rho_\bi\}$, e.g. we assume only finitely many $\rho_\bi\neq1$,
where $\delta=\delta(\lambda,\mathcal R)$ is a nonnegative real
number related to $\lambda$ and $\mathcal R$.

\paragraph{{\bf (3) Perturbation of the Sierpi\'{n}ski gasket.}}
Let $\mathscr S$ denote the Sierpi\'{n}ski gasket
(see~\cite{Falconer-2003}), which is the invariant set of
$\mathcal S=\{S_1,S_2,S_3\}$ in~$\mathbb R^2$, where
$S_1(z)=\frac12z+\frac14+\frac{\sqrt3}4\mi$, $S_2(z)=\frac12z$,
$S_3(z)=\frac12z+\frac12$ in $\mathbb C=\mathbb R^2$. Let
$\mathscr S_\bi:=S_\bi(\mathscr S)$ ($\bi\in\mathcal I(3)$).
Suppose
$A_\gamma=A\left(\gamma,1,\frac12+\frac{\sqrt3}2\mi\right)$
denotes the closed triangular region whose vertices are $\gamma$,
$1$ and $\frac12+\frac{\sqrt3}2\mi$ ($\gamma\in\mathbb C$).
Suppose $\Gamma=\{\gamma_\bi:\:\bi\in\mathcal I(4)\}$ satisfies
that $a\leqslant\gamma_\bi<1$ ($a$ is a fixed real number less
than~$1$) and that $\gamma_\bi\subseteq\gamma_\bj$ if
$\bi\supseteq\bj$. Let $G_\bi=S_\bi(A_{\gamma_\bi})$ and let
$\mathscr S(\Gamma)$ denote the nonempty compact set~$F_\bzero$
obtained by using the procedure in Remark~3.2.6. Then $\mathscr
S(\Gamma)$ is a $\delta$-perturbation of~$\mathscr S$. If $\Gamma$
satisfies some suitable condition, e.g.
$0\leqslant\gamma_\bi\leqslant b$, where $b$ is a fixed
nonnegative real number less than~$1$, then $\mathscr S(\Gamma)$
satisfies the open set condition~3.2.2 and thus $\mathscr
S(\Gamma)$ is an $s$-set, where $s=\log3/\log2$.

We may also perturb~$\mathscr S$ by following the way of
perturbing the von Koch curve~$\mathscr C$ above. Generally
suppose
$$
T_{r,\theta,t,\sigma,\rho}(z):=r\rho\,\me^{\mi\theta}\sigma(z)+t,
$$
where $0<r<1$, $0\leqslant\theta<2\pi$, $t\in\mathbb C$, $\rho>0$,
$\sigma(z)=z$ or $\sigma(z)=\overline{z}$. Let
$$
S_j=T_{r_j,\varphi_j,b_j,\sigma_j,1},
$$
where $0<r_j<1$, $0\leqslant\varphi_j<2\pi$, $b_j\in\mathbb C$,
$\sigma_j(z)=z$ or $\sigma_j(z)=\overline{z}$ ($j=1$, $\dots$,
$m$). Let $\mathfrak S=\mathfrak S(r_j,\varphi_j,b_j,\sigma_j)$
denote the invariant set of $\mathcal S=\{S_1,\dots,S_m\}$. Denote
$$
\mathcal
W=\{\{\theta_\bi\},\{t_\bi\},\{\sigma_\bi\},\{\rho_\bi\}\},
$$
where $0\leqslant\theta_\bi<2\pi$, $|t_\bi|\leqslant\lambda$,
$\sigma_\bi(z)=z$ or~$\overline{z}$, $\{\rho_\bi\}$ satisfies some
suitable restrictive condition~$\mathcal R$ (e.g. only finitely
many $\rho_\bi\neq1$). Suppose
$$
T_\bi:=T_{r_\bi,\theta_\bi,t_\bi,\sigma_\bi,\rho_\bi}(z),
$$
where $\bi=i_1\cdots i_k\in\mathcal I(m)$ and $r_\bi=r_{i_1}\cdots
r_{i_k}$, and
$$
K_\bi:=(T_{\bi|_1}\circ T_{\bi|_2}\circ\cdots\circ
T_{\bi|_k})(L_0),
$$
where $k$ is a nonnegative integer and $L_0$ is any fixed nonempty
compact set in~$\mathbb R^2$. Let
$$
K_\bi^{(p)}=\bigcup_{\bj\subseteq\bi,\;|\bj|=|\bi|+p} K_\bj
\quad(p=0,\,1,\,2,\,\dots).
$$
Then 
$$
h(K_\bi^{(p)},K_\bi^{(q)})\leqslant\delta_1(\lambda,\mathcal R)\,
r_\bi\,r_{\max}^p\, \mathrm{r}(\mathfrak S) \quad(q\geqslant p),
$$
where $r_{\max}=\max\{r_1,\dots,r_m\}$ and
$\delta_1(\lambda,\mathcal R)$ is a nonnegative real number
related to~$\lambda$ and~$\mathcal R$. Thus
$\{K_\bi^{(p)}\}_{p=1}^\infty$ approaches a nonempty compact set
in~$\mathbb R^2$ as $p\to+\infty$ for each $\bi\in\mathcal I(m)$,
denoted~$F_\bi$. It is easy to see
$$
F_\bi=\bigcup_{i=1}^m F_{\bi i}.
$$
We can also deduce
$$
\wth(F_\bi,\mathfrak S_\bi)\leqslant\delta_0(\lambda,\mathcal
R)\,r_\bi\,\mathrm{r}(\mathfrak S),
$$
where $\delta_0(\lambda,\mathcal R)$ is a nonnegative real number
related to~$\lambda$ and~$\mathcal R$. Therefore $\mathfrak
S(\mathcal W):=F_\bzero$ is a $\delta_0(\lambda,\mathcal
R)$-perturbation of~$\mathfrak S$ with a structure system
$\mathcal F=\{F_\bi:\:\bi\in\mathcal I(m)\}$.

We may also perturb $\mathscr K$, $\mathscr S$
 and~$\mathfrak S$ to get compact sets $\mathscr K'_n$, $\mathscr S'_n$
 and~$\mathfrak S'_n$
in~$\mathbb R^n$ ($n\geqslant2$) similarly.

\paragraph{{\bf 3.2.8.} {\it Remark}.}
(1) Let $P(\mathcal E)$ denote the probability of event~$\mathcal
E$. Suppose
\begin{align*}
\mathcal E=\{\mathcal W:\: \dim_H\mathfrak S(\mathcal
W)&=\dim_B\mathfrak S(\mathcal W)=s \\
 &\text{and}\ 0<\mathcal
H^s(\mathfrak S(\mathcal W))<+\infty\},
\end{align*}
 where
$r_1$, $\dots$, $r_m$ are fixed so that $s=\dim_S\mathfrak
S\leqslant2$ ($\sum_{j=1}^m r_j^s=1$), $\{\rho_\bi\}$ satisfies
some suitable restrictive condition~$\mathcal R$ or simply each
$\rho_\bi=1$, each~$\theta_\bi$ ($0\leqslant\theta_\bi<2\pi$),
each~$t_\bi$ ($|t_\bi|\leqslant\lambda$, $\lambda$ is a fixed
nonnegative real number) and each~$\sigma_\bi$ ($\sigma_\bi(z)=z$
or $\sigma_\bi(z)=\overline{z}$) possess some probability
distributions, for simplicity we assume that they are evenly
distributed. One question now arises: how much is~$P(\mathcal E)$?

Specially in Examples~3.2.7 for $\mathscr C\{a_\bi,b_\bi\}$
($a_0\leqslant a_\bi<b_\bi\leqslant b_0$), $\mathscr C'_n$,
$\mathscr K'_n$ and~$\mathscr S'_n$, which do not satisfy the
corresponding restrictive conditions that imply the open set
condition, each~$P(\mathcal E)$ seems to be~$1$.

(2) More generally we may perturb a self-similar set of~$\mathbb
R^n$ in~$\mathbb R^{n'}$ ($n'\geqslant n$) similarly to
perturbing~$\mathfrak S$.

(3) Let $\mathscr S'$ be a Sierpi\'{n}ski gasket which is obtained
from any triangle by following the method of getting the
Sierpi\'{n}ski gasket~$\mathscr S$ from a regular triangle. Then
$\mathscr S$ and~$\mathscr S'$ are $\delta$-perturbations of one
another. In fact we have the following more general result.

\paragraph{{\bf 3.2.9. Proposition.}}
{\it
 Let\/ $E$ and\/~$F$ be two invariant sets of IFSSs\/ $\{S_i:\:
i=1,\,\dots,\,m\}$ and\/ $\{T_i:\: i=1,\,\dots,\,m\}$
respectively. If\/ $\mathrm{Lip}\,S_i=\mathrm{Lip}\,T_i$ $(i=1$,
$\dots$, $m)$, then\/ $E$ and\/~$F$ are\/ $\delta$-perturbations
of one another.
 }
\par\noindent{\it Proof.}   It is easy to see that $S_\bi^{-1}\circ
T_\bi$ is an isometry. Let
$\mathrm{Lip}\,S_i=\mathrm{Lip}\,T_i=c_i$ ($i=1$, $\dots$, $m$).
Since $F_\bi=T_\bi(F)=S_\bi((S_\bi^{-1}\circ T_\bi)(F))$ we get
$$
\wth(F_\bi,E_\bi)=r_\bi\,\wth(E,F)=c_\bi\,\wth(E,F)=\alpha\,
c_\bi\,\mathrm{r}(E),
$$
where $\alpha=(\mathrm{r}(E))^{-1}\,\wth(E,F)$.
 \hfill $\Box$

\paragraph{{\bf 3.2.10.} {\it Remark}.}
(1) We easily see that a singleton is a $\delta$-perturbation of
any self-similar set which is not a singleton. So a singleton may
be considered as an extremely degenerate state of self-similar
sets.

(2) Let $C_\lambda$ be a $\lambda$-Cantor set which is determined
by $S_1(x)=\frac13x$, $S_2(x)=\frac13x+\frac23$ and
$S_\lambda(x)=\frac13x+\frac{\lambda}3$
($0\leqslant\lambda\leqslant2$) (see~\cite{Rao-Wen-1998}). Then by
Proposition~3.2.9 we know that $C_{{\lambda}'}$ is a
$\delta$-perturbation of~$C_\lambda$ with IFSS
$\{S_1,\,S_2,\,S_3\}$ ($0\leqslant\lambda'\leqslant2$).
From~\cite{Rao-Wen-1998} we see that $C_\lambda$
($0\leqslant\lambda\leqslant2$) have various Hausdorff dimensions.
Note that $C_0=C_2=\mathscr C$ and $C_1=[0,1]$.

(3) Normally we regard a self-similar set as an invariant set of
an IFSS which satisfies some kind of separation condition
(cf.~\cite[5.1]{Hutchinson-1981}). Now we may consider any
nonempty set~$F$ as a $\varDelta$-perturbation of any self-similar
set~$E$, where $\varDelta=\{\delta_\bi:\:\bi\in\mathcal I\}$.
However $\delta=\sup\{\delta_\bi:\:\bi\in\mathcal I\}$ is probably
equal to~$+\infty$ often.

(4) In fact in~(\ref{tag3.2-1}) the metric~$\wth$ may be replaced
by a simple quantity (the difference of diameters or radii) (refer
to~Proposition~2.6.3 and the proof of Theorem~3.2.4) and we still
have the same result as Theorem~3.2.4. That is rough, however, for
describing approximation degree of a fractal to a self-similar
set. When the perturbation is relatively small, we may imagine
that the quasi-self-similar set (defined in Definition~3.2.1)
visually and intuitively possesses approximate self-similarity.

(5) Similarly to the preceding discussion for self-similar sets we
may perturb some other kinds of fractals, e.g. a fractal with a
graph-directed construction (refer
to~\cite{Mauldin-Williams-1988}).

\paragraph{{\bf 3.3. Approximate self-similar sets.}}

According to~\cite{Hutchinson-1981} if a compact set can be
divided into a finite number of parts which are strictly similar
to the whole part~$F$ then all details have been determined and
each arbitrarily small part entirely reflects the whole. However,
if the strict similarity is not required, the determination of the
details will not exist any more. Now let us give the following
definition.

\paragraph{{\bf 3.3.1. Definition.}}
 Let $\mathcal F=\{F_\bi:\:\bi\in\mathcal I(\{m_j\})\}$ ($m_j$ is permitted to be $+\infty$)
  be a family of compact sets in~$\mathbb R^n$ satisfying
$$
F_\bi=\bigcup_{i=1}^{m_{k+1}}F_{\bi i}\quad(k=|\bi|).
$$
We call $\mathcal F$ a {\it structure system}\/ of~$F:=F_\bzero$.
Let $\varDelta=\{\delta_\bi:\:\bi\in\mathcal I(\{m_j\})\}$ be a
family of nonnegative real numbers. Let
$\delta:=\sup\{\delta_\bi:\: \delta_\bi\in\varDelta\}$.

 (1) If
 $$
 \widehat{h}(F_\bi,F)\leqslant\delta_\bi
 $$
for all~$\bi\in\mathcal I(\{m_j\})$, then $F$ is called a {\it
$\varDelta$-approximate self-similar set}. If $\delta<1$, then $F$
is also called a {\it $\delta$-approximate self-similar set}.

 (2) If
 $$
 \widehat{h}(F_{\bi i},F_\bi)\leqslant\delta_{\bi i}
 $$
for each~$\bi\in\mathcal I(\{m_j\})$, then $F$ is called a {\it
level-by-level $\varDelta$-approximate self-similar set}. If
$\delta<1$, then $F$ is also called a {\it level-by-level
$\delta$-approximate self-similar set}.

 (3) If
 $$
 \widehat{h}(F_{\bj},F_\bi)\leqslant\delta_{\bj}
 $$
for $\bi$, $\bj\in\mathcal I(\{m_j\})$ and $\bi\supset\bj$, then
$F$ is called a {\it uniformly}\/ ({\it level-by-level}\/) {\it
$\varDelta$-approximate self-similar set}. If $\delta<1$, then $F$
is also called a {\it uniformly}\/ ({\it level-by-level}\/) {\it
$\delta$-approximate self-similar set}.

\paragraph{{\bf 3.3.2.} {\it Remark}.}
(1) If $F$ is a $\delta$-approximate self-similar set
($0\leqslant\delta<\frac12$), then $F$ is a (uniformly)
level-by-level $2\delta$-approximate self-similar set. If $F$ is a
uniformly level-by-level $\varDelta$-approximate self-similar set,
then $F$ is a $\varDelta$-approximate self-similar set.

(2) Obviously if $\delta_i=0$ ($i=1$, $\dots$, $m_1$) then the
$\varDelta$-approximate self-similar set and (uniformly)
level-by-level $\varDelta$-approximate self-similar set~$F$ are
strictly self-similar. A compact set $F$ in $\mathbb R^n$ is
strictly self-similar if and only if $F$ is a $0$-approximate
self-similar set or a (uniformly) level-by-level $0$-approximate
self-similar set.

\paragraph{{\bf 3.3.3. Proposition.}}
{\it
 If\/ $F$ is a\/ $\varDelta$-perturbation of the self-similar set\/~$E$,
 where\/ $\varDelta=\{\delta_\bi\in[0,1):\:\bi\in\mathcal I(m)\}$,
 then\/
 $F$ is a\/ $2(\varDelta+\delta_\bzero)$-approximate self-similar
 set and is also a level-by-level\/ $\varDelta'$-approximate self-similar
 set, where\/ $2(\varDelta+\delta_\bzero)=\{2(\delta_\bi+\delta_\bzero):\:\bi\in\mathcal
 I(m)\}$ and\/ $\varDelta'=\{\delta'_{\bi i}=\delta_\bi+\delta_{\bi i}:\:\bi\in\mathcal
 I(m),\ i=1,\,\dots,\, m\}$. Specially,
 if\/ $\delta:=\sup\{\delta_\bi:\: \bi\in\mathcal
 I(m)\}<\frac14$, then\/ $F$ is a\/ $($uniformly\/$)$ level-by-level
 $4\delta$-approximate self-similar set.
 }
\par\noindent{\it Proof.}   Let $\mathcal F=\{F_\bi:\:\bi\in\mathcal
I(m)\}$ be a structure system of~$F$ so that
$$
\wth(F_\bi,E_\bi)\leqslant\delta_\bi\,c_\bi\,\mathrm{r}(E),
$$
where $E=E_\bzero$ is the invariant set of $\mathcal
S=\{S_\bi:\:\bi\in\mathcal I(m)\}$, $c_\bi=\mathrm{Lip}\,S_i$ and
$E_\bi=S_\bi(E)$. Then
$$
\wth(B,A_0)\leqslant\delta_\bi
$$
where $A_0=\frac1{\mathrm{r}(E)}E$ and
$B=\frac1{\mathrm{r}(E)}S_\bi^{-1}(F_\bi)$. Let
$B_0=\frac1{\mathrm{r}(B)}B$. Then
\begin{align*}
\widehat{h}(F_\bi,E)&=\widehat{h}(B,A_0)=\wth(B_0,A_0) \\
&\leqslant\wth(B_0,{\mathrm{r}(B)}\,B_0)+\wth(B,A_0) \\
&\leqslant|\mathrm{r}(B)-1|+\wth(B,A_0) \\
&\leqslant2\wth(B,A_0)\leqslant2\delta_\bi
\end{align*}
by Proposition~2.6.3. Therefore
$$
\widehat{h}(F_\bi,F)\leqslant\widehat{h}(F_\bi,E)+\widehat{h}(E,F)\leqslant
2(\delta_\bi+\delta_\bzero)
$$
and
\begin{equation}
\widehat{h}(F_\bi,F_\bj)\leqslant\widehat{h}(F_\bi,E)+\widehat{h}(E,F_\bj)\leqslant
2(\delta_\bi+\delta_\bj).\tag*{$\Box$}
\end{equation}


\paragraph{{\bf 3.3.4.} {\it Remark}.}
If $F$ is a $\delta$-quasi-self-similar set but $\delta$ is great
enough, then $F$ may not be a $\delta'$-approximate self-similar
set or level-by-level $\delta'$-approximate self-similar set for
some $\delta'\in[0,1)$.

\paragraph{{\bf 3.3.5. Definition.}}
 Let $\mathcal F^{l)}=\{F_\bi:\:\bi\in\mathcal I^{l)}(\{m_j\})\}$ ($l\in\mathbb P$)
 be a family of compact sets in~$\mathbb R^n$ satisfying
$$
F_\bi=\bigcup_{i=1}^{m_{k+1}}F_{\bi i}\quad(k=|\bi|),
$$
where $k=0$, $1$, $\dots$, $l-1$. We call $\mathcal F^{l)}$ a {\it
finite structure system}\/ of~$F:=F_\bzero$. Let
$\varDelta^{l)}=\{\delta_\bi:\:\bi\in\mathcal I^{l)}(\{m_j\})\}$
be a family of nonnegative real numbers which are less than~$1$.
Let $\delta:=\max\{\delta_\bi:\: \delta_\bi\in\varDelta^{l)}\}$.

 (1) If
 $$
 \widehat{h}(F_\bi,F)\leqslant\delta_\bi
 $$
for all~$\bi\in\mathcal I^{l)}(\{m_j\})$, then $F$ is called a
{\it $\varDelta^{l)}$-approximate self-similar set of
level~\/$l$}. $F$ is also called a {\it $\delta$-approximate
self-similar set of level~\/$l$}.

 (2) If
 $$
 \widehat{h}(F_{\bi i},F_\bi)\leqslant\delta_{\bi i}
 $$
for each~$\bi\in\mathcal I^{l-1)}(\{m_j\})$ and $i=1$, $\dots$,
$m_{k+1}$ ($k=|\bi|$), then $F$ is called a {\it level-by-level
$\varDelta^{l)}$-approximate self-similar set of level~\/$l$}. $F$
is also called a {\it level-by-level $\delta$-approximate
self-similar set of level~\/$l$}.

 (3) If
 $$
 \widehat{h}(F_{\bj},F_\bi)\leqslant\delta_{\bj}
 $$
for $\bi$, $\bj\in\mathcal I^{l)}(\{m_j\})$ and $\bi\supset\bj$,
then $F$ is called a {\it uniformly}\/ ({\it level-by-level}\/)
{\it $\varDelta^{l)}$-approximate self-similar set of
level~\/$l$}. $F$ is also called a {\it uniformly}\/ ({\it
level-by-level}\/) {\it $\delta$-approximate self-similar set of
level~\/$l$}.

\paragraph{{\bf 3.3.6.} {\it Remark}.}
In nature real objects rarely conform to Definition~3.3.1, but
they may conform to Definition~3.3.5.

\section{Comparison of fractals}\nopagebreak

{\it Only by comparing can one distinguish}. In Subsection~4.2 we
introduce some concepts to describe fractals using shape
differences by comparison. First we pose a problem which arises in
our life.

\paragraph{{\bf 4.1. Problem of shape vision error.}}

In the real world, errors always exist, which include matters of
human eyes. Here we suggest a problem concerning the ability of
man's visual sense.

We consider plane figures, which entirely get inside visual fields
of tested people and are not too far or too near from the tested
people, ignoring minor details. All the plane figures and
situations considered are as normal and simple as possible --- we
further make the following appointment: (i) The colors of figures
are black, the background is white and the brightness of light is
natural and moderate (of course, we may also consider effects of
these factors). (ii) Only usual Euclidean figures, such as
triangles, quadrilaterals (including rectangles, parallelograms,
trapezoids), polygons, ellipses, sectors and so on, are
considered. (iii) The figures are not too long and narrow and are
not too small or too large. Lengthes of sides of polygons have no
wide differences, etc. (iv) We only compare between figures
without any obvious distinctions. For example, we can consider
shape differences between an regular triangle (polygon) and other
triangles (polygons).

The problem is (let $A$ and~$B$ be two plane figures considered):
 How much is the critical value when $\widehat{h}(A,B)$ is
beyond it we can feel shapes of $A$ and $B$ are different and when
$h(A,B)$ is below it we feel shapes of $A$ and $B$ are the same?

\paragraph{{\it Remark}.}
(1) The critical value may be replaced by a small critical
interval.

(2) The results may be affected by specially appointed groups of
tested people to a certain degree. But we assume the tested people
are average.

(3) Complementary questions:  Whether or not are the results
affected by differences or ratios of radii (or diameters) of
figures and sizes or shapes of figures? Whether or not are the
results affected by distances between the figures and the tested
people?

\paragraph{{\bf 4.2. Atlases of fractals.}}

Let $0<\delta\leqslant\frac16$. If we look at Cantor's ternary
set~$\mathscr C$, which is the invariant set of~$\mathcal
S=\{S_1,\, S_2\}$ ($S_1(x)=\frac13x$ and
$S_2(x)=\frac13x+\frac23$) in~$\mathbb R$, we can easily find
$$
\widehat{h}(\mathscr C_\bi,C_k)=\widehat{h}(\mathscr
C,C_k)=\frac12\,\frac1{3^{k+1}}\leqslant\delta,
$$
where $\mathscr C_\bi=S_\bi(\mathscr C)$ ($\bi\in\mathcal I(2)$)
and $C_k=\mathcal S^k([0,1])$, when $k\geqslant k_0=\:
]\!-\log_3(6\delta)[$\,, where $]\alpha[$ denotes the smallest
integer more than or equal to~$\alpha$ ($\alpha\in\mathbb R$).
Thus we may consider $\mathscr C$ to possess a structure or form
of~$C_k$ for some $k\geqslant k_0$, if error~$\delta$ is
permitted. Now let us give the following

\paragraph{{\bf 4.2.1. Definition.}}
Let $\mathcal F=\{F_\bi:\:\bi\in\mathcal I(\{m_j\})\}$ be a
structure system of a compact set (fractal)~$F$ in~$\mathbb R^n$.
Assume $|F_\bi|\to0$ ($|\bi|\to+\infty$) and generally it is
required that $\mathcal F$ satisfies some kind of separation
condition, e.g. the open set condition, etc. Suppose $\mathcal
E^{(k)}:=\{E_j^{(k)}:\: j=1,\,\dots,\, p_{_{\scriptstyle k}}\}$
($k=0$, $1$, $2$, $\dots$) are families of compact sets
in~$\mathbb R^n$. Let $\varDelta^{(k)}=\{\delta_j^{(k)}:\:
j=1,\,\dots,\, p_{_{\scriptstyle k}}\}$, where
$0<\delta_j^{(k)}<1$ ($j=1,\,\dots,\, p_{_{\scriptstyle k}}$; and
$k=0$, $1$, $2$, $\dots$).

(1) If
$$
\widehat{h}(F_{\bi},E_{\bj_\bi})\leqslant\delta_{\bj_\bi}^{(k)}
$$
for each~$|\bi|\geqslant k$, where $\{\bj_\bi:\: \bi\in\mathcal
I(\{m_j\})\}=\{1,\,\dots,\, p_{_{\scriptstyle k}}\}$, then
$\mathcal E^{(k)}$ is called a {\it $k$-level\/
$\varDelta^{(k)}$-spline}\/ (or {\it $\delta^{(k)}$-spline}\/ if
$\delta_1^{(k)}=\cdots=\delta_{p_{_k}}^{(k)}=\delta^{(k)}$)
of\/~$\mathcal F$ (or~$F$) (or $E_1$, $\dots$, $E_{p_{_k}}$ are
called {\it $k$-level\/ $\varDelta^{(k)}$-splines}\/
of\/~$\mathcal F$ (or~$F$)) and we say $F$ possesses a {\it
$k$-th-level structure}\/~$\mathcal E^{(k)}$ (or {\it
structures}\/ $E_1$, $\dots$, $E_{p_{_k}}$) {\it of
error}\/~$\varDelta^{(k)}$ (or $\delta^{(k)}$ if
$\delta_1^{(k)}=\cdots=\delta_{p_{_k}}^{(k)}=\delta^{(k)}$). If
$F$ possesses a $0$-th-level structure~$\mathcal E^{(0)}$ of
error~$\varDelta^{(0)}$ (or $\delta^{(0)}$) then we also say $F$
possesses a {\it structure\/~$\mathcal E=\mathcal E^{(0)}$ of
error}\/~$\varDelta=\varDelta^{(0)}$ (or $\delta=\delta^{(0)}$)
and $\mathcal E$ is called a {\it $\varDelta$-spline}\/ (or {\it
$\delta$-spline}\/).

(2) If for any $\delta>0$, $F$ possesses a $k$-th-level
structure~$\mathcal E^{(k)}$, consisting of~$p_{_{\scriptstyle
k}}$ compact sets, of error~$\delta$, then $F$ is said to possess
a {\it $k$-th-level $p_{_{\scriptstyle k}}$-structure}. If
$p_{_{\scriptstyle k}}=1$ then $F$ is said to possess a {\it
single structure in $k$-th-level}. If $F$ possesses a $0$-th-level
$p_{_{\scriptstyle 0}}$-structure then $F$ is said to possess a
{\it $p$-structure}\/ ($p=p_{_{\scriptstyle 0}}$). If $p=1$ then
$F$ is said to possess a {\it single structure}.

(3) Let $\lambda>0$. If for all $\bi\in\mathcal I(\{m_j\})$ we
have $|F_\bi|\leqslant\lambda|F|$ then a $k$-th-level
structure~$\mathcal E^{(k)}$ of error~$\varDelta^{(k)}$ (or
$\delta^{(k)}$) of~$F$ is called a {\it $\lambda$-degree structure
of of error}\/~$\varDelta=\varDelta^{(k)}$ (or
$\delta=\delta^{(k)}$) and a $k$-level $\varDelta$-spline (or
$\delta$-spline) of~$F$ is called a {\it $\lambda$-degree
$k$-level\/ $\varDelta$-spline}\/ (or {\it $\delta$-spline}\/).


\paragraph{{\bf 4.2.2.} {\it Remark}.}
(1) A self-similar set~$F$ may be considered to possess a
$0$-spline~$F$. But now we make a convention that in~$\mathbb R^n$
the splines should be geometric patterns consisting of finite
formal Euclidean figures, which are called {\it Euclidean
patterns}. It is easy to see that a self-similar set still
possesses a single structure.

(2) The splines of fractals~$F$ are not unique. A standard of
searching for splines of~$F$ is trying to find splines of~$F$
which can help us to see and understand details of~$F$
approximately. A well-chosen spline of a fractal~$F$ is called an
{\it atlas}\/ of~$F$.

\paragraph{{\bf 4.2.3.} {\it Examples}.}
(1) Let $0<\delta\leqslant\frac12$. We consider the von Koch
curve~$\mathscr K$, see Example~3.2.7(2). Then
$$
\widehat{h}(S_\bi(\mathscr K),K^{(k)})=\widehat{h}(\mathscr
K,K^{(k)})\leqslant\frac12\,\frac1{3^k}\leqslant\delta,
$$
if $k\geqslant k_1=\: ]\!\log_3\frac{1}{2\delta}[$\,, where
$]\alpha[$ denotes the smallest integer more than or equal
to~$\alpha$ ($\alpha\in\mathbb R$). Hence $K^{(k)}$ is a
$\delta$-spline of~$\mathscr K$ for each $k\geqslant k_1$.

(2) Let $0<\delta\leqslant\frac14$. We consider the Sierpi\'{n}ski
gasket~$\mathscr S$, see Example~3.2.7(3). Let $A^{(k)}=\mathcal
S^{(k)}(A_0):=\bigcup_{\bi\in\mathcal I_k(3)} S_\bi(A_0)$ ($k=0$,
$1$, $2$, $\dots$). Then
$$
\widehat{h}(S_\bi(\mathscr S),A^{(k)})=\widehat{h}(\mathscr
S,A^{(k)})\leqslant\frac14\,\frac1{2^k}\leqslant\delta,
$$
if $k\geqslant k_2=\: ]\!\log_2\frac{1}{4\delta}[$\,. Hence
$A^{(k)}$ is a $\delta$-spline of~$\mathscr S$ for each
$k\geqslant k_2$.

\paragraph{{\bf 4.2.4.} {\it Remark}.}
(1) If $F$ is a $\delta$-approximate self-similar set
($0\leqslant\delta<1$), then there is a
$(\delta+\varepsilon)$-spline~$E$ of~$F$ for any $\varepsilon>0$
satisfying $\delta+\varepsilon<1$, where $E$ is an Euclidean
pattern. But $F$ may not possess any single structures. And a
$\delta$-quasi-self-similar set ($\delta>0$) may not possess any
single structures either (cf.~Remark~4.2.2(1)).

(2) The construction object in a graph directed construction
(refer to~\cite{Mauldin-Williams-1988}) may possess a
$p$-structure ($p\leqslant n$).

\paragraph{}

In order to distinguish the simplity and complexity of details of
a fractal we give the following

\paragraph{{\bf 4.2.5. Definition.}}
If a compact set~$F$ possesses a structure system $\mathcal
F=\{F_\bi:\:\bi\in\mathcal I\}$ and there exist a finite number of
compact sets $E_1$, $\dots$, $E_q$ such that for each
$\bi\in\mathcal I$ there exists $\bj_\bi\in\{1,\,\dots,\,q\}$
satisfying
$$
\widehat{h}(F_\bi,E_{\bj_\bi})\to0\quad(|\bi|\to+\infty),
$$
then $F$ is said to {\it approach finite structure}.

\paragraph{{\bf 4.2.6.} {\it Remark}.}
A self-similar set and a graph directed construction object
approach finite structure, but a $\delta$-quasi-self-similar set
and a $\delta$-approximate self-similar set may not approach
finite structure.

\paragraph{{\bf 4.2.7. Fractal indices.}}
 Let $F$ be a compact set (fractal) in~$\mathbb R^n$ and
 let $\mathcal F=\{F_\bi:\:\bi\in\mathcal I\}$ be a structure
system of~$F$. Assume $\delta$, $\lambda>0$ and $k$ is a
nonnegative integer.

(1) Denote
\begin{align*}
N(\delta,\lambda;\,F):=\min\{\sharp\,E_\delta(\lambda):\:
&E_\delta(\lambda) \text{ is a $\lambda$-degree} \\
 &\text{ structure of~$F$ of error~$\delta$} \}
\end{align*}
($\sharp\,A$ denotes the cardinal number of~$A$). Then
$N(\delta,\lambda;\,F)$ is a decreasing function of~$\lambda$. Let
$$
N_\delta(F)=N(\delta,F):=\sup_{\lambda>0}\,N(\delta,\lambda;\,F).
$$
Then
$$
N_\delta(F)=\lim_{\lambda\to0^+}\,N(\delta,\lambda;\,F).
$$
The faster the growth of $N_\delta(F)$ on~$\frac1\delta$ is, the
more complex the detail of~$F$ is; and the slower the growth is,
the simpler the detail is. So we call $N_\delta(F)$ the {\it
fractal $\delta$-index}\/ ({\it index function}\/) of~$F$
($(N_\delta(F))^{-1}$ is called the {\it $\delta$-self-similarity
index} ({\it function}\/) of~$F$) and call the (upper, lower)
growth order of $N_\delta(F)$ on~$\frac1\delta$
 the ({\it upper}, {\it lower}\/) {\it fractal order}\/ of~$F$,
where the {\it upper}\/ and {\it lower growth order}\/ are
$$
\overline{\rho}(F):=\limsup_{\delta\to0^+}\frac{N_\delta(F)}{-\log\delta}\quad
\text{and}\quad\underline{\rho}(F):=\liminf_{\delta\to0^+}\frac{N_\delta(F)}{-\log\delta}
$$
respectively, and if $\overline{\rho}(F)=\underline{\rho}(F)$ then
the {\it growth order}\/ is $\rho(F)=\overline{\rho}(F)$.

(2) Denote
$$
N_\delta^{(k)}(\mathcal F):=\min\{\sharp\,E_\delta^{(k)}:\:
E_\delta^{(k)} \text{ is a $k$-th-level}
 \text{ $\delta$-spline of~$\mathcal F$} \}.
$$
Let
$$
\overline{N}_\delta(\mathcal F):=\limsup_{k\to+\infty}
N_\delta^{(k)}(\mathcal F)\quad
\text{and}\quad\underline{N}_\delta(\mathcal
F):=\liminf_{k\to+\infty} N_\delta^{(k)}(\mathcal F).
$$
Define
\begin{align*}
\overline{N}_\delta(F)&:=\min\{\overline{N}_\delta(\mathcal F):\:
\mathcal F
 \text{ is a structure system of~$F$} \}, \\
\underline{N}_\delta(F)&:=\min\{\underline{N}_\delta(\mathcal
F):\: \mathcal F
 \text{ is a structure system of~$F$} \} \\
 \intertext{and}
N_\delta^{(0)}(F)&:=\min\{N_\delta^{(0)}(\mathcal F):\: \mathcal F
 \text{ is a structure system of~$F$} \},
\end{align*}
which are called {\it upper}, {\it lower}\/ and {\it whole fractal
$\delta$-indices}\/ ({\it index functions}\/) of~$F$ respectively.
We also call $(\overline{N}_\delta(F))^{-1}$,
$(\underline{N}_\delta(F))^{-1}$ and $(N_\delta^{(0)}(F))^{-1}$
{\it upper}, {\it lower}\/ and {\it whole $\delta$-self-similarity
indices}\/ ({\it index functions}\/) of~$F$ respectively.

(3) Let $F^{(j)}$ ($j\in J$, $J$ is an index set) be fractals
in~$\mathbb R^n$. Denote $F_J:=\{F^{(j)}:\:j\in J\}$. Assume
$\mathcal F_j:=\{F_{\bi^{(j)}}^{(j)}:\:\bi^{(j)}\in \mathcal
I^{(j)}\}$ is a structure system of~$F^{(j)}$ ($j\in J$). Denote
$\mathcal F_J:=\{\mathcal F_j:\:j\in J\}$, called a {\it structure
system}\/ of~$F_J$. A family $\mathcal E_J^{(k)}=\{E_l^{(k)}:\:
l=1,\,2,\,\dots\,p_{_{\scriptstyle k}} \}$ of compact sets is
called a {\it $\delta$-spline of}\/~$\mathcal F_J$ (or~$F_J$) if
$$
\widehat{h}(F_{\bi^{(j)}}^{(j)},E_{l(j,\bi^{(j)})}^{(k)})\leqslant\delta
$$
for $|\bi^{(j)}|\geqslant k$ ($\bi^{(j)}\in\mathcal I^{(j)}$,
$j\in J$), where
$$
\{l(j,\bi^{(j)})=l_{\bi^{(j)}}^{(j)}:\:j\in
J,\,\bi^{(j)}\in\mathcal I^{(j)},\,|\bi^{(j)}|\geqslant
k\}=\{1,\,2,\,\dots,\,p_{_{\scriptstyle k}}\}.
$$
Assume $p_{_{\scriptstyle k}}$ is the smallest one such that
$\mathcal E_J^{(k)}$ is a $\delta$-spline of~$\mathcal F_J$.
Suppose
\begin{align*}
\{E_{l_t}^{(k)}:\:&t=1,\,2,\,\dots,\,q_{_{\scriptstyle k}};\;
1\leqslant
l_1<l_2<\cdots<l_{q_{_{\scriptstyle k}}}\leqslant p_{_{\scriptstyle k}}\} \\
 =\{E_l^{(k)}:\:&1\leqslant
l\leqslant p_{_{\scriptstyle k}},\text{ and for each $j\in J$
} \\
&\text{ there exists $\bi^{(j)}\in\mathcal I{(j)}$ such that
$\widehat{h}(F_{\bi^{(j)}}^{(j)},E_l^{(k)})\leqslant\delta$}\}
\end{align*}
is a common $\delta$-spline of~$\mathcal F_J$. Denote
$$
\gamma_{_{\scriptstyle k}}(\mathcal
F_J,\delta):=\frac{q_{_{k}}}{p_{_{k}}}.
$$
Let
$$
\overline{\gamma}(\mathcal F_J,\delta):=\limsup_{k\to+\infty}
\gamma_{_{\scriptstyle k}}(\mathcal F_J,\delta)\quad
\text{and}\quad\underline{\gamma}(\mathcal
F_J,\delta):=\liminf_{k\to+\infty} \gamma_{_{\scriptstyle
k}}(\mathcal F_J,\delta).
$$
Define
\begin{align*}
\overline{\gamma}(F_J,\delta)&:=\max\{\overline{\gamma}(\mathcal
F_J,\delta):\: \mathcal F_J \text{ is a structure system of~$F_J$} \}, \\
\underline{\gamma}(F_J,\delta)&:=\max\{\underline{\gamma}(\mathcal
F_J,\delta):\: \mathcal F_J \text{ is a structure system of~$F_J$} \} \\
 \intertext{and}
{\gamma}_0(F_J,\delta)&:=\max\{{\gamma}_0(\mathcal F_J,\delta):\:
\mathcal F_J \text{ is a structure system of~$F_J$} \},
\end{align*}
which are called {\it upper}, {\it lower}\/ and {\it whole
$\delta$-similarity indices}\/ ({\it similarity index
functions}\/) of~$F_J$ (or $F^{(j)}$ ($j\in J$)) respectively.
If\/
$\overline{\gamma}(F_J,\delta)=\underline{\gamma}(F_J,\delta)$
then it is called the {\it $\delta$-similarity index}\/ ({\it
similarity index function}\/), denoted ${\gamma}(F_J,\delta)$.
Define
\begin{align*}
\overline{\gamma}(F_J)&:=\limsup_{\delta\to0^+}\overline{\gamma}(F_J,\delta), \\
\underline{\gamma}(F_J)&:=\limsup_{\delta\to0^+}\underline{\gamma}(F_J,\delta) \\
 \intertext{and}
 {\gamma}_{_{\scriptstyle 0}}(F_J)&:=\limsup_{\delta\to0^+}{\gamma}_{_{\scriptstyle 0}}(F_J,\delta),
\end{align*}
which are called {\it upper}, {\it lower}\/ and {\it whole
similarity indices}\/ of~$F_J$ (or $F^{(j)}$ ($j\in J$))
respectively. If the above superior limits are changed into
inferior limits then they are called {\it upper}, {\it lower}\/
and {\it whole strong similarity indices}\/ of~$F_J$ (or $F^{(j)}$
($j\in J$)) respectively.

It is obvious to see that
$0\leqslant\underline{\gamma}(F_J,\delta)\leqslant\overline{\gamma}(F_J,\delta)\leqslant1$,
$0\leqslant{\gamma}_{_{\scriptstyle 0}}(F_J,\delta)\leqslant1$,
$0\leqslant\underline{\gamma}(F_J)\leqslant\overline{\gamma}(F_J)\leqslant1$
and $0\leqslant{\gamma}_{_{\scriptstyle 0}}(F_J)\leqslant1$.

\paragraph{{\bf 4.2.8.} {\it Example}.}
If $F$ is a self-similar set then $N_\delta(F)=1$. If
$F=\bigcup_{i=1}^m F_i$ is the construction object in a graph
directed construction, then $N_\delta(F)\leqslant m$, and usually
when $\delta$ is small enough we have $N_\delta(F)=m$.

\paragraph{{\bf 4.2.9.} {\it Remark}.}
In general, computing $N_\delta(F)$,
$\overline{\gamma}(F_J,\delta)$, $\underline{\gamma}(F_J,\delta)$,
${\gamma}_0(F_J,\delta)$, $\overline{\gamma}(F_J)$,
$\underline{\gamma}(F_J)$ and ${\gamma}_0(F_J)$ is a very
difficult job. But for some $\delta>0$ and a structure system
$\mathcal F=\{F_\bi:\:\bi\in\mathcal I\}$ of~$F$, computing
$\overline{\gamma}(\mathcal F_J,\delta)$,
$\underline{\gamma}(\mathcal F_J,\delta)$ and ${\gamma}_0(\mathcal
F_J,\delta)$ may be a piece of operable work.

\paragraph{{\bf 4.3. Some examples and remarks.}}

A cookie-cutter set~$E$ 
 is a quasi-self-similar set in
the sense that every small piece of~$E$ can be uniformly expanded
to a standard size and then mapped quasi-isometrically back
into~$E$ and that $E$ can also be quasi-isometrically contracted
to any small part of~$E$ (see \cite{Falconer-1989},
\cite[Chapter~4]{Falconer-1997} and \cite{McLaughlin-1987}). Hence
$E$ is an $s$-set with $s=\dim_H E=\dim_B E$)(see
\cite[Corollary~4.6]{Falconer-1997}). A lot of research work has
been done in this aspect, for example, connecting with
thermodynamic formalism, we may refer to \cite{Bedford-1991},
\cite{Bowen-1975}, \cite{Falconer-1997}, \cite{Ruelle-1978},
\cite{Ruelle-1982} and \cite{Sinai-1972}, etc.

\paragraph{{\bf 4.3.1.}}
First let us consider an example which was given in
\cite[Section~4.1]{Falconer-1997}: Let $F$ be a cookie-cutter  set
that is an invariant set of~$g_1$ and~$g_2:\:[0,1]\to[0,1]$, where
$$
g_1(x)=\frac13x+\frac1{10}x^2\quad\text{and}\quad
g_2(x)=\frac13x+\frac23-\frac1{10}x^2,
$$
which is a nonlinear perturbation of Cantor's ternary set. Of
course we may consider a more general case:
$$
g_1(x)=\frac13x+ax^\alpha\quad\text{and}\quad
g_2(x)=\frac13x+\frac23-bx^\beta,
$$
where $a$ and~$b$ are two small positive numbers and $\alpha$,
$\beta>1$. We pose the following questions:
\begin{enumerate}
\item[(1)] Does $F$ approach finite structure?
\item[(2)] If we change $g_1$ and~$g_2$ above into
$$
g_1(x)=\frac13x+\varphi(x)\quad\text{and}\quad
g_2(x)=\frac13x+\frac23-\psi(x),
$$
where $\varphi$ and~$\psi$ are chosen such that the new $g_1$
and~$g_2$ are also contraction maps from $[0,1]$ to itself and the
invariant set~$F$ of $g_1$ and~$g_2$ approaches finite structure,
then what features should $\varphi$ and~$\psi$ possess? (Obviously
if $\varphi(x)=\psi(x)\equiv0$ then $F$ becomes Cantor's ternary
set, which approaches finite structure.)
\end{enumerate}

\paragraph{{\bf 4.3.2.}}
Consider the logistic map
$$
f(x)=\lambda x(1-x),
$$
where $\lambda$ is a positive constant. It is an important one
dimensional dynamic system, which has been deeply and
systematically studied. There is a known universal constant, the
Feigenbaum constant, for example (refer to
\cite[Section~13.2]{Falconer-2003}, \cite{Feigenbaum-1978} and
\cite{Feigenbaum-1979}, etc).

If $\lambda>2+\sqrt5$, then $f:\:[0,a]\cup[1\!-\!a,1]\to[0,1]$,
where $a=\frac12-\sqrt{\frac14-\frac1\lambda}$, gives a
cookie-cutter set, denoted~$E_\lambda$ (refer to
\cite[Section~4.1]{Falconer-1997} and
\cite[Section~13.2]{Falconer-2003}). A question similar
to~4.3.1(1) may be raised:
\begin{enumerate}
\item[(1)] Do $E_\lambda$ ($\lambda>2+\sqrt5$) approach finite structure?
\end{enumerate}

When $\lambda=\lambda_\infty\approx3.570$, the
attractor~$E_{\lambda_\infty}$ is a set of Cantor type, whose
Hausdorff dimension is about~0.538 (see
\cite[Section~13.2]{Falconer-2003}). One more similar question may
be mentioned:
\begin{enumerate}
\item[(2)] Does $E_{\lambda_\infty}$ approach finite structure?
\end{enumerate}

\paragraph{{\bf 4.3.3.}}
Iterating rational functions (generally, meromorphic functions) in
the complex plane~$\mathbb C$ a large number of fractals can
emerge. Specially a quadratic polynomial
$$
f_c(z)=z^2+c\quad(c\in\mathbb C)
$$
may bring about a colorful dynamic system, which has been one of
central issues in the research of complex analytic dynamic
systems, and many splendid and deep results have been found. Here
we only mention a few of references: \cite{Beardon-1991},
\cite{Bergweiler-1993}, \cite{Blanchard-1984},
\cite{Carleson-Gamelin-1993}, \cite{Lyubich-1997},
\cite{McMullen-1994}, \cite{McMullen-1998}, \cite{Milnor-1999},
\cite{Morosawa-Nishimura-Taniguchi-Ueda-2000},
\cite{Peitgen-Richter-1986}, \cite{Shishikura-1998},
\cite{Steinmetz-1993}, \cite{Sullivan-1985}, \cite{Tan-2000}, etc.

Let $\mathcal M$ denote the Mandelbrot set and $\mathcal J_c$
denote the Julia set of~$f_c$. It is known that Julia sets of
hyperbolic rational functions are quasi-self-similar in the sense
mentioned above in this subsection (see
\cite[Theorem~8.6]{Blanchard-1984} and
\cite[p.~742]{Sullivan-1983}) (it is not always true for all
$\mathcal J_c$, see \cite{Jarvi-1997}). In \cite{Tan-1990} Lei Tan
obtained a kind of similarity between $\mathcal M$ and $\mathcal
J_c$ for Misiurewicz points~$c$.

If $|c|>\frac{5+2\sqrt6}4$, then $\mathcal J_c$ is a cookie-cutter
set in plane. A similar question may still be asked:
\begin{enumerate}
\item[(1)] Do $\mathcal J_c$ ($|c|>\frac{5+2\sqrt6}4$) approach finite structure?
\end{enumerate}

Let $F$ be $\mathcal M$, $\mathcal J_c$ or a piece of them and let
$F_J$ be a family of fractals chosen from $\mathcal M$, $\mathcal
J_c$ or pieces of them. Now our questions are:
\begin{enumerate}
\item[(2)] Find the $\delta$-self-similarity index~$(N_\delta(F))^{-1}$ of~$F$;
\item[(3)] Find upper, lower and whole
$\delta$-similarity indices $\overline{\gamma}(F_J,\delta)$,
$\underline{\gamma}(F_J,\delta)$ and ${\gamma}_0(F_J,\delta)$
 and upper, lower and whole similarity indices
$\overline{\gamma}(F_J)$, $\underline{\gamma}(F_J)$ and
${\gamma}_0(F_J)$ of~$F_J$.
\end{enumerate}
 The same work may be done in~4.3.1 and~4.3.2.

\section{Notes on tilings, patterns, packings and crystals}\nopagebreak

The research on tilings, patterns and packings has been
extensively carried out, which has a long history and was also
motivated by Hilbert's 18th problem and theories about the
structure of solid matter. From \cite{Brass-Moser-Pach-2005},
\cite{Croft-Falconer-Guy-1994} and \cite{Goodman-O'Rourke-2004} we
can see numerous problems on the areas still remain open. For
tilings and patterns we may refer to
\cite{Grunbaum-Shephard-1980}, \cite{Grunbaum-Shephard-1987} and
\cite{Schattschneider-Senechal-2004}, etc. For packings we may
refer to \cite{Borozky-2004}, \cite{Conway-Sloane-1999},
\cite{FejesToth-2004}, \cite{Rogers-1964} and \cite{Zong-1999},
etc. For (mathematical) crystallography we may refer to
\cite{Axel-Gratias-1995}, \cite{Kittel-2005}, \cite{Senechal-1995}
and \cite{Senechal-2004}, etc. For some tilings the tiles may be
fractals, see e.g. \cite{Bandt-1991} and \cite{Mandelbrot-1982}.
The rigidity of tilings was considered in \cite{Kenyon-1992b}. For
self-similar and self-affine tilings we refer to e.g.
\cite{Kenyon-1990}, \cite{Kenyon-1992a}, \cite{Kenyon-1994},
\cite{Kenyon-1996}, \cite{Kenyon-Li-Strichartz-Wang-1999},
\cite{Lagarias-Wang-1996}, \cite{Lagarias-Wang-1996a},
\cite{Lau-Rao-2003}, \cite{Solomyak-1997},
\cite{Strichartz-Wang-1999} and \cite{Thurston-1989}, etc.

\paragraph{{\bf 5.1. Tilings with quasi-prototiles.}}
It is known that a {\it tiling}\/ (of~$\mathbb R^n$) is a
countable family of closed sets (usually {\it bodies}, which are
bounded and are closures of their interiors) in~$\mathbb R^n$,
whose union is the whole space and whose interiors are pairwise
disjoint (refer to~\cite{Grunbaum-Shephard-1987}
and~\cite{Schattschneider-Senechal-2004}). A {\it tile}\/ is an
element of a tiling. A body~$T$ {\it tiles}\/~$\mathbb R^n$ means
that there is a tiling of~$\mathbb R^n$ whose tiles are congruent
copies of~$T$.

\paragraph{{\bf 5.1.1. Tilings with quasi-prototiles.}}
Let $\mathcal T=\{T_i:\:i\in I\}$ be a tiling. Let $\mathcal P=
\{P_j:\:j\in J\}$ be a set of at most countable nonempty compact
subsets (usually bodies) of~$\mathbb R^n$ and $\{\delta_j\}_{j\in
J}$ be a (countable or finite) sequence of nonnegative real
numbers. Suppose for each $i\in I$ there exists $j\in J$ such that
$$
\wth(T_i,P_j)\leqslant\delta_j\,\mathrm{r}(P_j),
$$
where we assume every $P_j$ ($j\in J$) is taken at least once or
more strictly any $P_j$ ($j\in J$) can not be lost. Then  we call
$\mathcal T$ a {\it tiling with\/ $\{\delta_j\}$-quasi-prototile
types}\/ $P_j$ ($j\in J$) or {\it with quasi-prototile types}\/
$\langle P_j,\delta_j\rangle$ ($j\in J$). We may also say that
$\mathcal T$ is a {\it $\{\delta_j\}$-quasi-tiling with prototile
types}\/ $P_j$ ($j\in J$) or {\it prototile types}\/ $P_j$ ($j\in
J$) {\it admits a\/ $\{\delta_j\}$-quasi-tiling}\/~$\mathcal T$. A
tiling with a single quasi-prototile type is called a {\it
monohedral quasi-tiling}. A tiling with $k$ quasi-prototile types
is called a {\it $k$-hedral quasi-tiling}. If $\delta_j=\delta$
for all $j\in J$, then $\mathcal T$ is called a {\it tiling with\/
$\delta$-quasi-prototile type}\/ ({\it set}\/)~$\mathcal P$ or
{\it $\delta$-quasi-prototile types}\/ $P_j$ ($j\in J$) or is
called a {\it $\delta$-quasi-tiling with prototile types}\/ $P_j$
($j\in J$). If $P_j\in\mathcal T$ ($j\in J$) then we also call
(quasi-)prototile types $P_j$ ($j\in J$) (\nolinebreak{\it
quasi-}\nolinebreak){\it prototiles}. Usually we may not
distinguish between (\nolinebreak quasi-\nolinebreak)prototile
types and (\nolinebreak quasi-\nolinebreak)prototiles and usually
we assume the cardinal number $\sharp\nolinebreak\,J=k$ of~$J$ is
finite.

\paragraph{{\bf 5.1.2. Quasi-symmetry groups of tilings.}}
Let $\mathcal T$ be a tiling and $\lambda\geqslant0$. If the
isometry~$\varphi$ of~$\mathbb R^n$ satisfies that for any $i\in
I$ there exist $i'$ and $i''\in I$ such that
\begin{align*}
h(\varphi(T_i), T_{i'})&\leqslant\lambda\,\mathrm{r}(T_i) \\
\intertext{and}
 h(T_i,\varphi(T_{i''}))&\leqslant\lambda\,\mathrm{r}(T_i),
\end{align*}
then $\varphi$ is called a {\it $\lambda$-quasi-symmetry}\/
of~$\mathcal T$. A group~$\mathcal G$ consisting of
$\lambda$-quasi-symmetries of~$\mathcal T$ is called a {\it
$\lambda$-quasi-symmetry group}\/ of~$\mathcal T$. Note that a
$0$-quasi-symmetry of~$\mathcal T$ is a symmetry of~$\mathcal T$
and a $0$-quasi-symmetry group of~$\mathcal T$ is a symmetry group
of~$\mathcal T$.

\paragraph{{\bf 5.1.3. Quasi-self-similar tilings}
 \rm{(cf. \cite{Kenyon-1990} and \cite{Schattschneider-Senechal-2004}).}}

Recall that a {\it hierarchical tiling}\/ is a tiling whose tiles
(called {\it level-$0$} tiles) can be composed into larger tiles,
called {\it level-$1$} tiles, whose level-$1$ tiles can be
composed into {\it level-$2$} tiles, and so on ad infinitum (see
\cite[p.\,\nolinebreak66]{Schattschneider-Senechal-2004}).

(1) Let $\mathcal Q= \{Q_i:\:i\in I\}$ be a set of at most
countable nonempty compact subsets (usually bodies) of~$\mathbb
R^n$ and let $\lambda_{ij}\geqslant0$ ($i\in I$ and $j\in\mathbb
N$) ($\mathbb N$ denotes the set of all nonnegative integers). We
define a {\it $\{\lambda_{ij}\}$-quasi-self-similar tiling with a
quasi-prototype}\/ ({\it set}\/)~$\mathcal Q$ as a hierarchical
tiling~$\mathcal T$ satisfying that for any level-$j$ tile~$T_j$
of~$\mathcal T$ there exists a quasi-prototype~$Q_i\in\mathcal Q$
such that
$$
\widehat{h}(T_j,Q_i)\leqslant\lambda_{ij}.
$$
We also call $\mathcal T$ a {\it quasi-self-similar tiling with
a\/ $\{\lambda_{ij}\}$-quasi-prototype}\/ $\mathcal Q$. If
$\lambda_{ij}=\lambda$ for all $i\in I$ and $j\in\mathbb N$, then
we say that $\mathcal T$ is a {\it $\lambda$-quasi-self-similar
tiling with a quasi-prototype}\/ $\mathcal Q$ or a {\it
quasi-self-similar tiling with a\/ $\lambda$-quasi-prototype}\/
$\mathcal Q$. Usually we assume $\sharp\nolinebreak\,I=k$ is
finite.

(2) Let $\lambda_j\geqslant0$ ($j\in\mathbb P$). If a hierarchical
tiling~$\mathcal T$ satisfies that for any level-$j$ tile~$T_j$
($j\in\mathbb P$) and each level-$(j\!-\!1)$
tile~$T_{j-1}\subseteq T_j$ such that
$$
\widehat{h}(T_j,T_{j-1})\leqslant\lambda_j\quad(j\in\mathbb P),
$$
then $\mathcal T$ is called a {\it level-by-level\/
$\{\lambda_j\}$-quasi-self-similar tiling}. If $\lambda_j=\lambda$
for all $j\in\mathbb P$, then we say that $\mathcal T$ is a {\it
level-by-level\/ $\lambda$-quasi-self-similar tiling}.

(3) Let $\mathcal Q= \{Q_i:\:i\in I\}$ be a set of at most
countable nonempty compact subsets (usually bodies) of~$\mathbb
R^n$. Let $\lambda_{ij}\geqslant0$, $c_{ij}\geqslant1$ ($i\in I$
and $j\in\mathbb P$) and $\sup\{c_{ij}:\:j\in\mathbb P\}>1$ ($i\in
I$). If a hierarchical tiling $\mathcal T$ satisfies that for any
level-$j$ tile~$T_j$ ($j\in\mathbb P$) there exists
$Q_i\in\mathcal Q$ and a similitude~$\varphi_{ij}$ of lipschitz
constant $c_{ij}$ such that
$$
\wth(T_j,\varphi_{ij}(Q_i))\leqslant\lambda_{ij}
c_{ij}\,\mathrm{r}(Q_i),
$$
then $\mathcal T$ is called a {\it
$\{\lambda_{ij}\}$-quasi-self-similar tiling with a
quasi-prototype}\/ $\mathcal Q$ or a {\it quasi-self-similar
tiling with a\/ $\{\lambda_{ij}\}$-quasi-prototype\/ $\mathcal Q$
of ratio}\/ $\{c_{ij}\}$ ({\it $\lambda$-quasi-self-similar
tiling}\/ if $\lambda_{ij}=\lambda$ for all $i\in I$ and
$j\in\mathbb P$) (of {\it ratio}\/~$\{c_i\}$ if $c_{ij}=c_i^j$ for
all $i\in I$ and $j\in\mathbb P$, or of {\it ratio}\/~$c$ if
furthermore $c_i=c$ for all $i\in I$). Usually we assume
$\sharp\nolinebreak\,I=k$ is finite.

(4) If there is a quasi-self-similar tiling~$\mathcal T$ with a
 single $\{\lambda_j\}$-quasi-prototype $Q$ (of ratio $\{c_j\}$), then $Q$ is called
 a {\it rep $\{\lambda_j\}$-quasi-tile type}\/ ({\it $k$-rep
 $\{\lambda_j\}$-quasi-tile type}\/ if every level-$j$ title consists of $k$ level-$(j\!-\!1)$
 tiles, $j\in\mathbb N$ and $k\geqslant2$ is a natural number
 independent of $j\in\mathbb N$) (of {\it ratio}\/~$\{c_j\}$).
 If $Q\in\mathcal T$ then the ($k$-)rep $\{\lambda_j\}$-quasi-tile type $Q$
 is also called a ($k$-){\it rep $\{\lambda_j\}$-quasi-tile}.
 Usually we do not distinguish between these two notions.
 Similar statements can be made for the case $\lambda_j=\lambda$ for all $j\in\mathbb N$.

\paragraph{{\bf 5.1.4.} {\it Example}.}
Let
$$
S_{jk}:=\{(x,y):\:j\leqslant x\leqslant j+1,\ k\leqslant
y\leqslant k+1\},
$$
where $j$, $k\in\mathbb Z$ ($\mathbb Z$ denotes the set of all
integers). Then $\mathcal S=\{S_{jk}:\:j,\,k\in\mathbb Z\}$ is a
tiling with a prototile $S_{00}$, which is also considered as the
tiling from partitioning the plane~$\mathbb R^2$ into squares by
lines $L'_k$ and~$L''_j$ ($j$, $k\in\mathbb Z$), where $L'_t$ is
line $y=t$ and $L''_t$ is line $x=t$ ($t\in\mathbb R$).

Let $0\leqslant\delta<\frac12$. Assume that $C'_k$ is a Jordan
curve between $L'_{k-\delta}$ and~$L'_{k+\delta}$ and $C''_j$ is a
Jordan curve between $L''_{j-\delta}$ and~$L''_{j+\delta}$ such
that $C'_k\cap C''_j$ is a singleton for each pair $j$,
$k\in\mathbb Z$. Let $\mathcal T$ be a tiling from
partitioning~$\mathbb R^2$ into pieces by curves~$C'_k$ and~$C'_j$
($j$, $k\in\mathbb Z$). Then $\mathcal T$ is a monohedral tiling
with a $\delta$-quasi-prototile~$S_{00}$. It is easy to see that
translations $\tau(j,k):\:z\mapsto z+j+k\,\mi$ ($j$, $k\in\mathbb
Z$) and rotations $\sigma_m:\:z\mapsto\me^{^{\scriptstyle
\frac{m\pi\mi}2}}z$ ($m=0$, $1$, $2$, $3$) ($z\in\mathbb C=\mathbb
R^2$) are $\lambda$-quasi-symmetries of~$\mathcal T$ and the
transformation group~$\mathcal G$ generated by
$$
\{\sigma_m,\,\tau(j,k):\:m=0,\,1,\,2,\,3;\; j,\,k\in\mathbb Z\}
$$
is a $\lambda$-quasi-symmetry group of~$\mathcal T$, where
$\lambda=\frac{4\delta}{1-2\delta}$. We also see that $\mathcal T$
is a $\{\frac{\delta}{2^{^{j-1}}}\}$-quasi-self-similar tiling
with quasi-prototype $S_{00}$ of ratio~$2$ and any closed plane
domain~$P$ containing square $S'=\{(x,y):\:\delta\leqslant
x\leqslant1-\delta,\;\delta\leqslant y\leqslant1-\delta\}$ and
contained in square $S''=\{(x,y):\:-\delta\leqslant
x\leqslant1+\delta,\;-\delta\leqslant y\leqslant1+\delta\}$ is a
$4$-rep $\{\lambda_j\}$-quasi-tile
($\lambda_j=\frac{(2+2^{^{1-j}})\delta}{1-2\delta}$) or $4$-rep
$\lambda$-quasi-tile ($\lambda=\frac{4\delta}{1-2\delta}$) of
ratio $2$.

The above example can be easily generalized to cases of space and
$k$-hedral quasi-tilings.

\paragraph{{\bf 5.1.5. Quasi-isohedral tilings and quasi-anisohedral tiles}}
(cf. \cite{Schattschneider-Senechal-2004}). Let $\delta$,
$\delta_0$, $\delta_1$, $\delta_2\geqslant0$.

(1) {\it $\delta$-quasi-transitive action}. We say that a
transformation group~$\mathcal G$ {\it acts
$\delta$-quasi-transitively }\/ on a family $\mathcal
A=\{A_i:\:i\in I\}$ of subsets of~$\mathbb R^n$ if given any $i$,
$j\in I$ there exists $g_{ij}\in\mathcal G$ such that
$$
h(A_i,g_{ij}(A_j))\leqslant\delta,
$$
and on the other hand, for any $A\in\mathcal A$ and any
$g\in\mathcal G$ there exists $B\in\mathcal A$ such that
$$
h(B,g(A))\leqslant\delta.
$$
If $\mathcal G$ acts $\delta$-quasi-transitively on~$\mathcal A$,
then $\mathcal A$ is called a {\it $\delta$-quasi-orbit}\/
of~$\mathcal G$.

(2) A {\it $(\delta_1,\delta_2)$-quasi-isohedral tiling}\/ is a
tiling whose $\delta_1$-quasi-\linebreak[4]symmetry group acts
$\delta_2$-quasi-transitively on its tiles. Hence an isohedral
tiling is a $(0,0)$-quasi-isohedral tiling.

(3) A {\it $(\delta_0,\delta_1,\delta_2)$-quasi-anisohedral
tile}\/ ({\it type}\/) is a prototile (type) that admits at least
one monohedral $\delta_0$-quasi-tiling but no
$(\delta_1,\delta_2)$-quasi-isohedral tilings.

\paragraph{{\bf 5.1.6.} {\it Remark}.}
(1) In~5.1.1, 5.1.2 and~5.1.3 we may replace $\mathrm{r}(\cdot)$
by~$|\cdot|$ to get similar definitions.

(2) We may similarly define a quasi-tiling of $W\subseteq\mathbb
R^n$ and a quasi-symmetry~$\varphi$ (quasi-symmetry group) of a
tiling of~$W$.

(3) Since a quasi-tiling is also a tiling, some related results
for a tiling, such as the normality lemma (see
\cite[3.2.2]{Grunbaum-Shephard-1987} and
\cite[p.\,\nolinebreak55]{Schattschneider-Senechal-2004}), still
hold for a quasi-tiling (under suitable conditions).

\paragraph{{\bf 5.1.7. Some related questions.}}
We may consider {\it extensions of some relative classical results
to the case of quasi-tilings}. Below we still pose a few of
questions, some of which are related to Hilbert's eighteenth
problem (see \cite[Section~4.1]{Brass-Moser-Pach-2005},
\cite{Grunbaum-Shephard-1980}, \cite{Hilbert-1902} and
\cite[Sections~1.5 and~1.7]{Senechal-1995}, etc).

We remark here that there has existed a relative question
mentioned in \cite[p.\,\nolinebreak497]{Grunbaum-Shephard-1987}
(see also \cite[Problem~4 in Section~4.1]{Brass-Moser-Pach-2005}).

 Let $\delta$, $\delta_0$, $\delta_1$, $\delta_2\geqslant0$.

 (1) If there exists a $(\delta_0,\delta_1,\delta_2)$-quasi-anisohedral
tile, then what is the relation of $\delta_0$, $\delta_1$ and
$\delta_2$?

(2) According to \cite[p.\,\nolinebreak955 and
p.\,\nolinebreak956]{Grunbaum-Shephard-1980} it is more hopeless
to determine all $\delta$-quasi-prototiles of monohedral tilings
in~$\mathbb R^n$ for $\delta>0$. However, because of this reason,
we may find more $\delta$-quasi-prototiles of monohedral tilings
if $\delta$ is greater (assume we do not require a whole list is
shown). Perhaps we can find a large majority of
$\delta$-quasi-prototiles of monohedral tilings with the aid of
computers.

As $\delta\geqslant0$ gets smaller, the class of
$\delta$-quasi-prototiles of monohedral tilings  gets smaller.
Obviously prototiles of monohedral tilings are
$\delta$-quasi-prototiles of monohedral tilings and
$0$-quasi-prototiles of monohedral tilings are prototiles of
monohedral tilings. We make the suggestion above for considering
problems because prototiles of monohedral tilings are strict and
exact objects but $\delta$-quasi-prototiles of monohedral tilings
are freer relatively.

(3) (i) For a given set~$T$, such as a tetrahedron, a pentagon,
etc., determine $\delta_0\geqslant0$ so that when
$\delta>\delta_0$, there exists~$T_\delta$ satisfying
$\wth(T_\delta,T)\leqslant\delta$ and $T_\delta$ is a prototile of
a monohedral tiling but when $0\leqslant\delta<\delta_0$ there do
not exist any~$T_\delta$ satisfying
$\wth(T_\delta,T)\leqslant\delta$
 such that $T_\delta$ admit monohedral tilings.

 (ii) For a given set~$T$, such as a tetrahedron, a pentagon,
etc., determine $\delta_0\geqslant0$ so that $T$ admits
$\delta$-quasi-tilings as $\delta>\delta_0$ but $T$ admit no
$\delta$-quasi-tilings as $0\leqslant\delta<\delta_0$.

(iii) Do there exist any bodies~$T$ in~$\mathbb R^n$ such that for
some $\delta>0$ any~$T_\delta$ satisfying
$\wth(T_\delta,T)\leqslant\delta$ admit no monohedral tilings
whereas $T$ admit monohedral $\delta$-quasi-tilings?

(iv) Do there exist any bodies~$T$ in~$\mathbb R^n$ such that $T$
admit monohedral $\delta$-quasi-tilings  for any $\delta>0$
whereas $T$ admit no monohedral tilings? If this kind of sets~$T$
do exist then $T$ should be strange sets and probably possess
fractal boundaries. Among these sets more strange ones, if exist,
are positive answers of the following question.

(v) Do there exist any bodies~$T$ in~$\mathbb R^n$ such that for
any $\delta>0$ any~$T_\delta$ satisfying
$\wth(T_\delta,T)\leqslant\delta$ admit no monohedral tilings
whereas $T$ admit monohedral $\delta$-quasi-tilings?

\paragraph{{\bf 5.2. Quasi-patterns.}}

For the notion of a (mono-motif) pattern (in plane) we refer to
\cite[Chapter~5]{Grunbaum-Shephard-1987}. We now extend this
notion and as an example we also give a result, which is an
extension of a basic proposition in patterns.

At first, for two families $\{M_i:\:i\in I\}$ and $\{N_j:\:j\in
J\}$ (denoted $\{M_i\}$ and $\{N_j\}$) ($I$ and $J$ are two index
sets) of sets in~$\mathbb R^n$, we define
\begin{align*}
h(\{M_i\},\{N_j\}):=\sup_{i\in I,\,j\in
J}\left\{\right.\inf\{&h(M_i,N_{j'}):\:j'\in J\}, \\
&\inf\{h(M_{i'},N_j):\:i'\in I\}\left.\right\}.
\end{align*}

\paragraph{{\bf 5.2.1. Quasi-symmetry groups.}}
Let $\delta\geqslant0$. Let $\mathfrak P=\{M_i:\:i\in I\}$ be a
family of sets in~$\mathbb R^n$. A $\delta$-({\it quasi-}){\it
symmetry}\/ $\varphi$ of~$\mathfrak P$ is an isometry such that
$$
h(\{\varphi(M_i)\},\{M_i\})\leqslant\delta.
$$
 A $\delta$-({\it quasi-}){\it symmetry
group}\/ of $\mathfrak P$ is a group $\mathcal G$ consisting of
$\delta$-quasi-isometries of $\mathfrak P$.

\paragraph{{\bf 5.2.2. Quasi-patterns.}}
Let $\delta\geqslant0$, $\delta_{ij}\geqslant0$ and
$\delta_i\geqslant0$ ($i$, $j\in I$). Let $\mathfrak
P=\{M_i:\:i\in I\}$ be a nonempty family of nonempty subsets
of~$\mathbb R^n$ and $\mathcal G$ is a $\delta$-symmetry group.
Assume
\begin{enumerate}
\item[(i)] $M_i$ ($i\in I$) are pairwise disjoint;
\item[(ii)] Given any pairs $i$, $j\in I$ there exist
$g_{ij}\in\mathcal G$ so that
$$
h(M_i,g_{ij}(M_j)))\leqslant\delta_{ij}
$$
and
$$
\mathcal G'=\{g_{ij}:\:i,\,j\in I\}
$$
is a subgroup of~$\mathcal G$.
\end{enumerate}
Then $\mathfrak P$ is called a ({\it monomotif}\/)
$\langle\delta,\{\delta_{ij}\}\rangle$-({\it quasi}-){\it
pattern}\/ and each $M_i$ ($i\in I$) is called a {\it motif
quasi-copy}\/ or {\it motif-$i$ copy}\/ (or ({\it motif}\/) {\it
copy-$i$}) of $\mathfrak P$. A nonempty set $M\subseteq\mathbb
R^n$ satisfying
$$
\wth(M_i,M)\leqslant\delta_i\quad(i\in I)
$$
is called a $\delta_i$-({\it quasi}-){\it motif}\/ of~$\mathfrak
P$. Note that $\mathcal G'$ is not necessarily unique. Let
$\mathcal G_0$ is a maximal element of~$\{\mathcal G'\}$. Then we
call $\mathcal G_0$ a
$\langle\delta,\{\delta_{ij}\}\rangle$-(\nolinebreak{\it
quasi}-\nolinebreak){\it symmetry group}\/ of~$\mathfrak P$. If
$\delta_{ij}=\delta_0$ for all $i$, $j\in I$ then $\mathfrak P$ is
called a ({\it monomotif}\/)
$\langle\delta,\delta_0\rangle$-(\nolinebreak{\it
quasi}-\nolinebreak){\it pattern}\/ and a ({\it monomotif}\/)
$\delta$-(\nolinebreak{\it quasi}-\nolinebreak){\it pattern}\/ if
$\delta=\delta_0$ moreover.

\paragraph{{\bf 5.2.3. Discrete conditions.}}
We say that a quasi-pattern $\mathfrak P=\{M_i:\:i\in I\}$ is {\it
discrete}\/ if the following conditions hold:
\begin{enumerate}
\item[(D1)] All $M_i$ ($i\in I$) are bounded and usually they are
also connected (if all~$M_i$ are connected then $\mathfrak P$ is
called a {\it connected}\ quasi-pattern).
\item[(D2)] For each $i\in I$ there is an open set~$E_i$ which
contains $M_i$ but $M_j\cap E_i=\emptyset$ for all $j\in I$ and
$j\ne i$.
\item[(D3)] The cardinal number
$\sharp\nolinebreak\,I\geqslant2$.
\end{enumerate}

We say that $\mathfrak P$ is {\it $d_0$-discrete}\/ if (D2) is
replaced by the following stronger condition:
\begin{enumerate}
\item[(D2$'$)] $d_0=\inf\{d(M_i,M_j):\:i,\,j\in I;i\ne j\}>0$.
\end{enumerate}

\paragraph{{\bf 5.2.4. Engulfing and subtending.}}
Let $\delta_1$, $\delta_2\geqslant0$. Let $\mathfrak
P=\{M_i:\:i\in I \}$ and $\mathfrak Q=\{N_i:\:i\in I \}$ be two
quasi-patterns with the same index set~$I$.

(i) We say that $\mathfrak Q$ {\it
$\langle\delta_1,\delta_2\rangle$-engulfs}\/~$\mathfrak P$ if the
following two conditions hold:
\begin{enumerate}
\item[(E1)] $N_i\supseteq M_i$ for each $i\in I$;
\item[(E2)] There exist a $\delta_1$-symmetry group $\mathcal G_1$
of~$\mathfrak P$ and a $\delta_2$-symmetry group $\mathcal G_2$
of~$\mathfrak Q$ such that $\mathcal G_2\supseteq\mathcal G_1$.
\end{enumerate}

(ii) We say that $\mathfrak P$ {\it
$\langle\delta_1,\delta_2\rangle$-subtends}\/~$\mathfrak Q$ if the
following two conditions hold:
\begin{enumerate}
\item[(S1)] $M_i\subseteq N_i$ for each $i\in I$;
\item[(S2)] There exist a $\delta_1$-symmetry group $\mathcal G_1$
of~$\mathfrak P$ and a $\delta_2$-symmetry group $\mathcal G_2$
of~$\mathfrak Q$ such that $\mathcal G_1\supseteq\mathcal G_2$.
\end{enumerate}

For generality we may use the following condition to replace (E1)
and~(S1) above:
\begin{enumerate}
\item[(ES)] $\wth(M_i,N_i)\leqslant\delta^{(i)}$,
\end{enumerate}
where $\{\delta^{(i)}:\:i\in I\}$ is a set of fixed nonnegative
real numbers. Then we say that $\mathfrak Q$ {\it
$\langle\delta_1,\delta_2,\{\delta^{(i)}\}\rangle$-engulfs}\/~$\mathfrak
P$ and $\mathfrak P$ {\it
$\langle\delta_1,\delta_2,\{\delta^{(i)}\}\rangle$-subtends}\/~$\mathfrak
Q$ respectively.

\paragraph{{\bf 5.2.5. Definition.}}
(1) Let $A$ be a nonempty set of~$\mathbb R^n$. By the {\it
infield}\/ of~$A$, denoted $\mathrm{In}(A)$, we mean the
complement of~$A^{\mathrm{c}}_\infty$, where
$A^{\mathrm{c}}_\infty$ is the unbounded component of~$\mathbb
R^n\setminus A$.

(2) Let $\mathfrak P=\{M_i:\:i\in I\}$ be a quasi-pattern. If
$\mathrm{In}(M_i)\cap\mathrm{In}(M_j)=\emptyset$ for all $i$,
$j\in I$ and $i\ne j$, then we say $\mathfrak P$ is {\it
separated}. If each $\mathrm{In}(M_i)$ ($i\in I$) contains a ball
of radius~$r_0$ and contained in a ball of radius~$R_0$, where
$r_0$ and~$R_0$ are two positive constants only related
to~$\mathfrak P$, then we say that $\mathfrak P$ is {\it
fine-distributed}.

(3) Let $A$ be a nonempty set of~$\mathbb R^n$ and
$\lambda\geqslant0$. If $\lambda$-neighborhood $\mathscr
N(\mathrm{In}(A),\lambda)$ of~$\mathrm{In}(A)$ is still simply
connected, then $A$ is said to be {\it $\lambda$-exterior open
topology-free}.

\bigskip
Some corresponding classical results about patterns may be
considered to  be extended to the case of quasi-patterns. For
example, we have the following result about quasi-patterns similar
to \cite[5.1.1]{Grunbaum-Shephard-1987}.

\paragraph{{\bf 5.2.6. Proposition.}}
{\it
 Let\/ $\delta\geqslant0$ and\/ $d_0>4\delta$. Suppose\/
$\mathfrak P=\{M_i:\:i\in I\}$ is a\/ $d_0$-discrete connected
$\delta$-pattern in\/ $\mathbb R^2$ with a $\delta$-symmetry
group\/ $\mathcal G_1$. Then $\mathfrak P$ can be\/
$\langle\delta,4\delta+2\varepsilon_0\rangle$-engulfed by a
separated and fine-distributed connected\/
$\langle4\delta+2\varepsilon_0,0\rangle$-pattern\/ $\mathfrak Q$
$(\varepsilon_0$ is any fixed number satisfying\/
$0<\varepsilon_0<\frac{d_0-4\delta}2)$. If there is a motif
quasi-copy\/ $M_{i_0}\in\mathfrak P$ which is\/ $\lambda$-exterior
open topology-free, where\/
$\lambda=5\delta+3\varepsilon_0+\varepsilon$ and\/ $\varepsilon$
is any sufficiently small positive number, then\/ $\mathfrak P$
can be\/ $\langle\delta,4\delta+2\varepsilon_0\rangle$-engulfed by
a fine-distributed open disk\/
$\langle4\delta+2\varepsilon_0,0\rangle$-pattern\/ $\mathfrak Q$}
 ({\it all motif quasi-copies of\/ $\mathfrak Q$ are topological
disks}).
\par\noindent{\it Proof.}   Take $r=\delta+\varepsilon_0$, where
$0<\varepsilon_0<\frac{d_0-4\delta}2$. Suppose
$$
N_{i_0}:=\mathcal N_r(M_{i_0}).
$$
 Let $g_{ji_0}\in\mathcal G_1$ satisfy
\begin{equation}
h(g_{ji_0}(M_{i_0}),M_j)\leqslant\delta \label{tag5.2-1}
\end{equation}
for $j\in I$ and $j\ne i_0$. Suppose
\begin{equation}
N_j:=g_{ji_0}(N_{i_0}). \label{tag5.2-2}
\end{equation}
Then each~$N_j$ ($j\in I$) is connected. Below we show that
$$
\mathfrak Q:=\{N_i:\:i\in I\}
$$
is a separated fine-distributed
$\langle4\delta+2\varepsilon_0,0\rangle$-pattern and $\mathfrak Q$
$\langle\delta,4\delta+2\varepsilon_0\rangle$-engulfs~$\mathfrak
P$.

First, we prove $N_i\cap N_j=\emptyset$ for $i$, $j\in I$ and
$i\ne j$. Let $0<\varepsilon<\frac{d_0-4\delta-2\varepsilon_0}3$.
Select suitable $n_i\in N_i$, $n_j\in N_j$, $n'_{i_0}$,
$n''_{i_0}\in N_{i_0}$, $m'_{i_0}$, $m''_{i_0}\in M_{i_0}$,
$m_i\in M_i$ and $m_j\in M_j$ so that
\begin{align*}
d(N_i,N_j)&>d(n_i,n_j)-\varepsilon \\
 &\geqslant d(m_i,m_j)-d(n_i,m_i)-d(n_j,m_j)-\varepsilon \\
 &\geqslant d(M_i,M_j)-[d(g_{ii_0}(n'_{i_0}),g_{ii_0}(m'_{i_0}))+d(g_{ii_0}(m'_{i_0}),m_i)] \\
 &\qquad-[d(g_{ji_0}(n''_{i_0}),g_{ji_0}(m''_{i_0}))+d(g_{ji_0}(m''_{i_0}),m_j)]-\varepsilon \\
 &\geqslant d_0-2[r+(\delta+\varepsilon)]-\varepsilon
 =d_0-4\delta-2\varepsilon_0-3\varepsilon>0.
\end{align*}
For $i\in I$, $g\in\mathcal G_1$ and $\varepsilon>0$, since
$gg_{ii_0}\in\mathcal G_1$, we can take a suitable $j\in I$ such
that
$$
h(gg_{ii_0}(M_{i_0}),M_j)<\delta+\varepsilon.
$$
Thus
\begin{align*}
 h(g(N_i),N_j)&\leqslant h(gg_{ii_0}(N_{i_0}),gg_{ii_0}(M_{i_0}))+h(gg_{ii_0}(M_{i_0}),M_j) \\
 &\qquad+h(M_j,g_{ji_0}(M_{i_0}))+h(g_{ji_0}(M_{i_0}),g_{ji_0}(N_{i_0})) \\
 &<r+(\delta+\varepsilon)+\delta+r=4\delta+2\varepsilon_0+\varepsilon.
\end{align*}
By the same reasoning we can also get
$$
h(N_i,g(N_j))<4\delta+2\varepsilon_0+\varepsilon
$$
for a suitable $j\in I$. Therefore we have
\begin{equation}
 h(\{g(N_i)\},\{N_i\})\leqslant 4\delta+2\varepsilon_0. \label{tag5.2-3}
\end{equation}

For $i$, $j\in I$, let $g=g_{ji_0}g_{ii_0}^{-1}$. Then
$g\in\mathcal G_1$ and
$$
g(N_i)=g_{ji_0}g_{ii_0}^{-1}g_{ii_0}(N_{i_0})=g_{ji_0}(N_{i_0})=N_j.
$$
Consequently
$$
h(g(N_i),N_j)=0.
$$

It is obvious that $\mathfrak Q$ is fine-distributed. Since $N_i$
($i\in I$) are disjoint from and isometric to each other, it
follows that the bounded components $C_{k_i}^{(i)}$ ($k_i=1$, $2$,
$\dots$) of the complement of~$N_i$ are disjoint from each~$N_j$
and from the bounded components $C_{k_j}^{(j)}$ ($k_j=1$, $2$,
$\dots$) of the complement of~$N_i$ ($j\in I$, $j\ne i$). Hence
$\mathfrak Q$ is separated.

From (\ref{tag5.2-1}) and~ (\ref{tag5.2-2}) it follows that
$M_j\subseteq N_j$ for all $j\in I$. We may choose $\mathcal G_2$
to be a maximal element of~$\{\mathcal G'\}$ satisfying the
definition of $\langle4\delta+2\varepsilon_0,0\rangle$-patterns
according to~5.2.2 and including~$\mathcal G_1$.

Finally let $M_{i_0}$ be $\lambda$-exterior open topology free for
$\lambda=5\delta+3\varepsilon_0+\varepsilon$, where $\varepsilon$
is any sufficiently small positive number, and let
$N_i^*=\mathrm{In}(N_i)$ ($i\in I$). Then by~(\ref{tag5.2-3}) we
can obtain
$$
h(\{g(N_i^*)\},\{N_i^*\})\leqslant 4\delta+2\varepsilon_0.
$$
Therefore $\mathfrak Q^*=\{N_i^*:\:i\in I\}$ is an open disk
$\langle4\delta+2\varepsilon_0,0\rangle$-pattern. \hfill $\Box$

\paragraph{{\bf 5.2.7. Corollary.}}
{\it
 Let\/ $\delta\geqslant0$ and\/ $d_0>4\delta$. Suppose\/
$\mathfrak P$ is a\/ $d_0$-discrete connected $\delta$-pattern
in\/~$\mathbb R^2$. Then
 $\mathfrak P$ is locally finite $($i.e.,
 set $\{M:\:M\in\mathfrak P,\,M\cap D\ne\emptyset\}$ is finite for any
 disk~$D)$.
 }

\paragraph{{\bf 5.3. Quasi-packings.}}

\paragraph{{\bf 5.3.1. Packing density.}}
We recall that a family $\boldsymbol P=\{P_i:\:i\in I\}$ ($I$ is
an at most countable index set) of compact sets with nonempty
interiors is said to form a {\it packing}\/ in a domain
$\varOmega\subseteq\mathbb R^n$ if $\bigcup_{i\in
I}P_i\subseteq\varOmega$ and no two members of~$\boldsymbol P$
have an interior point in common. Suppose a bounded domain~$D$ and
members $P_i$ ($i\in I$) of~$\boldsymbol P$ are Jordan-measurable.
The density of the packing~$\boldsymbol P$ relative to~$D$ is
defined as
$$
\mathsf{d}(\boldsymbol P,D):=\frac{\mathrm{Vol}(P_i\cap
D)}{\mathrm{Vol}(D)}.
$$
The {\it upper}\/ and {\it lower densities}\/ of~$\boldsymbol P$
(in~$\mathbb R^n$) are
\begin{align*}
 \overline{\mathsf{d}}(\boldsymbol P)&:=\limsup_{r\to+\infty}\mathsf{d}(\boldsymbol P,B(r)) \\
\intertext{and}
 \underline{\mathsf{d}}(\boldsymbol P)&:=\liminf_{r\to+\infty}\mathsf{d}(\boldsymbol P,B(r))
\end{align*}
respectively, where $B(r):=\{x\in\mathbb R^n:\:d(x,o)<r\}$ ($o$
denotes the origin). If these two numbers are the same, then it is
called the {\it density}\/ of~$\boldsymbol P$ (in~$\mathbb R^n$),
denoted $\mathsf{d}(\boldsymbol P)$. Let $P$ be an $n$-dimensional
compact set with nonempty interior. Then the {\it packing
density}\/ of~$P$ is the largest density of a packing of congruent
copies of~$P$ in~$\mathbb R^n$.

For the notion of packings we refer to
\cite[1.1]{Brass-Moser-Pach-2005} and \cite{FejesToth-2004},
 etc.

\paragraph{{\bf 5.3.2.} {\it Example.}}
Let $\varepsilon\geqslant0$. Let $P$ be the dodecagon
$A_0A_1\cdots A_{11}$, where the vertexes are: $A_0(0,0)$,
$A_1(2,0)$, $A_2(2,4)$, $A_3(3,4)$, $A_4(3,0)$, $A_5(5,0)$,
$A_6(5,5)$, $A_7(3+\varepsilon,5)$, $A_8(3+\varepsilon,8)$,
$A_9(2-\varepsilon,8)$, $A_{10}(2-\varepsilon,5)$, $A_{11}(0,5)$.
We consider a packing of congruent copies of~$P$ in the plane in
two cases: (i) $\varepsilon>0$ and (ii) $\varepsilon=0$. We will
find the packings in the two cases have an obvious difference,
even if $\varepsilon$ is sufficiently small in case~(i). The
packing densities have an obvious jump even if the difference of
$\varepsilon$ and~$0$ is very very tiny.

\bigskip
 This example, along with
experience of life, tells us that sometimes if we ``squeeze''
bodies we may obtain much more dense packings. Thus we may give a
corresponding notion below.

\paragraph{{\bf 5.3.3. Quasi-packing density.}}
Let $P\subseteq\mathbb R^n$ be a compact set with nonempty
interior and $\delta\geqslant0$. If $Q\subseteq\mathbb R^n$
satisfies
$$
\wth(Q,P)\leqslant\delta
$$
then $Q$ is called a $\delta$-({\it quasi-}){\it copy}\/ of~$P$.

We define the $\delta$-({\it quasi-}){\it packing density}\/
$\mathsf{d}_\delta(P)$ of~$P$ as the largest density of a packing
of $\delta$-copies with nonempty interiors of~$P$ in~$\mathbb R^n$
(this packing is called a $\delta$-({\it quasi-}){\it packing}\/
of~$P$ (in~$\mathbb R^n$)), i.e.,
$$
\mathsf{d}_\delta(P):=\sup\{\mathsf{d}(\boldsymbol
P):\:\boldsymbol P \text{ is a $\delta$-packing of~P in~$\mathbb
R^n$}\}.
$$

If $\mathsf{d}_\delta(P)$ as a function of $\delta$ is not
continuous at $\delta=\delta_0$, then $P$ is said to possess the
{\it collapse property}\/ at $\delta=\delta_0$, $\delta_0$ is
called a {\it collapse value}\/ of~$P$ and
$$
\varDelta_{\delta_0}(P):=\lim_{\delta\to\delta_0^+}\mathsf{d}_\delta(P)
-\lim_{\delta\to\delta_0^-}\mathsf{d}_\delta(P)
$$
is called the {\it collapse quantity}\/ of~$P$ at
$\delta=\delta_0$.

In Example~5.3.2 it is easy to see that $P$ possesses the collapse
property at $\delta=\frac{\varepsilon}2$.

\paragraph{{\bf 5.3.4. Some related questions.}}
(1) Give a few more ``natural'' examples than Example~5.3.2 (we
may see some ``natural'' strange phenomena in~\cite{Zong-1996}).

(2) How much is $\sup\{\varDelta_{\delta_0}(P):\:P$  is a compact
connected set with nonempty interior in~$\mathbb R^n$ which
possesses the collapse property at $\delta=\delta_0\}$?

(3) Do there exist any compact connected sets~$P$ with nonempty
interiors in~$\mathbb R^n$ which have infinitely many collapse
values with $0$ as an accumulation point? If the answer is
positive, then may all these collapse values of some~$P$ have an
uncountable cardinal?

\paragraph{{\bf 5.4. Approximate crystals.}}
For nearly two hundred years the internal structure of a crystal
was considered periodic. The discovery of X-ray diffraction and
Laue's experiment (by W. Friedrich and P. Knipping) supported this
viewpoint. A (classical) crystal may be defined as the union of a
finite number of orbits of a crystallographic group. In this
subsection we refer to \cite{Senechal-1995}) and
\cite{Senechal-2004}.

After the announcement of the discovery of crystals with
icosahedral symmetry in the year 1984 the classical viewpoint has
been extended.
 In 1992 the Commission on Aperiodic
Crystals of the International Union of Crystallography proposed as
a working definition: a {\it crystal}\/ is a solid with an
essentially discrete diffraction pattern.
 In \cite{Senechal-2004} a ({\it generalized}\/) {\it crystal}\/
 is defined as a Delone set $\varLambda$ with
 nontrivial~$\varLambda_d$. Here a {\it Delone set}\/ is an
 {\it $(r,R)$ system}\/ ($r$, $R>0$), i.e., a set $\varLambda$ of
 points in $\mathbb R^n$ that is $r$-discrete and {\it relatively
 dense}\/ (every sphere of radius $R$ contains at least one point
 of $\varLambda$).

On the other hand, in the case of noncrystals, because of their
disorder, the method of radial distribution functions is applied
to deal with this situation.

Below let us try to consider the problem from another way.

\paragraph{{\bf 5.4.1. Approximate crystals.}}
Let $\delta\geqslant0$. A {\it $\delta$-approximate crystal}\/
({\it $\delta$-crystal}\/) is defined as a {\it
$\delta$-dot-pattern}\/ ($\delta$-pattern consisting of
singletons) with a crystallographic group as its $\delta$-symmetry
group (called a {\it $\delta$-symmetry crystallographic group}\/).

Let $\mathfrak P=\{P_i:\:i\in I\}$ be a $\delta$-crystal and let
$\mathfrak G$ denote the set of all its $\delta$-symmetry
crystallographic groups. Let
$$
\lambda(\mathcal G):=\frac{\inf\{h(\{P_i\},\{g(P):\:g\in\mathcal
G\}):\:P\in\mathbb R^n\}}{d_{\mathcal G}},
$$
where
$$
d_{\mathcal G}:=\inf\{d(g_1(P),g_2(P)):\:g_1,\,g_2\in\mathcal
G,\,P\in\mathbb R^n,\,g_1(P)\ne g_2(P)\},
$$
for $\mathcal G\in\mathfrak G$. It is easy to see $d_{\mathcal
G}\ne0$. Define
\begin{equation}
\lambda=\lambda_{\mathfrak G}:=\inf\{\lambda(\mathcal
G):\:\mathcal G\in\mathfrak G\}. \label{tag5.4-1}
\end{equation}
We call $\mathfrak P$ a {\it $\lambda\!\sim$approximate crystal}
or {\it $\lambda\!\sim$crystal}.

If there exists $\mathcal G\in\mathfrak G$ and $P\in\mathbb R^n$
such that for each $i\in I$ there exists $g\in\mathcal G$
satisfying
\begin{equation}
d(P_i,g(P))=\lambda_{\mathfrak G}\,d_{\mathcal G},
\label{tag5.4-2}
\end{equation}
then $\mathfrak P$ is called a {\it strict
$\lambda\!\sim$crystal}.

Now we give simple and obvious definitions of an approximate
crystal in two ways.

(1) {\it A quasi-tiling model}. Let $\mathcal T=\{T_i:\:i\in I\}$
be a monohedral tiling model of a crystal~$\boldsymbol C$ with a
prototile~$P_0$ and $\delta\geqslant0$. Assume $\mathcal
T_\delta=\{T'_i:\:i\in I\}$ is a quasi-tiling whose tiles are
still polyhedra obtained by perturbing~$\mathcal T$ such that
\begin{equation}
\wth(T'_i,P_0)\leqslant\delta\mathrm{r}(P_0),\label{tag5.4-3}
\end{equation}
where $\mathrm{r}(P_0)$ is the radius of $P_0$ (see
Definition~2.6.2). If a non-crystal $\boldsymbol C_\delta$ is
obtained by perturbing the crystal $\boldsymbol C$ and
$\boldsymbol C_\delta$ possesses tiling model $\boldsymbol
C_\delta$, then $\boldsymbol C_\delta$ is called a {\it
$\delta\!\sim\!\boldsymbol C$-crystal}.

(2) {\it A discrete point-set model}. Let $\mathfrak
C=\{C_i:\:i\in I\}$ be a dot-pattern model of a normal crystal and
$\lambda\geqslant0$. Let $\mathfrak P=\{P_i:\:i\in I\}$ be a
discrete set of points. Suppose
\begin{equation}
\wth(\mathfrak P,\mathfrak C)\leqslant\lambda d_{\mathfrak
C},\label{tag5.4-4}
\end{equation}
where
\begin{equation}
d_{\mathfrak C}=\min\{d(C_i,C_j):\:i,\,j\in I,\,i\ne j\}.
\label{tag5.4-5}
\end{equation}
Then $\mathfrak P$ is called a {\it $\lambda\!\sim\!\mathfrak
C$-crystal}. If moreover the equality in~(\ref{tag5.4-4}) always
holds, i.e.,
\begin{equation}
\wth(\mathfrak P,\mathfrak C)=\lambda\, d_{\mathfrak
C},\label{tag5.4-6}
\end{equation}
then $\mathfrak P$ is called a {\it strict
$\lambda\!\sim\!\mathfrak C$-crystal}.

\paragraph{{\bf 5.4.2. Some related questions.}}
Given a substance $\boldsymbol S$ ignoring its internal structure,
let $\Sigma$ be the set of all approximate crystals consisting of
the same substance $\boldsymbol S$. Suppose $\mathfrak C$ is the
crystal consisting of $\boldsymbol S$. Since a $0\!\sim\!\mathfrak
C$-crystal is a normal crystal, it possesses a discrete
diffraction pattern. We may imagine that when $\lambda>0$ is
sufficiently small a $\lambda\!\sim\!\mathfrak C$-crystal may
still produces a diffraction pattern. Now we ask the following
questions.

(1) How much is $\lambda_{\boldsymbol
S}:=\sup\{\lambda:\:\mathfrak P\in\Sigma$ is a (well-distributed
random strict) $\delta\!\sim\!\mathfrak C$-crystal, where
$0\leqslant\delta\leqslant\lambda$, which produces an essentially
discrete diffraction pattern$\}$?

(2) For different kinds of substances $\boldsymbol S$, are
$\lambda_{\boldsymbol S}$ the same or different?

(3) Furthermore, find the laws of changes of the properties (in
mechanics, acoustics, heat, optics, electricity and magnetism,
etc.) of a (well-distributed random strict)
$\lambda\!\sim\!\mathfrak C$-crystal depending on $\lambda$.

\paragraph{{\bf 5.4.3.} {\it Remark}.}
(1) In the research if we weaken the demand of ``strict'' in the
strict $\lambda\!\sim$crystal or strict $\lambda\!\sim\!\mathfrak
C$-crystal, e.g. we require ``strict'' in nanometer or over scales
but ignore ``strict'' in angstrom scales (for nanometer solid
materials) the problem perhaps becomes easier.

(2) Since the scanning tunneling microscope (STM) was invented,
Richard P. Feynman's imagination that we could arrange the atoms
one by one the way we want them has seen the dawn of its
realization. The questions that we ask above are just some more
concrete ones related to Feynman's question: {\it what would the
properties of materials be if we could really arrange the atoms
the way we want them}? (see \cite{Feynman-1960})

(3) Experts in different research fields may put forward different
related questions. For example, one may pose similar questions for
substances in liquid or gaseous states or for organic compounds,
etc.

(4) In fact we could ask a few of questions in this subsection,
some of which might seem abrupt or even unreasonable. For example,
we may ask the question: if we had calculated out the critical
values (interval) of the shape vision error (see Section~4) and
$\lambda_{\boldsymbol S}$, what might we say about them?

(5) As examples, we define several concepts for further
consideration. Let $\mu>0$, $\lambda\geqslant0$ and
$0\leqslant\lambda_1\leqslant\lambda_2$.

(i) Let $\mathfrak C=\{C_i:\:i\in I\}$ be a dot-pattern model of a
normal crystal. Let $\mathfrak P=\{P_i:\:i\in I\}$ be a discrete
set of points. Suppose (\ref{tag5.4-4}) or (\ref{tag5.4-6}) is now
changed into
\begin{equation}
\lambda_1 d_{\mathfrak C}\leqslant\wth(\mathfrak P,\mathfrak
C)\leqslant\lambda_2 d_{\mathfrak C},\label{tag5.4-7}
\end{equation}
where $d_{\mathfrak C}$ has been defined in (\ref{tag5.4-5}). Then
$\mathfrak P$ is called a {\it
$\langle\lambda_1,\lambda_2\rangle\!\sim\!\mathfrak C$-crystal}.

(ii) Let $\varLambda$ be a {\it uniformly discrete}\/ set of
points in~$\mathbb R^n$, i.e.,
$$
d_{\varLambda}:=\inf\{d(P_1,P_2):\:P_1,\,P_2\in\varLambda\}>0.
$$
 If
$$
\wth(\mu\varLambda,\varLambda)\leqslant\lambda\,d_{\varLambda},
$$
then $\varLambda$ is said to possess
$\langle\mu;\lambda\rangle$(-{\it approximate}\/) {\it inflation
symmetry}\/ (cf. \cite{Senechal-2004}).

(iii) If in~(i) we change (\ref{tag5.4-7}) into
$$
\lambda_1 d_{\mathfrak C}\leqslant\wth(\mu\mathfrak P,\mathfrak
C)\leqslant\lambda_2 d_{\mathfrak C},
$$
then $\mathfrak P$ is called a {\it
$\langle\mu;\lambda_1,\lambda_2\rangle\!\sim\!\mathfrak
C$-crystal}. Similarly we may change (\ref{tag5.4-2}) into
$$
\left|\frac{d(P_i,g(P))}{d_{\mathcal G}}-\lambda_{\mathfrak G}
\right|\leqslant\delta,
$$
where $\delta\geqslant0$ and $\lambda_{\mathfrak G}$ is the one
defined in (\ref{tag5.4-1}), and change (\ref{tag5.4-3}) into
$$
\delta_1\mathrm{r}(P_0)\leqslant\wth(T'_i,P_0)\leqslant\delta_2\mathrm{r}(P_0)
$$
or more generally
$$
\delta_1\mathrm{r}(P_0)\leqslant\wth(\mu
T'_i,P_0)\leqslant\delta_2\mathrm{r}(P_0),
$$
 where $0\leqslant\delta_1\leqslant\delta_2$, to obtain
 corresponding definitions.
In the same manner we may give more general notions than the ones
in quasi-tilings and quasi-patterns.

(6) In this paper if metric $\wth$ is replaced by $\overline{h}$
and $\widehat{h}$ is replaced by $\check{h}$, corresponding
concepts can be defined and corresponding results and questions
can be considered similarly.

\bigskip

 Department of Mathematics, Shanghai Jiao Tong University

  800 Dongchuan Road, Shanghai 200240, China

  \textit{e-mail:}\, jyyu@sjtu.edu.cn

\end{document}